\documentclass[a4paper,11pt]{article}

\usepackage{graphicx}
\usepackage{amsmath}
\def\gleq{{\ \mathop{=}\limits_{>}^{<} \ }}

\begin{document}

\title{\bf Increasing power of the test through pre-test
\\- a robust method}

\author{Rossita M Yunus\footnote{on leave from Institute of Mathematical Sciences, Faculty of Sciences, University of Malaya, Malaysia.}\ \  and Shahjahan Khan \\
Department of Mathematics and Computing \\ Australian Centre for Sustainable Catchments  \\
University of Southern Queensland \\ Toowoomba, Q 4350,
AUSTRALIA.\\
{\em Emails:yunus@usq.edu.au and khans@usq.edu.au} }
\date{}
\maketitle

\begin{abstract} This paper develops robust test procedures for
testing the intercept of a simple regression model when it is
\textit{apriori} suspected that the slope has a specified value.
Defining unrestricted test (UT), restricted test (RT) and pre-test
test (PTT) corresponding to the unrestricted (UE), restricted
(RE), and preliminary test estimators (PTE) in the estimation
case, the M-estimation methodology  is used to formulate the
M-tests and derive their asymptotic power functions. Analytical
and graphical comparisons of the three tests are obtained by
studying the power functions with respect to size and power of the
tests. It is shown that PTT achieves a reasonable dominance over
the others asymptotically.
\end{abstract}



\section{Introduction}
\setcounter{equation}{0}
In recent years many researchers have contributed to the
estimation of one parameter in the presence of uncertain prior
information on the value of another parameter. In general,
inclusion of non-sample prior information improves the quality of
inference. In spite of plethora of work in the area of improved
estimation using non-sample prior information (c.f. Saleh, 2006),
very little attention has been paid on the testing of parameters
in the presence of uncertain prior information. It may be a normal
expectation that testing of one parameter after pre-testing on
another would improve the performance of the ultimate test in the
sense of better power and size of the ultimate test. In this
paper, this improvement is achieved by using a robust test, namely
score type M-test defined along the line of M-estimation
methodology.

If the underlying distribution is known and the assumed model
holds, the statistical test that offers the most powerful test is
the classical likelihood ratio test (LRT). However, this
parametric test is generally very non-robust, even to small
departures from the assumed distribution (c.f. Huber, 1981, p.264,
Jur$\check{e}$ckov$\acute{a}$ and Sen, 1996, p.408). Several
robust tests are suggested in literature to tackle the problem.
For example, Huber (1981, p.264) suggests some censored likelihood
ratio type test and shows the test possesses a minimax property.
The test however does not work out conveniently for composite null
hypothesis against composite alternative testing problem (c.f.
Jur$\check{e}$ckov$\acute{a}$ and Sen 1996 p.407). The Rao's
efficient score statistic could be tempted on the least favorable
distribution to obtain a robust test statistic. Sen (1982)
suggests a score type statistic by replacing the robust Rao score
test statistic by an M-statistic. The score-type M-test has
several advantages. The construction of the score-type M-test
needs less parameters to estimate than the robust LRT suggested by
Schrader and Hettmansperger (1980) yet both are equiefficient
(c.f. Sen, 1982) and are applicable to composite hypothesis
testing.

The properties of unrestricted estimator (UE), restricted
estimator (RE) and pre-test estimator (PTE) have been investigated
by many authors (Khan and Saleh, 2001, Khan et al., 2002). Most of
the studies are based on normal or t-models and the results are
non-robust. In the studies, the PTE (a linear combination of UE
and RE) possesses a small quadratic risk when the distance
parameter are large and too close to zero, that makes it the best
choice over the other two estimators. Instead of least squares
(LS) and maximum likelihood (ML) estimators, the properties of UE,
RE and PTE are also studied in the framework of general robust
estimators, explicitly, M-estimators. As such, a robust estimator
namely the preliminary test M-estimator (PTME) are proposed for
linear models (Sen and Saleh, 1987). In this paper, three tests
correspond to the UE, RE and PTE are defined. They are
unrestricted test (UT), restricted test (RT) and pre-test test
(PTT). The study of the properties on these tests formulated using
robust statistics is unavailable in the statistics literature. 

The properties of size of the pre-test as well as the power of the
test followed by pre-test have been studied in parametric cases
(Bechhofer, 1951, Bozivich et al., 1956). After almost three
decades, the effect of pre-test (on slope) on the size and power
of the ultimate test (on the intercept) are investigated for
rank-based nonparametric cases by Saleh and Sen (1982). However,
there are some limited discussions in investigating the power of
the PTT discussed in the paper. To author's knowledge, no research
has been done in the investigation of the performance (size and
power) of the ultimate test following a pre-test in linear models
that is formulated using the score-type M-test. Since the
M-estimation method is more popular compared to the other robust
methods, it is incomplete to ignore the study of the performance
characteristics of the power function after pre-testing based on
M-test. M-estimation is known for its flexibility and well defined
for a variety of models for which MLE is also defined (Huber,
1981, p.43, Jur$\check{e}$ckov$\acute{a}$ and Sen, 1996 p.80). In
this paper, the study on the power of test after pre-testing using
M-test is considered for a simple linear regression model.

Consider a simple regression model of $n$ observable random
variables, $X_i,\; i=1,\ldots,n$
\begin{equation}X_i
= \theta + \beta c_i + e_i, \label{regress}
\end{equation}
where the errors $e_i$'s are from an unspecified symmetric and
continuous distribution function, $F_i, \; i=1,\ldots,n,$ the
$c_i$'s are known real constants of the explanatory variable and
$\theta$ and $\beta$ are the unknown intercept and slope
parameters respectively.

We wish to test the significance of the intercept parameter under
various conditions on the slope parameter. Basically, testing the
intercept of a simple linear regression depends on the knowledge
on the slope. We may have the case when the slope is unspecified.
For this case, the slope is treated as a nuisance parameter in
testing the significance of the intercept and we refer the test as
the unrestricted test (UT). We may also have the case when the
value of the slope is specified (say zero) and thus we return to
the situation of testing the location parameter. The test on the
intercept after specifying the value of slope is defined as the
restricted test (RT).  Besides these two cases, if the value of
the slope is suspected to be close to 0 (or any specified value),
a natural action is to remove the uncertainty of the suspicious
value of the slope by performing a test on the slope before
testing the intercept. If the null hypothesis of the pre-test is
rejected, then we use the UT, otherwise the RT. This is in line
with definition of the PTE in the estimation problem. For the
final case, the ultimate test following a pre-test is defined as
the pre-test test (PTT). Obviously the preliminary test (on the
slope) affects the power and size of the ultimate test (on the
intercept). To simplify,
\begin{itemize}
\item{\em{The unrestricted test}} : Test function $\phi_n^{UT}$ is
designed for testing $H_0^\star:\theta=0$ when $\beta$ is
unspecified, \item{\em{The restricted test}}: The test function
$\phi_n^{RT}$ is designed for testing $H_0^\star:\theta=0$ when
$\beta$ is 0 (specified) and \item{\em{The pre-test test}} : The
test function $\phi_n^{PTT}$ is designed for testing
$H_0^\star:\theta=0$ following a pre-test on the slope.
\end{itemize}

The objectives of the current research are to propose a robust
test statistic based on M-statistic to formulate the asymptotic
power functions for testing the intercept after pre-testing on
slope and to carry out investigations on the asymptotic properties
of this power function.

Along with some preliminary notions, the method of M-estimation is
presented in Section 2. In Section 3, three statistical tests
concerning testing on the intercept, namely, the UT, RT and PTT
are proposed for the three different cases mention earlier.
Further, the asymptotic distributions of the test statistics are
derived in Section 4. In Section 5, the asymptotic power functions
of the tests are derived. Section 6 is devoted to the analytical
results comparing the asymptotic power functions of the UT, RT and
PTT while the investigation of the power functions through an
illustrative example is presented in Section 7. The final section
presents discussions and concluding remarks.

\section{The M-estimation}
\setcounter{equation}{0} Given an absolutely continuous function
$\rho: \Re \rightarrow \Re$, M-estimator of $\theta$ and $\beta$
is defined as the values of $\theta$ and $\beta$ that minimize the
objective function
\begin{equation}
\sum_{i=1}^n \rho(X_i -\theta - \beta c_i). \label{eq:rho}
\end{equation} M-estimator of $\theta$ and $\beta$ can also be defined as the solutions of the
system of equations,
\begin{equation}
\begin{array} {ccccc}
\sum_{i=1}^n \psi_\theta(X_i) &=&\sum_{i=1}^n \psi(X_i -\theta - \beta c_i)&=&0, \\
\sum_{i=1}^n \psi_\beta(X_i) &=&\sum_{i=1}^n c_i \psi(X_i-\theta -
\beta c_i)&=&0. \label{eq:psi}
\end{array}
\end{equation}
If $\rho$ is differentiable with partial derivatives $\psi_\theta
=\partial{\rho}/\partial{\theta}$ and $\psi_\beta
=\partial{\rho}/\partial{\beta}$, then the M-estimators that
minimize the function in (\ref{eq:rho}) are the solutions to the
system (\ref{eq:psi}). On the contrary, the M-estimators obtained
from solving system (\ref{eq:psi}) may not minimize equation
(\ref{eq:rho}) (c.f. Caroll and Rupert, 1988 p.210). The system of
equations (\ref{eq:psi}) may have more roots, while only one of
them leads to a global minimum of (\ref{eq:rho}).
Jur$\check{e}$ckov$\acute{a}$ and Sen (1996) have given proof that
there exists at least one root of (\ref{eq:psi}) which is a
$\sqrt{n}$ - consistent estimator of $\theta$ and $\beta$ under
some conditions [c.f. p.215 - 224]. The $\psi$ function is
decomposed into the sum
\[\psi=\psi_a +\psi_c+\psi_s,\]
where \begin{description} \item{(a)} $\psi_a$ is absolutely
continuous function with absolutely continuous derivative.
\item{(b)} $\psi_c$ is a continuous, piecewise linear function
with knots at $\mu_1,\ldots, \mu_k$, that is, constant in a
neighborhood of $\pm \infty$ and hence its derivative is a step
function $\psi_c^\prime(z) = \alpha_v,\; \mu_v<z<\mu_{v+1},\;
v=0,1,\ldots,k$ where $\alpha_0,\ldots, \alpha_k \;\epsilon \;
\Re,\; \alpha_0=\alpha_k=0$ and $-\infty=\mu_0 <\mu_1
<\ldots<\mu_{k+1}=\infty$. We assume that $f(z)=\frac{dF(z)}{dz}$
is bounded in neighborhoods of $S_{\mu_1},\ldots,S_{\mu_k}.$
\item{(c)} $\psi_s$ is a nondecreasing step function,
$\psi_s(z)=\lambda_v,\;q_v <z \leq q_{v+1},v=1,\ldots,m$ where
$-\infty=q_0<q_1<\ldots<q_{m+1}=\infty$ and $-\infty < \lambda_0
<\lambda_1 <\ldots <\lambda_m<\infty.$ We assume that
$0<f(z)=(d/dz)F(z)$ and $f^\prime(z)=(d^2/dz^2)F(z)$ are bounded
in neighborhoods of $S_{q_1},\ldots,S_{q_m}.$
\end{description} The asymptotic result under conditions M1 to M5
of Jur$\check{e}$ckov$\acute{a}$ and Sen (1996, p.217) is used in
this paper. Further assume that all $\psi_a$, $\psi_c$ and
$\psi_s$ are nondecreasing and skew symmetric that is $\psi_j(-x)=
-\psi_j(x),\; j=1,2,3.$ Let $F$ be symmetric about 0, so that
\[ \int_{-\infty}^\infty \psi(x)dF(x)=0.\]
Assume that \begin{equation} \sigma_0^2 = \int_{-\infty}^\infty
\psi^2(x)dF(x). \end{equation} Following
Jur$\check{e}$ckov$\acute{a}$ and Sen (1996, p.217), two cases are
considered: \begin{description} \item(i) if $\psi_s=0$ then
\begin{equation} \gamma =
\int_{-\infty}^\infty (\psi_a^\prime(x)+ \psi_c^\prime(x)) f(x)dx.
\end{equation}  \item(ii) if $\psi_a=\psi_c=0$, then
\begin{equation} \gamma = \sum(\lambda_v - \lambda_{v-1})
f(S_{q_v}). \end{equation} \end{description} Further assume that
$\sigma_0$ and $\gamma$ are both positive and finite quantities.
Let the distribution function, $F$ be continuous and symmetric
about zero and have finite Fisher information,
\begin{equation}
I(f)=\int_{-\infty}^{\infty} \{ f'(x)/f(x) \}^2 dF(x),
\end{equation} where $f'(x) = (d/dx)f(x)=(d^2/dx^2)F(x).$
Assume that
\begin{description}
\item (i) there exists finite constants $\bar{c}$ and $C^\star
(>0) $ such that
\begin{equation}
\lim_{n \rightarrow ^\infty} \bar{c}_n = \bar{c}  \;\; \;
\mbox{and} \; \; \; \lim_{n \rightarrow \infty}  n^{-1}
{C_n^\star}^2 = {C^\star}^2 \label{limcncstar} \end{equation} with
\begin{equation} \bar{c}_n = n^{-1} \sum_{i=1}^n c_i  \;\; \;
\mbox{and}  \;\; \; {C_n^\star}^2 = \sum_{i=1}^n c_i^2 -
n\bar{c}_n^2
\end{equation} both exist.

 \item(ii) the $c_i$'s are all bounded, so that by (i),
 \begin{equation} \max_{1 \leq i \leq n} (c_i -
 \bar{c}_n)^2/{C_n^\star}^2
 \rightarrow 0, \;\;\; \mbox{as} \;\; n \rightarrow \infty.
 \label{limmax} \end{equation}
\end{description}
Let $\psi : \Re \rightarrow \Re$ be nondecreasing and skew
symmetric score function. For any real numbers $a$ and $b$,
consider the statistics below
\begin{equation}
M_{n_1}(a,b) = \sum_{i=1}^n \psi(X_i-a-bc_i), \label{mn1}
\end{equation}
\begin{equation}
M_{n_2}(a,b) = \sum_{i=1}^n c_i \psi(X_i-a-bc_i). \label{mn2}
\end{equation} Let $\tilde{\beta}$ be the constrained M-estimator
of $\beta$ when $\theta=0$, that is, $\tilde{\beta}$ is the
solution of $M_{n_2}(0,b)=0$ and it may be conveniently be
expressed as
\begin{equation}
\tilde{\beta} = [ \underbrace{\sup \{b:M_{n_2}(0,b) >0
\}}_{b_1}\;+\;  \underbrace{\inf\{b:M_{n_2}(0,b) <0 \}}_{b_2}]/2.
\label{eq:betasupinfL}
\end{equation}
Any value $b_1 < b < b_2$ can serve as the estimate of $M_n(0,b).$
However $\tilde{\beta}$, the centroid (and the median) of the
interval $[b_1,\;b_2]$ achieves the smallest maximum bias among
all translation invariant functionals (Huber, 1981 p.75), hence it
is a robust estimator with optimum robustness properties.
Similarly, let $\tilde{\theta}$ be the constrained M-estimator of
$\theta$ when $\beta=0$, that is, $\tilde{\theta}$ is the solution
of $M_{n_1}(a,0)=0$ and conveniently be expressed as
\begin{equation}
\tilde{\theta}= [ \sup \{ a:M_{n_1}(a,0) >0 \} \;+\; \inf
\{a:M_{n_1}(a,0) <0 \} ]/2. \label{eq:thetasupinfL}
\end{equation} By the same argument as above, $\tilde{\theta}$ is a robust
estimator.

The preliminaries notations and assumptions in this section are
used to develop the tests and construct the asymptotic power
functions of the UT, RT and PTT. In the next section, the UT and
RT are proposed. The PT (testing on slope) is introduced and the
PTT is constructed.

\section{The UT, RT and PTT}
\setcounter{equation}{0} Sen (1982) shows that the asymptotic
distribution of
\begin{equation} n^{-1/2}M_{n_2}(\tilde{\theta},0) \rightarrow N(0,\sigma_0^2
{C^\star}^2) \label{asmn2} \end{equation} under
$H_0^{(1)}:\beta=0.$ The consistency of $[S_n^{(3)}]^2=n^{-1}
\sum_{i=1}^n \psi^2(X_i-\tilde{\theta})$ as an estimator of
$\sigma_0^2$ follows from Jur$\check{e}$ckov$\acute{a}$ and Sen (1981). Hence, a test statistic
$A_n=M_{n_2}(\tilde{\theta},0)[C_n^\star\; S_n^{(3)}]^{-1}$ is
proposed by Sen (1982). The advantage of this test statistic
(score-type M-test) is that it does not require the computation of
the M-estimates or the estimation of functional $\gamma.$

By the same way, it is easy to show that the asymptotic distribution
of
\begin{equation} n^{-1/2}M_{n_1}(0,\tilde{\beta})
\rightarrow N(0,\sigma_0^2 {C^\star}^2/ \{{C^\star}^2
+\bar{c}^2\}) \label{asmn1} \end{equation} under
$H_0^\star:\theta=0.$ By the same token, the consistency of
$[S_n^{(1)}]^2=n^{-1} \sum_{i=1}^n \psi^2(X_i-\tilde{\beta}c_i)$
as an estimator of $\sigma_0^2$ follows.

These two asymptotic distributions results are useful to construct
suitable test in formulating the asymptotic power function for
testing the intercept after pre-testing. We are primarily
concerned with statistical tests for the parameter $\theta$ as
well as $\beta.$ In essence we need to consider four test
functions correspond to the four proposed tests.

\subsection{The unrestricted test (UT)}

If $\beta$ is unspecified, the designated test function is
$\phi_n^{UT}$ with the null hypothesis $H_0^\star: \theta =0$. The
testing for $\theta$ involves the elimination of the nuisance
parameter $\beta$. It follows that $M_{n_2}(0,b)$ is decreasing if
$b$ is increasing (Jur$\check{e}$ckov$\acute{a}$ and Sen, 1996,
p.85) and under local hypothesis, $H_0^{(1)}: \beta =0$,
$M_{n_2}(0,0)$ has expectation 0. Then let  \[ \tilde{\beta} =
(\sup\{b:M_{n_2}(0,b)
>0\} + \inf \{b:M_{n_2}(0,b) < 0 \})/2. \] Then $\tilde{\beta}$
is a translation invariant and robust estimator of $\beta$.

We consider the test statistic $T_n^{UT}=M_{n_1}(0,\tilde{\beta})$
where under $H_0^\star$, as $n \rightarrow \infty,$
\[ \frac{T_n^{UT}}{ \sqrt{C_n^{(1)}{S_n^{(1)}}^2}} \rightarrow
N(0,1)\] with $C_n^{(1)} = n - n^2 \bar{c}_n^2/\sum c_i^2=
n{C_n^\star}^2/({C_n^\star}^2+n\bar{c}_n^2)$ and
$[S_n^{(1)}]^2=\sum \psi^2(x_i-\tilde{\beta}c_i)/n.$

\subsection{The restricted test (RT)}

If $\beta =0$, the designated test function is $\phi_n^{RT}$ for
testing the null hypothesis $H_0^\star: \theta =0$. The proposed
test statistic is $T_n^{RT} = M_{n_1}(0,0).$ Note that for large
$n$, under $H_0:\theta=0,\beta=0$,
\begin{equation} n^{-1/2} T_n^{RT} = n^{-1/2} M_{n_1}(0,0) \rightarrow N(0,\sigma_0^2), \label{asmn} \end{equation}  where
$\sigma_0^2 = \int_{-\infty}^{\infty} \psi^2(x)dF(x)$ (see Sen,
1982, eq 3.7).

\subsection{The pre-test test (PTT)}

In this section, test on slope is proposed first and followed by
the construction of the ultimate test for testing intercept.

For the preliminary test on the slope, the test function,
$\phi_n^{PT}$ is designed to test the null hypothesis $H_0^{(1)}:
\beta =0$. The proposed test statistic is
$T_n^{PT}=M_{n_2}(\tilde{\theta},0)$ where \[ \tilde{\theta} = (
\sup\{a:M_{n_1}(a,0)>0\} + \inf \{a:M_{n_1}(a,0) < 0 \})/2 \] is a
robust estimator. Under $H_0^{(1)}$, as $n \rightarrow \infty,$
\[\frac{T_n^{PT}}{ \sqrt{C_n^{(3)}{S_n^{(3)}}^2}} \rightarrow N(0,1),\] where
$C_n^{(3)} = \sum c_i^2 -n \bar{c}_n^2 = {C_n^\star}^2 $ and $
[S_n^{(3)}]^2=\sum \psi^2(x_i-\tilde{\theta})/n.$

The consistency of $[S_n^{(1)}]^2$, $[S_n^{(2)}]^2=\sum
\psi^2(x)/n$ and $[S_n^{(3)}]^2$ as estimators of $\sigma_0^2$
follows by law of large number (Jur$\check{e}$ckov$\acute{a}$ and
Sen, 1981).

Now, we are in a position to formulate a test function
$\phi_n^{PTT}$ to test $H_0^\star:\theta =0$ following a
preliminary test on $\beta.$ First, we consider the case where all
of $\phi_n^{(j)},\;j=1,2,3$ are one-sided test. Let us choose
positive numbers $\alpha_j\;(0<\alpha_j <1)$ and real values
$\ell_{n,\alpha_j}^{(j)},\;j=1,2,3,$ such that for large sample
size,
\begin{equation}
P[T_n^{UT} > \ell_{n,\alpha_1}^{UT}|H_0^\star:\theta =0] =
\alpha_1, \label{alpha2} \end{equation}
\begin{equation} P[T_n^{RT}
> \ell_{n,\alpha_2}^{RT}|H_0:\theta=0,\beta=0 ] = \alpha_2,
\label{alpha1} \end{equation}
\begin{equation}
P[T_n^{PT} > \ell_{n,\alpha_3}^{PT}|H_0^{(1)}:\beta=0] = \alpha_3,
\label{alpha3} \end{equation} where $\ell_{n,\alpha_j}^{(j)}$ is
the critical value of $T_n^{(j)}$ at the $\alpha_j$ level of
significance. Let $\Phi(x)$ be the standard normal cumulative
distribution function, then
\begin{equation} \Phi(\tau_\alpha)=1-\alpha,\;\;\mbox{for}\;\;0<\alpha<1.\label{g} \end{equation} Using
equations
(\ref{asmn}), (\ref{alpha1}) and (\ref{g}), we obtain
\[   \begin{array}{rcl} \vspace{10pt} 1-\alpha_2 &=& P[T_n^{RT}
\leq \ell_{n,\alpha_2}^{RT}] \\ &=&\vspace{10pt} P \left
[n^{-1/2}T_n^{RT}/ \sqrt{{S_n^{(2)}}^2} \leq n^{-1/2}
\ell_{n,\alpha_2}^{RT}/\sqrt{{S_n^{(2)}}^2} \; \right ]\\ &
\rightarrow & P\left [n^{-1/2}T_n^{RT}/ \sqrt{\sigma_0^2} \leq
n^{-1/2} \ell_{n,\alpha_2}^{RT}/ \sqrt{\sigma_0^2} \;\right ]=
\Phi(n^{-1/2}
\ell_{n,\alpha_2}^{RT}/ \sigma_0 ),\\
\end{array} \] where $ {S_n^{(2)}}^2 = \sum \psi^2(x_i)/n.$
Thus as $n \rightarrow \infty$ we have \begin{equation} \frac
{n^{-1/2} \ell_{n,\alpha_2}^{RT}}{\sqrt{{S_n^{(2)}}^2}}
\rightarrow \tau_{\alpha_2} =\frac
{n^{-1/2}\ell_{n,\alpha_2}^{RT}}{\sqrt{\sigma_0^2}}\;\;
\mbox{(say).}\label{tau1}
\end{equation}
By the same way, using (\ref{asmn1}) and (\ref{alpha2}), we
observe that as $n \rightarrow \infty$,
\begin{equation} \frac
{n^{-1/2} \ell_{n,\alpha_1}^{UT}}{\sqrt{ {S_n^{(1)}}^2
C_n^{(1)}/n}} \rightarrow \tau_{\alpha_1} =\frac
{n^{-1/2}\ell_{n,\alpha_1}^{UT}}{\sqrt{\sigma_0^2 {C^\star}^2 /
({C^\star}^2+\bar{c}^2)}} \;\;\mbox{(say)},\label{tau2}
\end{equation}
 where $ {S_n^{(1)}}^2 = \sum \psi^2(x_i - \tilde{\beta}c_i)/n,$
and $C_n^{(1)}=n - n^2\bar{c}_n^2/\sum c_i^2.$ Also by
(\ref{asmn2}) and (\ref{alpha3}), as $n \rightarrow \infty,$
\begin{equation} \frac
{n^{-1/2} \ell_{n,\alpha_3}^{PT}}{\sqrt{{S_n^{(3)}}^2
{C_n^\star}^2 /n}} \rightarrow \tau_{\alpha_3} =\frac
{n^{-1/2}\ell_{n,\alpha_3}^{PT}}{\sqrt{\sigma_0^2 {C^\star}^2}  }
\;\; \mbox{(say),} \label{tau3}
\end{equation}
 where $ {S_n^{(3)}}^2 = \sum \psi^2(x_i - \tilde{\theta})/n,$
and ${C_n^\star}^2= \sum c_i^2 - n\bar{c}_n^2.$

Now we may write
\begin{equation}
\phi_n^{PTT} = I \left [ (T_n^{PT} \leq \ell_{n,\alpha_3}^{PT},
T_n^{RT}
> \ell_{n,\alpha_2}^{RT}) \;\;\mbox{or} \;\;(T_n^{PT}
>
\ell_{n,\alpha_3}^{PT}, T_n^{UT} > \ell_{n,\alpha_1}^{UT} ) \right
]
\end{equation}
as the test function for testing $H_0^\star:\theta=0$ after a
pre-test on $\beta.$ Note that $I(A)$ stands for the indicator
function of the set $A.$ It takes value 1 if $A$ occurs, otherwise
0. The function enables us to define the power of the test
$\phi_n^{PTT}$, that is given by
\begin{eqnarray} \Pi_n^{PTT}(\theta) &=& E(\phi_n^{PTT}|\theta) \nonumber\\ &=& P [T_n^{PT} \leq \ell_{n,\alpha_3}^{PT},
T_n^{RT} > \ell_{n,\alpha_2}^{RT}|\theta]+ P[T_n^{PT}
> \ell_{n,\alpha_3}^{PT}, T_n^{UT} >
\ell_{n,\alpha_1}^{UT}|\theta]. \label{eqpow} \end{eqnarray} In
general, the power of the test $\phi_n^{PTT}$ depends on
$\alpha_1,\alpha_2,\alpha_3, \theta,n$ as well as $\beta.$ Note
that the size of the ultimate test $\alpha_n^{PTT}$ is a special
case of the power of the test when $\theta=0.$ Since the nuisance
parameter $\beta$ is unknown, but, suspected to be close to 0, it
is of interest to study the dependence of both $\alpha_n^{PTT}$
and $\Pi_n^{PTT}(\theta)$ on $\beta$ (close to 0).

\section{Asymptotic distribution of $T_n^{UT}$, $T_n^{RT}$ and $T_n^{PT}$}
\setcounter{equation}{0}
This section is devoted to the asymptotic
distribution theory of statistics involved in proposing the PTT.
The asymptotic joint distributions of $\left[T_n^{UT},
T_n^{PT}\right]$ and $\left [T_n^{RT}, T_n^{PT}\right ]$ are
derived. The results are used in the next section in the
construction of the power function of the UT, RT and PTT.

Let $\{K_n\}$ be a sequence of alternative hypotheses, where
\begin{equation}
K_n: (\theta, \beta)=(n^{-1/2}\lambda_1, n^{-1/2}\lambda_2),
\label{Kn}
\end{equation}
with $\lambda_1,\lambda_2$ are (fixed) real numbers.

Interested readers are referred to Jur$\check{e}$ckov$\acute{a}$
(1977), Sen (1982) and Jur$\check{e}$ckov$\acute{a}$ and Sen
(1996, p.221) for the following asymptotic properties:
\begin{description}
\item{(i)} under $H_0:\theta = 0, \beta = 0,$ as $n$ grows large,
\begin{equation} n^{-1/2} \left (
\begin{array} {c} M_{n_1}(0,0)\\  M_{n_2}(0,0)\end{array} \right )
\rightarrow N_2 \left ( \left [\begin{array} {c} 0 \\0 \end{array} \right ], \sigma_0^2 \left ( \begin{array} {cc} 1 & \bar{c} \\
\bar{c} & {C^\star}^2+\bar{c}^2 \\ \end{array} \right )   \right
), \label{eq:3.7S1982}
\end{equation}
where $N_2(\cdot\;,\cdot\;)$ represents a bivariate normal
distribution with appropriate parameters.
 \item{(ii)} under $H_0:\theta = 0, \beta =
0,$
\begin{equation}\begin{array} {c} \vspace{10pt}
\sup \{n^{-1/2}| M_{n_1}(a,b) - M_{n_1}(0,0) + n\gamma (a + b
\bar{c})|: \\|a| \leq n^{-1/2}K,|b| \leq n^{-1/2}K \}
\rightarrow_p 0
\end{array} \label{eq:juri1}
\end{equation} as $n \rightarrow \infty$ and $K$ is a positive constant.

\item{(iii)} under $H_0:\theta = 0, \beta = 0,$
\begin{equation}
\begin{array} {c} \vspace{10pt} \sup \{  n^{-1/2}|M_{n_2}(a,b) - M_{n_2}(0,0) + n\gamma \{a\bar{c} + b
({C^\star}^2+\bar{c}^2)\} |: \\ \hspace{150pt}|a| \leq n^{-1/2}K,
|b| \leq n^{-1/2}K \} \rightarrow_p 0
 \end{array} \label{eq:juri2}
\end{equation} as $n \rightarrow \infty$ and $K$ is a positive constant.
The above convergence is in probability, means the sequences of
random variables converges in probability to a fix value (0).
\end{description}

An important concept that dominates the asymptotic theory of
statistics is the \textit{contiguity} of probability measures
(Jur$\check{e}$ckov$\acute{a}$ and Sen, 1996, p.61). Contiguity
arguments are a technique to obtain the limit distribution of a
sequence of statistics under the alternative hypothesis from a
limiting distribution under the null hypothesis (c.f. van der
Vaart, 1998 p.85). Let $\{P_n\}$ and $\{Q_n\}$ be two sequence of
probability measures defined in a measure spaces $(\Omega_n, B_n,
\mu_n).$ In the LeCam's first lemma (H$\acute{a}$jek et al., 1999,
p.251), if \[ \log L_n \stackrel{D}\rightarrow
N(-\frac{1}{2}\;\sigma^2,\sigma^2)\;\;(\mbox{under}\;\;
\{P_n\}),\] then $\{Q_n\}$ is contiguous to $\{P_n\}.$ Here the
likelihood ratio statistic $L_n$ is given by \[L_n =\left \{
\begin{array} {ccc} q_n/p_n &\mbox{for}&p_n>0 \\ 1 & \mbox{for}&p_n=q_n=0 \\ \infty  & \mbox{for}& 0=p_n<q_n, \end{array}
\right.
\] where $\{p_n, q_n \}$ are the sequence of simple hypothesis densities. In the LeCam's third lemma (H$\acute{a}$jek et al., 1999, p.257), if
 \[ \left [\begin{array}{c}T_n \\ \log L_n \end{array} \right ] \stackrel{D}\rightarrow
N_2 \left ( \left [\begin{array} {c} \mu_1 \\ \mu_2 \end{array} \right ], \left [ \begin{array} {cc} \sigma_{11}& \sigma_{12} \\
\sigma_{21} & \sigma_{22} \\ \end{array} \right ] \right )\;\;(\mbox{under}\;\; \{P_n\}),\]
where $T_n$ is a statistic with $\mu_2=-\frac{1}{2} \sigma_{22}$, then
\[T_n\stackrel{D}\rightarrow N(\mu_1+\sigma_{12},\;\sigma_{11})\;\;(\mbox{under}\;\;
\{Q_n\}).\] The LeCam's second lemma (H$\acute{a}$jek et al.,
1999, p.253) gives conditions when $\log L_n
\stackrel{D}\rightarrow N(-\frac{1}{2}\;\sigma^2,\sigma^2)$.

The concept of contiguity is more popular in R-estimation (rank
statistic) than in M-estimation. However, Sen (1982) uses the
contiguity of probability measures under
$H_n:\beta=n^{-1/2}\lambda$ to those under $H_0^\prime:\beta=0$ to
find the asymptotic distribution of $n^{-1/2}\left
[M_{n_1}(\theta,0),M_{n_2}(\theta,0) \right]$ under $H_n.$ In this
paper, the contiguity concept is utilized to find the asymptotic
distributions of statistics $ n^{-1/2}[\;T_n^{RT}, T_n^{PT}]$ and
$n^{-1/2}[\;T_n^{UT}, T_n^{PT}]$ under $K_n.$

\subsection{Asymptotic distribution of $ n^{-1/2} T_n^{RT}$ and $n^{-1/2} T_n^{PT}$}

 Following Jur$\check{e}$ckov$\acute{a}$
and Sen (1996, p.259), let $\{P_n\}$ and $\{Q_n\}$ denote the
probability distributions with the densities $p_n=\prod_{i=1}^n
f(X_i)$ and $q_n=\prod_{i=1}^n f(X_i-n^{-1/2}\lambda_1 - t
c_{ni})$ of the null hypothesis $H_0$ and the alternative
hypothesis $K_n,$ respectively, where
$c_{ni}=(c_i-\bar{c}_n)/C_n^\star,\;\;i=1,\ldots,n.$ Note that
under (\ref{regress}), (\ref{limcncstar}), (\ref{limmax}) and
(\ref{Kn}), the contiguity of the sequence of probability measures
under $\{K_n\} $ to those under $H_0$ follows from LeCam's first
and second lemmas (H$\acute{a}$jek et al., 1999, Chapter 7). We
are interested in the asymptotic distribution of the joint
statistics $\left [n^{-1/2}T_n^{RT}, n^{-1/2}T_n^{PT}\right ].$
Here convergence of $\left [n^{-1/2}T_n^{RT},
n^{-1/2}T_n^{PT}\right ] + \Upsilon \rightarrow [0,0]$ under $H_0$
implies $\left [n^{-1/2}T_n^{RT}, n^{-1/2}T_n^{PT} \right] +
\Upsilon \rightarrow [0,0]$ under $\{K_n\}$ since the probability
measures under $\{K_n\}$ is contiguous to that of under ${H_0}$
(c.f. Saleh, 2006, p.44). Here, $\Upsilon$ is a known vector.

Under $H_0:\theta = 0, \; \beta = 0,$ with relation to
(\ref{eq:juri1}) and (\ref{eq:juri2}),
\begin{eqnarray} \left [\begin{array} {cc}
n^{-1/2} M_{n_1}(0,0) \\
n^{-1/2} M_{n_2}(\tilde{\theta},0)
\end{array}  \right ]
-
\left [ \begin{array} {c} n^{-1/2} M_{n_1}(0,0)  \\ n^{-1/2} M_{n_2}(0,0) \\
\end{array} \right ] + \left [ \begin{array} {c}
0 \\
n^{1/2} \gamma \tilde{\theta} \bar{c}  \\
\end{array} \right]
 \rightarrow_p \left [\begin{array} {c} 0 \\ 0 \end{array} \right
 ]
\label{eq:maj}\end{eqnarray} or equivalently
\begin{equation}
n^{-1/2} M_{n_1}(0,0) = n^{-1/2}M_{n_1}(0,0)+ o_p(1)
\end{equation} and
\begin{equation}
n^{-1/2} M_{n_2}(\tilde{\theta},0) = n^{-1/2}M_{n_2}(0,0) -n^{1/2}
\gamma \tilde{\theta} \bar{c} + o_p(1). \label{eq:thetatilde}
\end{equation}Note also that under $H_0:\theta = 0, \; \beta = 0,$
\begin{equation} n^{-1/2}M_{n_1}(\tilde{\theta},0) = n^{-1/2} M_{n_1}(0,0) -
n^{1/2} \gamma \tilde{\theta}  + o_p(1). \label{41} \end{equation}
Recalling definition (\ref{eq:thetasupinfL}), the equation
(\ref{41}) reduces to
\begin{equation} n^{-1/2} M_{n_1}(0,0) = n^{1/2} \gamma \tilde{\theta} + o_p(1),
\end{equation} and hence equation (\ref{eq:thetatilde}) becomes
\begin{equation}
n^{-1/2} M_{n_2}(\tilde{\theta},0) = n^{-1/2}M_{n_2}(0,0) -
n^{-1/2} M_{n_1}(0,0)\bar{c} + o_p(1). \label{eq:r1}
\end{equation} Therefore, under $H_0$, equation (\ref{eq:maj}) becomes

\begin{eqnarray} \vspace{5pt} \left [ \begin{array} {c}
n^{-1/2} M_{n_1}(0,0) \\
n^{-1/2} M_{n_2}(\tilde{\theta},0) \nonumber\\
\end{array}  \right ]
-\left [ \begin{array} {c} n^{-1/2} M_{n_1}(0,0) \\ n^{-1/2} M_{n_2}(0,0) - n^{-1/2} M_{n_1}(0,0)\bar{c}\\
\end{array} \right ] \hspace{50pt}\\ \hspace{50pt}
\vspace{5pt} = \left [\begin{array} {c}
n^{-1/2} M_{n_1}(0,0) \\
n^{-1/2} M_{n_2}(\tilde{\theta},0)\\
\end{array}  \right ]
- \left [ \begin{array} {cc} 1 & 0 \\ - \bar{c} & 1 \\ \end{array} \right ] \left [ \begin{array} {c} n^{-1/2} M_{n_1}(0,0)\\ n^{-1/2} M_{n_2}(0,0) \\
\end{array} \right ]
 \rightarrow_p \left [\begin{array} {c} 0 \\ 0 \end{array} \right
 ]. \label{eq:bsr1} \end{eqnarray}
Now utilizing the contiguity of probability measures under
$\{K_n\}$ to those under $H_0$, the equation (\ref{eq:bsr1})
implies that
\[ \left [
\begin{array} {c}
n^{-1/2} M_{n_1}(0,0) \\
n^{-1/2} M_{n_2}(\tilde{\theta},0)\\
\end{array}  \right ] \] under $\{K_n\}$ is asymptotically equivalent to the
random vector \[ \left [ \begin{array} {cc} 1 & 0 \\ - \bar{c} & 1  \\ \end{array} \right ] \left [ \begin{array} {c} n^{-1/2} M_{n_1}(0,0)\\ n^{-1/2} M_{n_2}(0,0) \\
\end{array} \right ] \] under $H_0.$
But the asymptotic distribution of the above random vector under
$\{K_n\}$ is the same as
\[ \left [ \begin{array} {cc} 1 & 0 \\ - \bar{c} & 1 \\ \end{array} \right ]
\left [ \begin{array} {c} n^{-1/2} M_{n_1}(-n^{-1/2}\lambda_1,-n^{-1/2} \lambda_2)\\ n^{-1/2} M_{n_2}(-n^{-1/2}\lambda_1,-n^{-1/2}\lambda_2) \\
\end{array} \right ] \] under $H_0$ by the fact that the distribution of $M_{n_1}(a,b)$
under $\theta = a, \beta= b$ is the same as that of
$M_{n_1}(\theta -a, \beta-b)$ under $\theta =0, \beta =0,$ and
similarly for $M_{n_2}(0,0)$ (c.f. Saleh, 2006 p.332).

Note that under $H_0:\theta=0, \beta=0$, with relation to
(\ref{eq:juri1}) and (\ref{eq:juri2}),
\begin{equation}
\begin{array} {c}  \vspace{10pt} \left [\begin{array} {cc}
n^{-1/2} M_{n_1}(-n^{-1/2}\lambda_1,-n^{-1/2} \lambda_2)\\
n^{-1/2} M_{n_2}(-n^{-1/2}\lambda_1,-n^{-1/2}\lambda_2)
\end{array}  \right ]
-
\left [ \begin{array} {c} n^{-1/2} M_{n_1}(0,0)  \\ n^{-1/2} M_{n_2}(0,0) \\
\end{array} \right ] \hspace{50pt}\\ \hspace{50pt} - \left [ \begin{array} {c}
 \gamma (\lambda_1 + \lambda_2 \bar{c})\\
 \gamma \{\lambda_1\bar{c}+ \lambda_2({C^\star}^2 +\bar{c}^2)\}  \\
\end{array} \right]
 \rightarrow_p \left [\begin{array} {c} 0 \\ 0 \end{array} \right
 ].
\end{array}
\end{equation} Hence, by equation (\ref{eq:3.7S1982}), under $H_0$,
\begin{eqnarray} &&\left [ \begin{array} {c} n^{-1/2} M_{n_1}(-n^{-1/2}\lambda_1,-n^{-1/2} \lambda_2)\\ n^{-1/2} M_{n_2}(-n^{-1/2}\lambda_1,-n^{-1/2}\lambda_2) \\
\end{array} \right ] \nonumber \\ &&\;\;\;\;\;\;\;\;\;\rightarrow N_2 \left( \left ( \begin{array} {c}
\gamma (\lambda_1+ \lambda_2 \bar{c})  \\
\gamma \{\lambda_1 \bar{c} + \lambda_2 ({C^\star}^2+\bar{c}^2) \}  \\
\end{array} \right) ,\;\sigma_0^2
\left ( \begin{array} {cc} 1 & \bar{c} \\
\bar{c} & {C^\star}^2 + \bar{c}^2 \\ \end{array}  \right ) \right
). \label{eq:dis1} \end{eqnarray} Thus, the distribution of
\[ \left [ \begin{array}{c} n^{-1/2}T_n^{RT} \\ n^{-1/2}T_n^{PT} \end{array} \right ] =\left [
\begin{array} {c}
n^{-1/2} M_{n_1}(0,0) \\
n^{-1/2} M_{n_2}(\tilde{\theta},0)\\
\end{array}  \right ] \] under $\{K_n\}$
is bivariate normal with mean vector
\[ \left [ \begin{array} {cc} 1 & 0 \\ - \bar{c} & 1 \\
\end{array} \right ] \left [ \begin{array} {c}
 \gamma (\lambda_1+ \lambda_2 \bar{c})  \\
\gamma \{\lambda_1 \bar{c} + \lambda_2 ({C^\star}^2+\bar{c}^2) \}
\end{array} \right ] = \left [ \begin{array} {c}
 \gamma (\lambda_1+ \lambda_2 \bar{c})  \\
\gamma \lambda_2 {C^\star}^2
\end{array} \right ] \] and covariance matrix
\begin{equation} \left [ \begin{array} {cc} 1 & 0 \\ - \bar{c} & 1 \\
\end{array} \right ] \sigma_0^2
\left ( \begin{array} {cc} 1 & \bar{c} \\
\bar{c} & {C^\star}^2 + \bar{c}^2 \\ \end{array}  \right ) \left [ \begin{array} {cc} 1 & 0 \\ - \bar{c} & 1 \\
\end{array} \right ]' = \sigma_0^2 \left [ \begin{array} {cc} 1 & 0 \\ 0 & {C^\star}^2 \\
\end{array} \right ]. \label{joint1} \end{equation}
Since the two statistics $n^{-1/2}T_n^{RT}$ and $n^{-1/2}T_n^{PT}$
are uncorrelated, asymptotically, they are independently
distributed normal variables.

\subsection{Asymptotic distribution of $n^{-1/2}T_n^{UT}$ and $n^{-1/2}T_n^{PT}$}

Under $H_0:\theta = 0, \beta = 0,$ with relation to
(\ref{eq:juri1}) and (\ref{eq:juri2}), as $n \rightarrow \infty,$

\begin{eqnarray}  \vspace{10pt} \left [\begin{array} {cc}
n^{-1/2} M_{n_1}(0,\tilde{\beta}) \\
n^{-1/2} M_{n_2}(\tilde{\theta},0)
\end{array}  \right ]
-
\left [ \begin{array} {c} n^{-1/2} M_{n_1}(0,0)  \\ n^{-1/2} M_{n_2}(0,0) \\
\end{array} \right ]   \left [ \begin{array} {c}
n^{1/2} \gamma \tilde{\beta} \bar{c} \\
n^{1/2} \gamma \tilde{\theta} \bar{c} \\
\end{array} \right]
 \rightarrow_p \left [\begin{array} {c} 0 \\ 0 \end{array} \right
 ]
\label{eq:maj2}\end{eqnarray} or equivalently
\begin{equation}
n^{-1/2} M_{n_1}(0,\tilde{\beta}) =
n^{-1/2}M_{n_1}(0,0)-n^{1/2}\gamma \tilde{\beta} \bar{c}  + o_p(1)
\label{eq:betatilde2}\end{equation} and
\begin{equation}
n^{-1/2} M_{n_2}(\tilde{\theta},0) = n^{-1/2}M_{n_2}(0,0) -n^{1/2}
\gamma \tilde{\theta} \bar{c} + o_p(1). \label{eq:thetatilde2}
\end{equation}
Again under $H_0:\theta = 0, \beta = 0,$
\begin{equation} n^{-1/2}M_{n_2}(0,\tilde{\beta}) = n^{-1/2} M_{n_2}(0,0) -
n^{1/2} \gamma \tilde{\beta}({C^\star}^2+\bar{c}^2)  + o_p(1).
\end{equation} But by definition (\ref{eq:betasupinfL}), we have
\begin{equation} n^{-1/2} M_{n_2}(0,0) = n^{1/2} \gamma \tilde{\beta}({C^\star}^2+\bar{c}^2) + o_p(1),
\end{equation} so equation (\ref{eq:betatilde2}) becomes
\begin{equation}
n^{-1/2} M_{n_1}(0,\tilde{\beta}) = n^{-1/2}M_{n_1}(0,0) -
n^{-1/2} M_{n_2}(0,0)\bar{c}/({C^\star}^2+\bar{c}^2) + o_p(1).
\label{eq:r2}
\end{equation} Hence, under $H_0$, using equations (\ref{eq:r1}) and
(\ref{eq:r2}), equation (\ref{eq:maj2}) becomes

\begin{eqnarray}  \vspace{5pt} \left [ \begin{array} {c}
n^{-1/2} M_{n_1}(0,\tilde{\beta}) \\
n^{-1/2} M_{n_2}(\tilde{\theta},0)\\
\end{array}  \right ]
-\left [ \begin{array} {c} n^{-1/2}M_{n_1}(0,0) -
n^{-1/2} M_{n_2}(0,0)\bar{c}/({C^\star}^2+\bar{c}^2) \\ n^{-1/2} M_{n_2}(0,0) - n^{-1/2} M_{n_1}(0,0)\bar{c}\\
\end{array} \right ] \hspace{50pt} \nonumber\\ \hspace{1pt}
\vspace{5pt} = \left [\begin{array} {c}
n^{-1/2} M_{n_1}(0,\tilde{\beta}) \\
n^{-1/2} M_{n_2}(\tilde{\theta},0)\\
\end{array}  \right ]
- \left [ \begin{array} {cc} 1 & -\bar{c}/({C^\star}^2 + \bar{c}^2) \\ - \bar{c} & 1 \\ \end{array} \right ] \left [ \begin{array} {c} n^{-1/2} M_{n_1}(0,0)\\ n^{-1/2} M_{n_2}(0,0) \\
\end{array} \right ]
 \rightarrow_p \left [\begin{array} {c} 0 \\ 0 \end{array} \right
 ]. \label{eq:bsr2} \end{eqnarray}   Now by using the contiguity of probability measures under $\{K_n\}$
to those under ${H_0}$, the equation (\ref{eq:bsr2}) implies that
\[ \left [
\begin{array} {c}
n^{-1/2} M_{n_1}(0,\tilde{\beta}) \\
n^{-1/2} M_{n_2}(\tilde{\theta},0)\\
\end{array}  \right ] \] under $\{K_n\}$ is asymptotically equivalent to the
random vector \[ \left [ \begin{array} {cc} 1 & -\bar{c}/({C^\star}^2 + \bar{c}^2) \\ - \bar{c} & 1  \\ \end{array} \right ] \left [ \begin{array} {c} n^{-1/2} M_{n_1}(0,0)\\ n^{-1/2} M_{n_2}(0,0) \\
\end{array} \right ]. \]
But the asymptotic distribution of the above random vector under
$\{K_n\}$ is the same as
\[ \left [ \begin{array} {cc} 1 & -\bar{c}/({C^\star}^2 + \bar{c}^2) \\ - \bar{c} & 1  \\ \end{array} \right ]
\left [ \begin{array} {c} n^{-1/2} M_{n_1}(-n^{-1/2}\lambda_1,-n^{-1/2} \lambda_2)\\ n^{-1/2} M_{n_2}(-n^{-1/2}\lambda_1,-n^{-1/2}\lambda_2) \\
\end{array} \right ] \] under $H_0$ by the fact that the distribution of $M_{n_1}(a,b)$
under $\theta = a, \beta= b$ is the same as that of
$M_{n_1}(\theta -a, \beta-b)$ under $\theta =0, \beta =0$ and
similarly for $M_{n_2}(0,0)$ (c.f. Saleh, 2006 p.332).

Then it follows that by equation (\ref{eq:dis1}),
\[ \left [
\begin{array}{c} n^{-1/2}T_n^{UT} \\ n^{-1/2}T_n^{PT} \end{array} \right ]
=\left [
\begin{array} {c}
n^{-1/2} M_{n_1}(0,\tilde{\beta}) \\
n^{-1/2} M_{n_2}(\tilde{\theta},0)\\
\end{array}  \right ] \]
is bivariate normal with mean vector
\[ \left [ \begin{array} {cc} 1 & -\bar{c}/({C^\star}^2 + \bar{c}^2) \\ - \bar{c} & 1 \\
\end{array} \right ] \left [ \begin{array} {c}
 \gamma (\lambda_1+ \lambda_2 \bar{c})  \\
\gamma \{\lambda_1 \bar{c} + \lambda_2 ({C^\star}^2+\bar{c}^2) \}
\end{array} \right ] = \left [ \begin{array} {c}
 \gamma \lambda_1 {C^\star}^2/({C^\star}^2 + \bar{c}^2)  \\
\gamma \lambda_2 {C^\star}^2
\end{array} \right ] \] and covariance matrix
\begin{eqnarray}  \vspace{5pt} \left [ \begin{array} {cc} 1 & -\bar{c}/({C^\star}^2 + \bar{c}^2) \\ - \bar{c} & 1  \\
\end{array} \right ] \sigma_0^2
\left ( \begin{array} {cc} 1 & \bar{c} \\
\bar{c} & {C^\star}^2 + \bar{c}^2 \\ \end{array}  \right ) \left [ \begin{array} {cc} 1 & -\bar{c}/({C^\star}^2 + \bar{c}^2) \\ - \bar{c} & 1  \\
\end{array} \right ]' \nonumber\\
= \sigma_0^2 \left [ \begin{array} {cc} {C^\star}^2/({C^\star}^2+\bar{c}^2) & -\bar{c}\;{C^\star}^2/({C^\star}^2+\bar{c}^2) \\  -\bar{c}\;{C^\star}^2/({C^\star}^2+\bar{c}^2) & {C^\star}^2 \\
\end{array} \right ].\hspace{100pt}  \label{joint2} \end{eqnarray}
Clearly, the two test statistics $n^{-1/2}T_n^{UT}$ and
$n^{-1/2}T_n^{PT}$ are not independent, rather correlated.

\section{Asymptotic properties for UT, RT and PTT}
\setcounter{equation}{0} In this section, the asymptotic power
functions of $\phi_n^{UT},$ $\phi_n^{RT}$ and $\phi_n^{PTT}$ are
derived by using the results obtained in the previous sections.

Under $\{K_n\}$, the power function of $\phi_n^{PTT}$ is given by
\begin{eqnarray}  \Pi_n^{PTT}(\lambda_1,\lambda_2)=E(\phi_n^{PTT}|K_n)
&=&P[T_n^{PT} \leq \ell_{n,\alpha_3}^{PT},T_n^{RT} \geq
\ell_{n,\alpha_2}^{RT}|K_n] \nonumber\\
&& +\; P[T_n^{PT} \geq \ell_{n,\alpha_3}^{PT},T_n^{UT} \geq
\ell_{n,\alpha_1}^{UT}|K_n].
\end{eqnarray}
Note that
\begin{eqnarray}   && P[T_n^{PT} \leq
\ell_{n,\alpha_3}^{PT},T_n^{RT} >
\ell_{n,\alpha_2}^{RT}|K_n]\nonumber\\
&=& P\left [ \frac{n^{-1/2} T_n^{PT}-\gamma\lambda_2{C^\star}^2} {
\sqrt{{S_n^{(3)}}^2{C_n^\star}^2/n}} \leq \frac{n^{-1/2}
\ell_{n,\alpha_3}^{PT}-\gamma\lambda_2{C^\star}^2}{\sqrt{{S_n^{(3)}}^2{C_n^\star}^2/n}}
, \right. \nonumber \\
&& \left. \hspace{15pt}\frac{n^{-1/2} T_n^{RT}-\gamma(\lambda_1
+\lambda_2\bar{c})}{\sqrt{{S_n^{(2)}}^2}} > \frac{n^{-1/2}
\ell_{n,\alpha_2}^{RT}-\gamma(\lambda_1
+\lambda_2\bar{c})}{\sqrt{{S_n^{(2)}}^2}} \right ]\nonumber\\
&  \rightarrow & P\left [ \frac{n^{-1/2}
T_n^{PT}-\gamma\lambda_2{C^\star}^2}{\sqrt{\sigma_0^2
{C^\star}^2}} \leq \frac{n^{-1/2}
\ell_{n,\alpha_3}^{PT}-\gamma\lambda_2{C^\star}^2}{\sqrt{\sigma_0^2
{C^\star}^2}}
, \right. \nonumber\\
&& \left. \hspace{15pt}\frac{n^{-1/2} T_n^{RT}-\gamma(\lambda_1
+\lambda_2\bar{c})}{\sqrt{\sigma_0^2}} > \frac{n^{-1/2}
\ell_{n,\alpha_2}^{RT}-\gamma(\lambda_1
+\lambda_2\bar{c})}{\sqrt{\sigma_0^2}} \right ] \mbox{ as } n
\rightarrow \infty \label{eq61} \end{eqnarray} because the limit
of ${S_n^{(2)}}^2$ and ${S_n^{(3)}}^2$ are $\sigma_0^2$ and
${C_n^\star}^2/n \rightarrow C^\star$ as $n \rightarrow \infty.$

From equations (\ref{g}), (\ref{tau1}) and (\ref{tau3}) and
(\ref{joint1}), the probability statement in (\ref{eq61}) becomes
\[ \Phi(\tau_{\alpha_3}-\gamma\lambda_2
C^\star/\sigma_0)[1-\Phi(\tau_{\alpha_2}-\gamma(\lambda_1+\lambda_2
\bar{c})/\sigma_0)]. \] Note that $T_n^{RT}$ and $T_n^{PT}$ are
independent by equation (\ref{joint1}).

\noindent Define $d(q_1,q_2:\rho)$ to be the bivariate normal
probability integral for random variables $x$ and $y$,
\begin{equation}
d(q_1,q_2;\rho) = \frac{1}{2 \pi (1-\rho^2)^{1/2}}
\int_{q_1}^\infty \int_{q_2}^\infty \mbox{exp} \left \{
\frac{-(x^2 +y^2 - 2\rho xy)}{2(1-\rho^2)} \right \} dxdy,
\label{bi.int}
\end{equation} where $q_1,q_2$ are real numbers and $-1<\rho <1$. Here
$d(q_1,q_2;\rho)$ is the complement of the cumulative density
function of standard bivariate normal variable.

Since ${S_n^{(1)}}^2$ and ${S_n^{(3)}}^2$ both converge to
$\sigma_0^2,$
 and $C_n^{(1)}/n \rightarrow
{C^\star}^2/({C^\star}^2 + \bar{c}^2)$  as $n \rightarrow \infty,$
we observe that
\begin{eqnarray} && P[T_n^{PT} >
\ell_{n,\alpha_3}^{PT},T_n^{UT} >
\ell_{n,\alpha_1}^{UT}|K_n] \nonumber\\
&=& P\left [ \frac{n^{-1/2}
T_n^{PT}-\gamma\lambda_2{C^\star}^2}{\sqrt{{S_n^{(3)}}^2{C_n^\star}^2/n}}
> \frac{n^{-1/2}
\ell_{n,\alpha_3}^{PT}-\gamma\lambda_2{C^\star}^2}{\sqrt{{S_n^{(3)}}^2{C_n^\star}^2/n}}
, \right. \nonumber \\
&& \left. \hspace{15pt}\frac{n^{-1/2} T_n^{UT}-\gamma \lambda_1
{C^\star}^2/({C^\star}^2+\bar{c}^2)}{\sqrt{{S_n^{(1)}}^2C_n^{(1)}/n}}
> \frac{n^{-1/2} \ell_{n,\alpha_1}^{UT}-\gamma \lambda_1
{C^\star}^2/({C^\star}^2+\bar{c}^2)}{\sqrt{{S_n^{(1)}}^2C_n^{(1)}/n}}
\right ]\nonumber\\ &  \rightarrow & P\left [ \frac{n^{-1/2}
T_n^{PT}-\gamma\lambda_2{C^\star}^2}{\sqrt{\sigma_0^2
{C^\star}^2}}
> \frac{n^{-1/2}
\ell_{n,\alpha_3}^{PT}-\gamma\lambda_2{C^\star}^2}{\sqrt{\sigma_0^2
{C^\star}^2}}
, \right. \nonumber\\
&& \left. \hspace{15pt}\frac{n^{-1/2} T_n^{UT}-\gamma \lambda_1
{C^\star}^2/({C^\star}^2+\bar{c}^2)}{\sqrt{\sigma_0^2{C^\star}^2/({C^\star}^2+\bar{c}^2)}}
> \frac{n^{-1/2} \ell_{n,\alpha_1}^{UT}-\gamma \lambda_1
{C^\star}^2/({C^\star}^2+\bar{c}^2)}{\sqrt{\sigma_0^2{C^\star}^2/({C^\star}^2+\bar{c}^2)}}
\right ] \label{pt} \end{eqnarray} $\mbox{ as } n \rightarrow
\infty.$ Further, the equation (\ref{pt}) is written as
\[ d(\tau_{\alpha_3}-\gamma\lambda_2
C^\star/\sigma_0, \tau_{\alpha_1}-\gamma \lambda_1 \sqrt{
{C^\star}^2/({C^\star}^2 + \bar{c}^2)} /\sigma_0;
-\bar{c}/\sqrt{{C^\star}^2 + \bar{c}^2}\;)
\]
by using equations (\ref{tau2}), (\ref{tau3}), (\ref{joint2}) and
(\ref{bi.int}). Note that $T_n^{UT}$ and $T_n^{PT}$ are not
independent because of (\ref{joint2}).

Hence, the power function of $\phi_n^{PTT}$ for the PTT becomes
\begin{eqnarray} && \Pi_n^{PTT}(\lambda_1,\lambda_2)= E(\phi_n^{PTT}|K_n) \rightarrow \Pi^{PTT}(\lambda_1,\lambda_2) \nonumber \\
&&= \Phi(\tau_{\alpha_3}-\gamma\lambda_2
C^\star/\sigma_0)[1-\Phi(\tau_{\alpha_2}-\gamma(\lambda_1+\lambda_2
\bar{c})/\sigma_0)]\;+ \nonumber \\&&\;\;\;\;\;
d(\tau_{\alpha_3}-\gamma\lambda_2 C^\star/\sigma_0,
\tau_{\alpha_1}-\gamma \lambda_1 \sqrt{ {C^\star}^2/({C^\star}^2 +
\bar{c}^2)} /\sigma_0; -\bar{c}/\sqrt{{C^\star}^2 + \bar{c}^2}\;).
\label{pistar}
\end{eqnarray}
Similarly, the power function of $\phi_n^{RT}$ for the RT is given
by
\begin{eqnarray} \Pi_n^{RT}(\lambda_1,\lambda_2)&=&E(\phi_n^{RT}|K_n) \nonumber\\
&=& P[T_n^{RT} > \ell_{n,\alpha_2}^{RT}|K_n] \nonumber\\
&=& P \left [ \frac{n^{-1/2} T_n^{RT} - \gamma(\lambda_1+\lambda_2
\bar{c})}{\sqrt{{S_n^{(2)}}^2}} > \frac{n^{-1/2}
\ell_{n,\alpha_2}^{RT} -
\gamma(\lambda_1+\lambda_2 \bar{c})}{\sqrt{{S_n^{(2)}}^2}} \right] \vspace{5pt} \nonumber\\
&\rightarrow& P \left [ \frac{n^{-1/2} T_n^{RT} -
\gamma(\lambda_1+\lambda_2 \bar{c})}{\sqrt{\sigma_0^2}}
> \frac{n^{-1/2} \ell_{n,\alpha_2}^{RT} -
\gamma(\lambda_1+\lambda_2 \bar{c})}{\sqrt{\sigma_0^2}} \right]
\label{57} \end{eqnarray} since ${S_n^{(2)}}^2\rightarrow
\sigma_0^2.$ Combining equations (\ref{g}) and (\ref{tau1}), the
power function of $\phi_n^{RT}$ becomes
\begin{equation} \Pi^{RT}(\lambda_1,\lambda_2)=1-\Phi(\tau_{\alpha_2} -\gamma(\lambda_1+\lambda_2
\bar{c})/\sigma_0).\label{pi1} \end{equation} Finally, the power
function of $\phi_n^{UT}$ for the UT is obtained as
\begin{eqnarray}\Pi_n^{UT}(\lambda_1,\lambda_2)&=&E(\phi_n^{UT}|K_n)
= P[T_n^{UT} > \ell_{n,\alpha_1}^{UT}|K_n] \nonumber\\
&=& P \left [ \frac{n^{-1/2} T_n^{UT} - \gamma \lambda_1
{C^\star}^2/({C^\star}^2+\bar{c}^2)}{\sqrt{{S_n^{(1)}}^2C_n^{(1)}/n}}
> \frac{n^{-1/2} \ell_{n,\alpha_1}^{UT} - \gamma \lambda_1
{C^\star}^2/({C^\star}^2+\bar{c}^2)}{\sqrt{{S_n^{(1)}}^2C_n^{(1)}/n}} \right]\vspace{5pt}\nonumber \\
&\rightarrow& P \left [ \frac{n^{-1/2} T_n^{UT} - \gamma \lambda_1
{C^\star}^2/({C^\star}^2+\bar{c}^2)}{\sqrt{\sigma_0^2{C^\star}^2/({C^\star}^2+\bar{c}^2)}
}
> \frac{n^{-1/2} \ell_{n,\alpha_1}^{UT} -
\gamma \lambda_1
{C^\star}^2/({C^\star}^2+\bar{c}^2)}{\sqrt{\sigma_0^2{C^\star}^2/({C^\star}^2+\bar{c}^2)}
} \right] \nonumber\\
\end{eqnarray}
since ${S_n^{(1)}}^2\rightarrow \sigma_0^2.$ Further the power
function for the UT is written as
\begin{equation} \Pi^{UT}(\lambda_1,\lambda_2)=1-\Phi(\tau_{\alpha_1}-\gamma \lambda_1 \sqrt{
{C^\star}^2/({C^\star}^2 + \bar{c}^2)} \;/\sigma_0)
\label{pi2} \end{equation} using equations (\ref{g}) and (\ref{tau2}).

The asymptotic power functions for the UT, RT and PTT that are
derived using M-test in this section are found to have the same
form as that derived by using the rank statistic by Saleh and Sen
(1982) though the methodology of M-estimation and R-estimation is
different. Therefore, the investigation on the properties of the
power of the M-test is similar to the power of the test based on
rank statistic.

\section{Asymptotic comparison}
\setcounter{equation}{0} This section gives analytic asymptotic
comparison of the power functions of the UT, RT and PTT.

\noindent If we consider $\bar{c}=0$ in equation (\ref{pistar}),
\begin{eqnarray}
\Pi^{PTT}(\lambda_1,\lambda_2)&=&\Phi(\tau_{\alpha_3}-\gamma\lambda_2
C^\star/\sigma_0)[1-\Phi(\tau_{\alpha_2}-\gamma
\lambda_1/\sigma_0)]+\; \nonumber \\&&
[1-\Phi(\tau_{\alpha_3}-\gamma\lambda_2 C^\star/\sigma_0)]
[1-\Phi(\tau_{\alpha_1}-\gamma \lambda_1/\sigma_0)].
\label{eq.pistar1}
\end{eqnarray}
Letting $\alpha_1=\alpha_2=\alpha$ and from equations (\ref{pi1}),
(\ref{pi2}) and (\ref{eq.pistar1}), we observe that the power
functions for the UT, RT and PTT are the same, i.e.
\begin{eqnarray}
\Pi^{UT}(\lambda_1,\lambda_2)=\Pi^{RT}(\lambda_1,\lambda_2)=\Pi^{PTT}(\lambda_1,\lambda_2)=
1-\Phi(\tau_{\alpha}-\gamma(\lambda_1+\lambda_2
\bar{c})/\sigma_0). \label{same} \end{eqnarray} From equations
(\ref{pistar}) and (\ref{pi1}),
\begin{eqnarray}
&& \Pi^{RT}(\lambda_1,\lambda_2)-\Pi^{PTT}(\lambda_1,\lambda_2)
\nonumber \\ &=& 1-\Phi(\tau_{\alpha_2}-\gamma(\lambda_1+\lambda_2
\bar{c})/\sigma_0) - \Phi(\tau_{\alpha_3}-\gamma\lambda_2
C^\star/\sigma_0)[1-\Phi(\tau_{\alpha_2}-\gamma(\lambda_1+\lambda_2
\bar{c})/\sigma_0)]\; \nonumber
\\ && -\;d(\tau_{\alpha_3}-\gamma\lambda_2 C^\star/\sigma_0,
\tau_{\alpha_1}-\gamma \lambda_1 \sqrt{ {C^\star}^2/({C^\star}^2 +
\bar{c}^2)} /\sigma_0; -\bar{c}/\sqrt{{C^\star}^2 + \bar{c}^2}\;)
\nonumber \\&=& d(\tau_{\alpha_3}-\gamma\lambda_2
C^\star/\sigma_0, \tau_{\alpha_2}-\gamma(\lambda_1+\lambda_2
\bar{c})/\sigma_0; 0 )\nonumber
\\&&-d(\tau_{\alpha_3}-\gamma\lambda_2 C^\star/\sigma_0,
\tau_{\alpha_1}-\gamma \lambda_1 \sqrt{ {C^\star}^2/({C^\star}^2 +
\bar{c}^2)} /\sigma_0; -\bar{c}/\sqrt{{C^\star}^2 + \bar{c}^2}\;).
\label{eq.pistar2}
\end{eqnarray}
Letting $\alpha_1=\alpha_2=\alpha,$ $\bar{c}>0,$ $\lambda_2>0$ and
$\lambda_1+\lambda_2 \bar{c}
> \lambda_1 \sqrt{ {C^\star}^2/({C^\star}^2+ \bar{c}^2)}$,
\begin{description} \item {\textit{Result }(i):}\hspace{10pt}$\Pi^{RT}(\lambda_1,\lambda_2)
> \Pi^{PTT}(\lambda_1,\lambda_2)$ from equation (\ref{eq.pistar2})
and  \item{\textit{Result
}(ii):}\hspace{10pt}$\Pi^{RT}(\lambda_1,\lambda_2) >
\Pi^{UT}(\lambda_1,\lambda_2)$ from equations (\ref{pi1}) and
(\ref{pi2}). \end{description} On the contrary, taking
$\alpha_1=\alpha_2=\alpha,$ $\bar{c}<0,$ $\lambda_2>0$ and
$\lambda_1+\lambda_2 \bar{c} < \lambda_1 \sqrt{
{C^\star}^2/({C^\star}^2+ \bar{c}^2)}$,
\begin{description} \item{\textit{Result }(iii):}\hspace{10pt}$\Pi^{RT}(\lambda_1,\lambda_2)
< \Pi^{PTT}(\lambda_1,\lambda_2)$  from equation
(\ref{eq.pistar2}) and \item{\textit{Result }(iv):}
\hspace{10pt}$\Pi^{RT}(\lambda_1,\lambda_2) <
\Pi^{UT}(\lambda_1,\lambda_2)$ from equations (\ref{pi1}) and
(\ref{pi2}).  \end{description}  From equations (\ref{pi1}) and
(\ref{pi2}), when $\lambda_1=\lambda_2=0$ and
$\alpha_1=\alpha_2=\alpha,$ we find $\Pi^{RT}=\Pi^{UT}=\alpha.$
Failure to satisfy the conditions does not means \textit{Result
}(i) and \textit{Result }(iii) could not be obtained. But if
$\lambda_1=0,$ these conditions are always met. Hence, under
$H_0^\star:\theta=0,$ $\alpha^{RT} > \alpha^{PTT}$ and
$\alpha^{RT}>\alpha^{UT}=\alpha$ when $\bar{c}
> 0$ and $\lambda_2>0.$
Letting $\alpha_1=\alpha_2=\alpha,$ we write
\begin{eqnarray} &&
\Pi^{UT}(\lambda_1,\lambda_2)-\Pi^{PTT}(\lambda_1,\lambda_2) = A +
B, \nonumber \end{eqnarray} where $A= [1-\Phi(\tau_{\alpha}-\gamma
\lambda_1 \sqrt{ {C^\star}^2/({C^\star}^2 + \bar{c}^2)}
/\sigma_0)]-[1-\Phi(\tau_{\alpha}-\gamma(\lambda_1+\lambda_2
\bar{c})/\sigma_0)]$ and $B = d(\tau_{\alpha_3}-\gamma\lambda_2
C^\star/\sigma_0, \tau_{\alpha}-\gamma(\lambda_1+\lambda_2
\bar{c})/\sigma_0; 0 ) \nonumber
-d(\tau_{\alpha_3}-\gamma\lambda_2 C^\star/\sigma_0,
\tau_{\alpha}-\gamma \lambda_1 \sqrt{ {C^\star}^2/({C^\star}^2 +
\bar{c}^2)} /\sigma_0;$ $-\bar{c}/\sqrt{{C^\star}^2 +
\bar{c}^2}\;).  $ For $\bar{c}>0,$ then $\lambda_1+\lambda_2
\bar{c} \ge \lambda_1 \sqrt{ {C^\star}^2/({C^\star}^2+
\bar{c}^2)}$ and $\tau_{\alpha}-\gamma \lambda_1 \sqrt{
{C^\star}^2/({C^\star}^2 + \bar{c}^2)} \geq
\tau_{\alpha}-\gamma(\lambda_1+\lambda_2 \bar{c})/\sigma_0.$ Thus,
$ A = [1-\Phi_2]-[1-\Phi_1] \leq 0$ because $\Phi_1 \leq \Phi_2$
where $\Phi_1 = \Phi(\tau_{\alpha}-\gamma(\lambda_1+\lambda_2
\bar{c})/\sigma_0)$ and $\Phi_2= \Phi(\tau_{\alpha}-\gamma
\lambda_1 \sqrt{ {C^\star}^2/({C^\star}^2 + \bar{c}^2)}).$ We
observe three cases
\begin{eqnarray} \Pi^{UT}(\lambda_1,\lambda_2)-\Pi^{PTT}(\lambda_1,\lambda_2) \gleq 0
& \mbox{  if  } & B \gleq |A|, \nonumber
\end{eqnarray} In a special case, $\lambda_1 = 0=\lambda_2,$ $A=0$
and $B>0,$ thus, $\Pi^{UT}(0,0)
> \Pi^{PTT}(0,0).$

When $\bar{c}>0$ and $\lambda_2>0,$ the asymptotic size of the RT
is larger than both UT and PTT. For $\bar{c}>0$ and $\lambda_1=0,$
the size of the PTT may also be smaller than that of UT (when
$\lambda_2$ is small). Similarly, for $\bar{c}<0,$ $\alpha^{RT} <
\alpha$ and $\alpha^{RT} < \alpha^{PTT}$ while $\alpha^{PTT}$ is
more closer to $\alpha.$

Refer to equation ({\ref{pistar}), as $\alpha_3 \rightarrow 0$ and
$\tau_{\alpha_3}-\gamma \lambda_2 C^\star/\sigma_0 \rightarrow
\infty$, $\Phi(\tau_{\alpha_3}-\gamma \lambda_2 C^\star/\sigma_0)
\rightarrow 1$ and $d(\tau_{\alpha_3}-\gamma \lambda_2
C^\star/\sigma_0, \tau_{\alpha_1}-\gamma \lambda_1
\sqrt{{C^\star}^2/({C^\star}^2+\bar{c}^2)}/ \sigma_0;
-\bar{c}/\sqrt{{C^\star}^2+\bar{c}^2}) \rightarrow 0$ because one
of the lower limits is approaching infinity. Thus, we observe that
\begin{eqnarray} \Pi^{PTT}(\lambda_1,\lambda_2)\rightarrow
1-\Phi(\tau_{\alpha_2} -\gamma (\lambda_1+ \lambda_2
\bar{c})/\sigma_0)=\Pi^{RT}(\lambda_1,\lambda_2)\;\;\mbox{as}\;\;
\alpha_3 \rightarrow 0. \label{app1} \end{eqnarray} Whereas as
$\alpha_3 \rightarrow 1$ and $\tau_{\alpha_3}-\gamma \lambda_2
C^\star/\sigma_0 \rightarrow -\infty,$
$\Phi(\tau_{\alpha_3}-\gamma \lambda_2 C^\star/\sigma_0)
\rightarrow 0$ and \begin{eqnarray} && d (\tau_{\alpha_3}-\gamma
\lambda_2 C^\star/\sigma_0, \tau_{\alpha_1}-\gamma \lambda_1
\sqrt{{C^\star}^2/({C^\star}^2+\bar{c}^2)}/ \sigma_0 \;;
-\bar{c}/\sqrt{{C^\star}^2+\bar{c}^2}) \nonumber
\\&& \hspace{100pt} \rightarrow 1-\Phi(\tau_{\alpha_1}-\gamma
\lambda_1 \sqrt{{C^\star}^2/({C^\star}^2+\bar{c}^2)}/\sigma_0)
\nonumber
\end{eqnarray} because one of the lower limits is approaching
negative infinity. Thus, we observe that
\begin{eqnarray} \Pi^{PTT}(\lambda_1,\lambda_2)\rightarrow
1-\Phi(\tau_{\alpha_1}-\gamma \lambda_1
\sqrt{{C^\star}^2/({C^\star}^2+\bar{c}^2)}/\sigma_0)=\Pi^{UT}(\lambda_1,\lambda_2)
\; \; \mbox{as}\;\; \alpha_3 \rightarrow 1. \label{app2}
\end{eqnarray}
The analytical results in this section is accompanied with an
illustrative example in investigating the comparison of the power
of the tests discussed in the next section. The power of the tests
at any value other than $\theta=0$ is also considered in the
example to study the behavior of the power functions corresponds
to the probabilities of type I and type II errors. Moreover, the
study of relationship between the level of significance for the
PTT and the nominal size of the PT as well as the nominal sizes of
the UT and RT are explored.

\section{Illustrative Example - Power Comparison}
\setcounter{equation}{0} The asymptotic power functions for the
 UT, RT and PTT are compared in this section. Under $\{K_n\},$ we
note
\begin{description}
\item{(i)} $\Pi^{UT}(\lambda_1,\lambda_2)$ as the asymptotic power
function for testing $H_0^\star:\theta=0$ when $\beta$ is assumed
to be undefined in the construction of the test statistic
$T_n^{UT}$,

\item{(ii)} $\Pi^{RT}(\lambda_1,\lambda_2)$ as the asymptotic
power function for testing $H_0^\star:\theta=0$ when $\beta$ is
assumed to be zero in the construction of the test statistic
$T_n^{RT}$ and

\item{(iii)}
 $\Pi^{PTT}(\lambda_1,\lambda_2)$ as the asymptotic
power function  for testing $H_0^\star:\theta=0$ after pre-testing
$H_0^{(1)}:\beta=0$.
\end{description}

For this illustrative example, the random errors of the simple
linear model are generated from Normal distribution with mean $0$
and variance $1$. The sample size is $n=1000.$ Three sets of
values: 0 and 1 with 50\% for each for the first set, $-1$ and 0
with 50\% for each for the second set and $-1$ and 1 with 50\% for
each for the third set are considered as the values of the
regressor $c_i,\; i=1,2,\ldots,1000.$ These values guarantee
$\bar{c}>0$, $\bar{c}=0$ and $\bar{c}<0$ respectively to the sets
of regressors.

In this example, the $\psi$ function is taken as Huber $\psi$
function (Hoaglin et al., 1983, p.366, Wilcox, 2005, p.77), is
defined as
\[ \psi_h(u_i) =\left
\{ \begin{array} {ccc} -k & \mbox{if }&u_i <-k \\
u_i & \mbox{if }& |u_i| \leq k \\
k & \mbox{if }& u_i > k, \\ \end{array} \right. \] where
$u_i=X_i-\theta - \beta c_i.$ As suggested in many reference books
(Wilcox, 2005, p.76), the value of $k=1.28$ is chosen because
$k=1.28$ is the 0.9 quantile of a standard normal distribution,
there is a 0.8 probability that a randomly sampled observations
will have a value between $-k$ and $k$ (Wilcox, 2005, p.76). The
estimate for $\sigma_0$ is taken to be $\sum \psi(u)^2/n.$ The
estimate for $\gamma$  is $\sum \psi^\prime(u)/n$ (Caroll and
Rupert, 1988, p.212) where
\[\psi^\prime(u)=\left
\{ \begin{array} {ccc} 0 & \mbox{if }&u <-1.28 \\
1 & \mbox{if } &|u| \leq 1.28 \\
0 & \mbox{if }& u > 1.28. \\ \end{array} \right. \] The
$\Pi^{PTT}$, $\Pi^{RT}$ and $\Pi^{UT}$ are calculated using the
formulas given by equations (\ref{pistar}), (\ref{pi1}) and
(\ref{pi2}). The R-package (mvtnorm) is used in computing the
bivariate Normal probability integral.

\begin{figure}
\begin{center}
\includegraphics[width=0.35\textwidth]{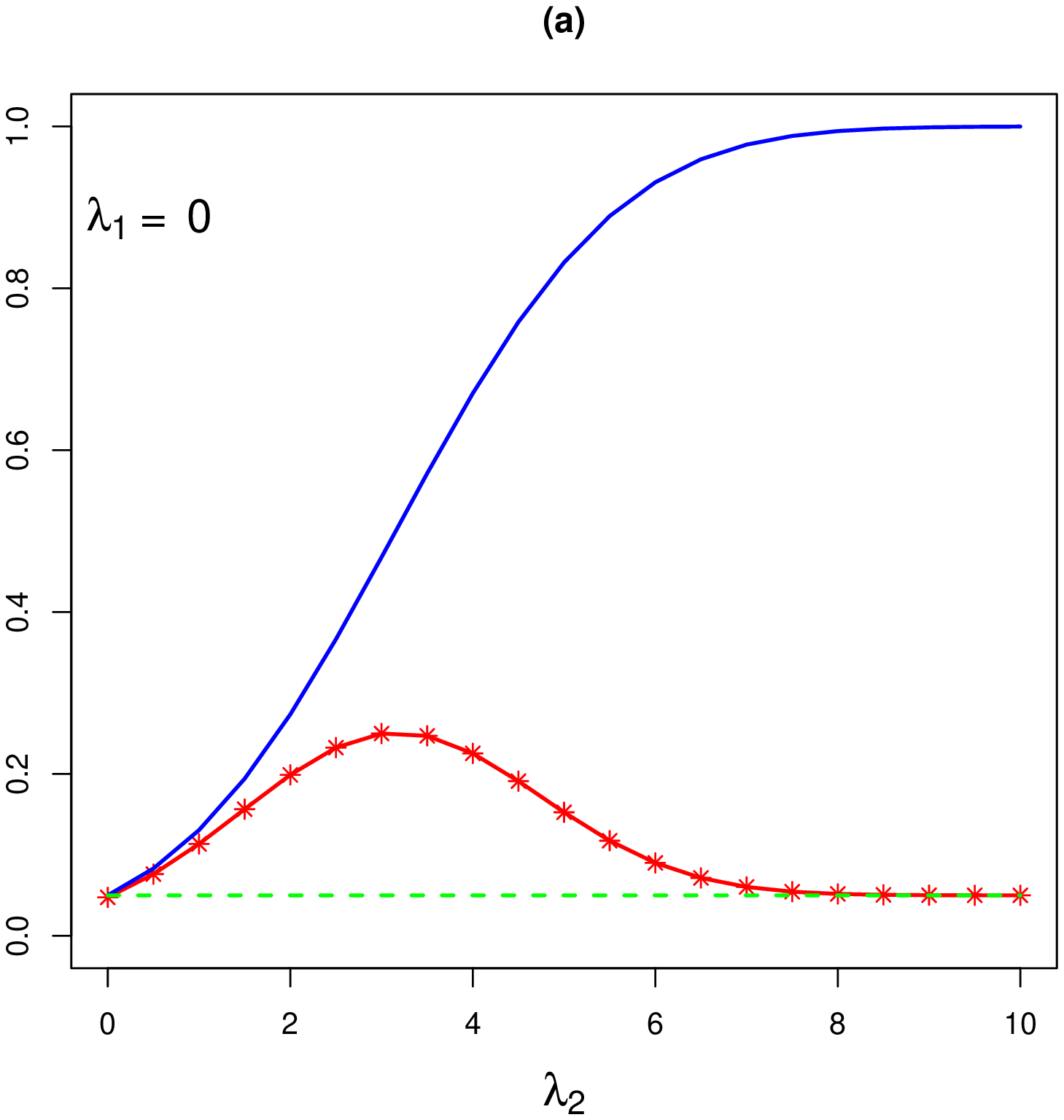} 
\includegraphics[width=0.35 \textwidth]{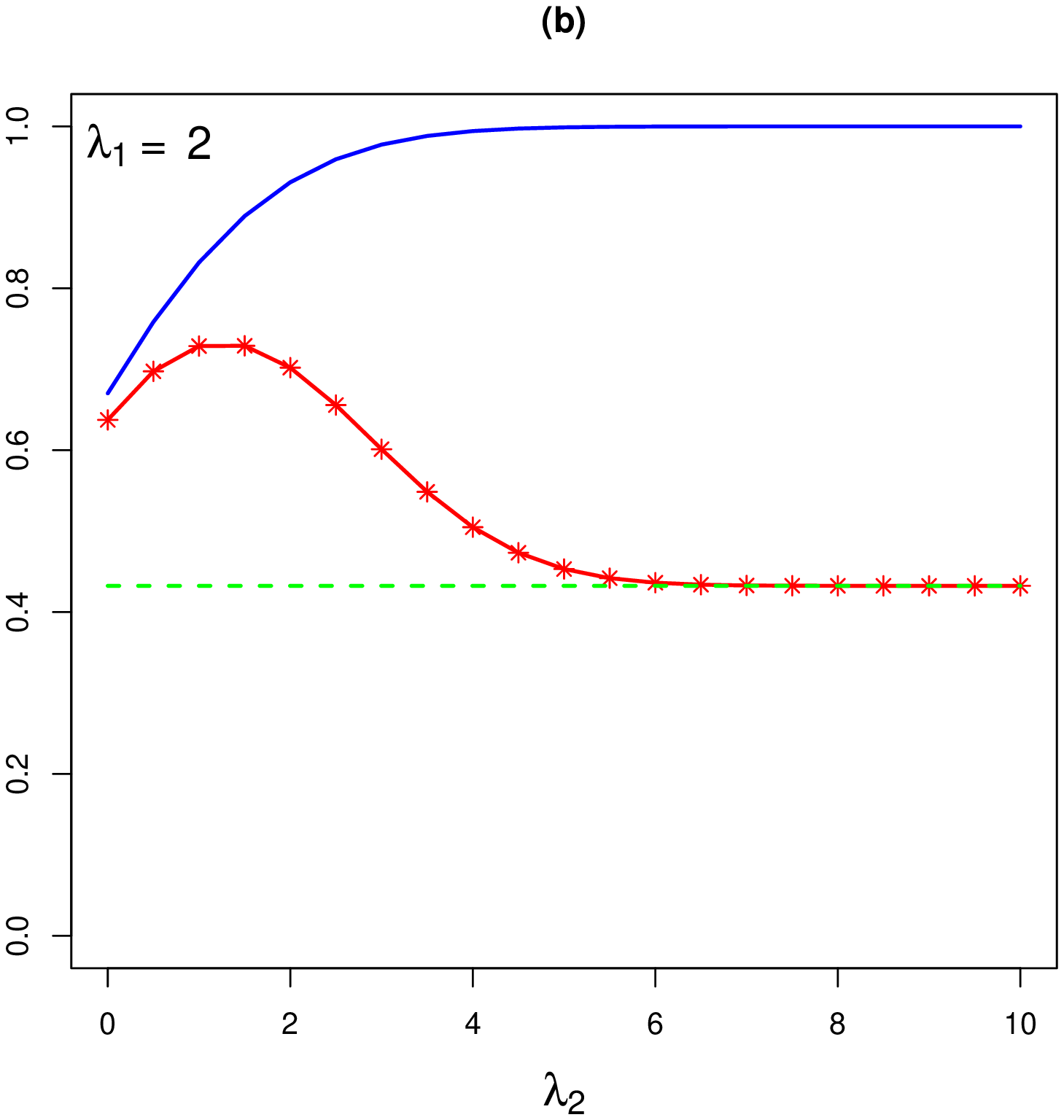} 
\includegraphics[width=0.35\textwidth]{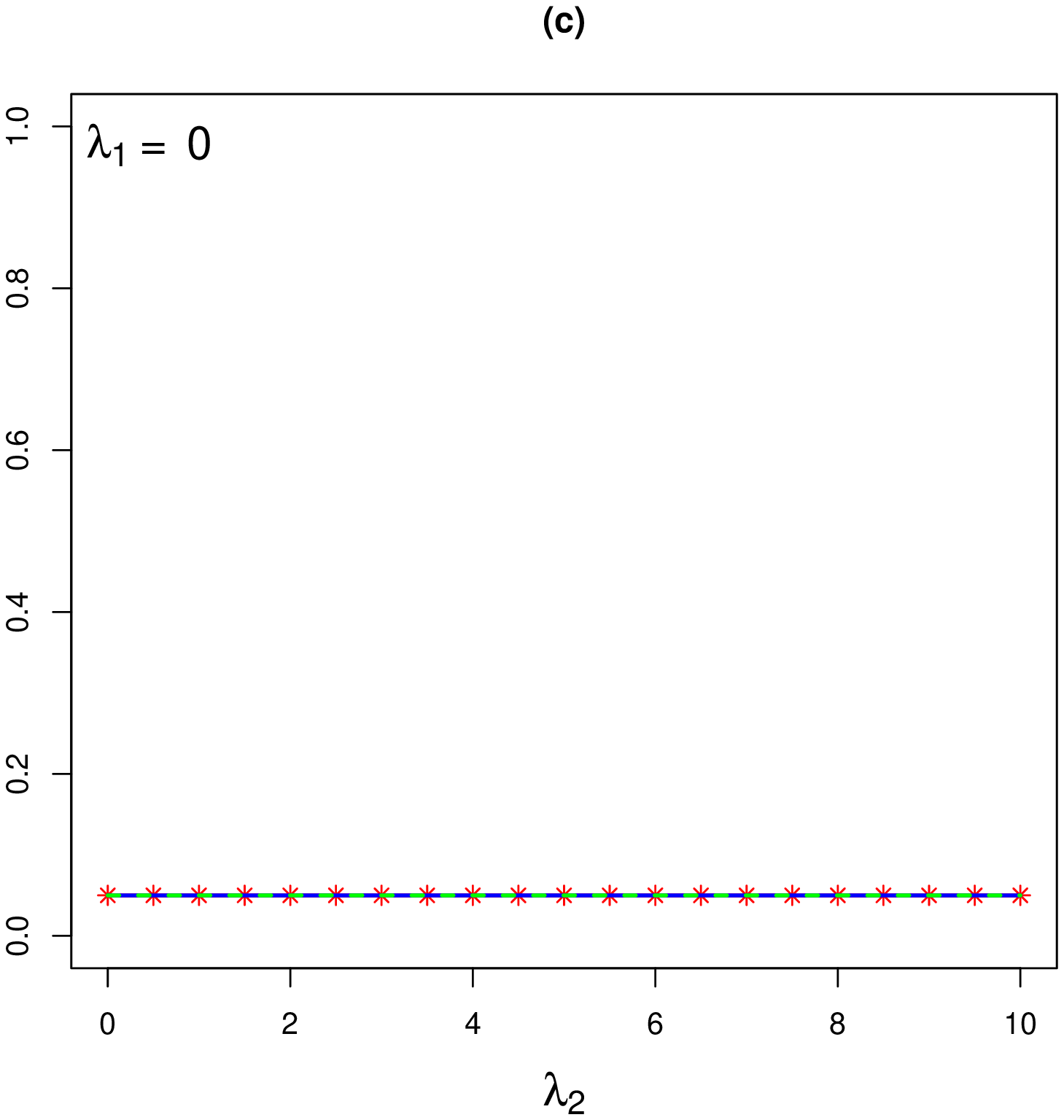}    `
\includegraphics[width=0.35 \textwidth]{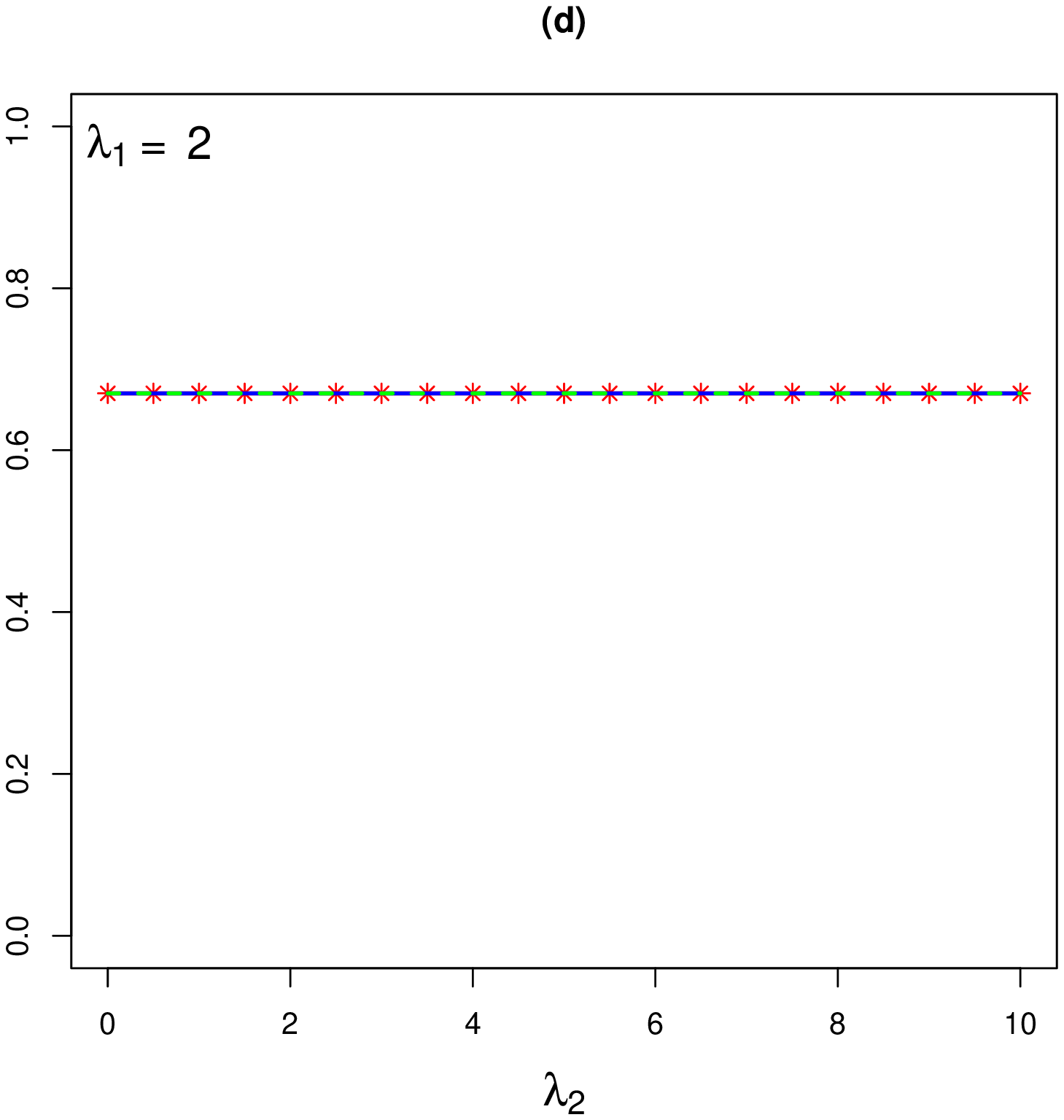}
\includegraphics[width=0.35 \textwidth]{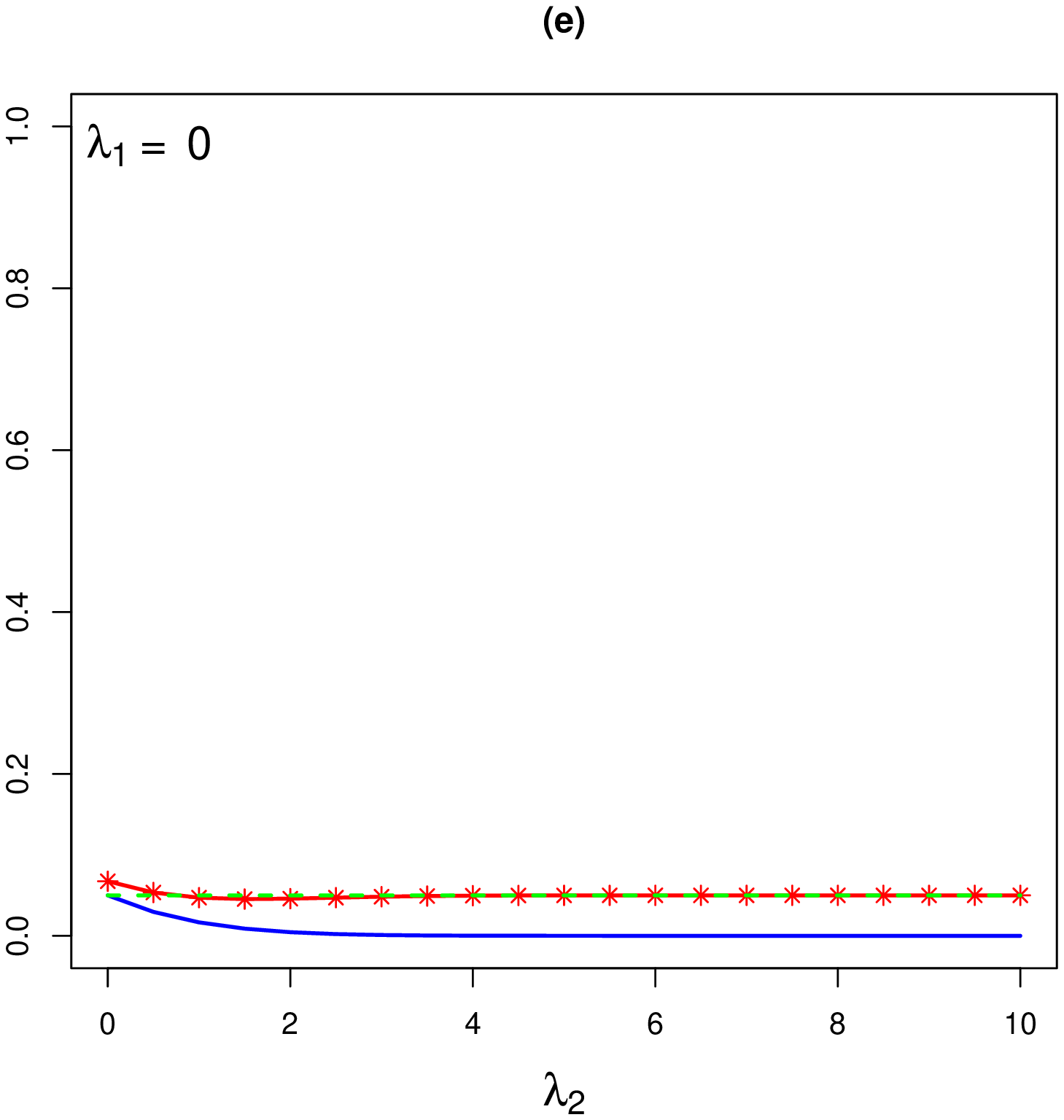}
\includegraphics[width=0.35 \textwidth]{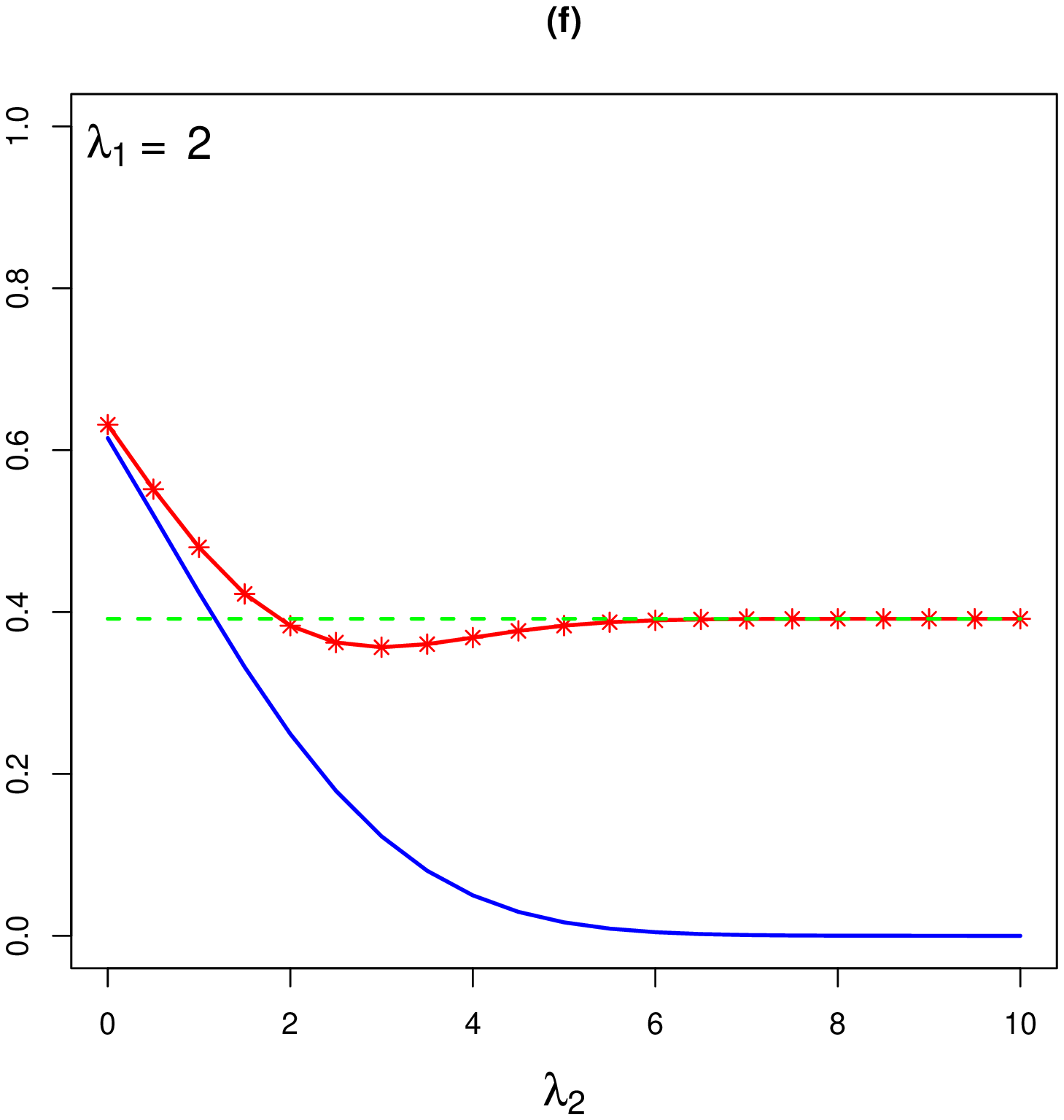}
\caption{Graphs of power functions as a function of \(\lambda_2\)
for selected values of \(\lambda_1\) and
\(\alpha_1=\alpha_2=\alpha_3=\alpha=0.05\). Dotted line, solid
line and line with star represent \(\Pi^{UT}
(\lambda_1,\lambda_2)\), \(\Pi^{RT}(\lambda_1,\lambda_2) \) and
\(\Pi^{PTT} (\lambda_1,\lambda_2)\) respectively. Graphs (a) and
(b) are for \(\bar{c}>0 \),
(c) and (d) are for \(\bar{c}=0 \) and (e) and (f) are for \(\bar{c}<0\). } 
\label{fig1}
\end{center}
\end{figure}

In Figure 1, the power functions for the UT, RT and PTT are
plotted against $\lambda_2$ at two values of $\lambda_1.$ Here
$\lambda_1=0$ is chosen to study the asymptotic sizes of the tests
and we desire the size of a particular test to be small so that
the probability of type I error is small. Since we also wish to
get small value of probability of type II error, the power of the
test at $\lambda_1=2$ is considered. An acceptable power function
of the test is the one that is small when the null hypothesis is
true but large when $\lambda_1$ differs much from $\theta=0.$ The
first set of regressors is used to plot Figures 1(a) and 1(b). As
$\lambda_2$ grows larger, $\Pi^{RT}(0,\lambda_2)$ approaches 1.
However, $\Pi^{PTT}(0,\lambda_2)$, after an initial increase,
drops and converges to the nominal size $\alpha=0.05$ as
$\lambda_2$ grows larger. Thus, the asymptotic size (with very
small $\lambda_1$) of $\phi_n^{PTT}$ is close to $\alpha$ for
small $\lambda_2$ and large $\lambda_2$, while for moderate values
of $\lambda_2$ it is somehow larger than $\alpha$ but lesser than
that of $\Pi^{RT}(0,\lambda_2)$. The $\Pi^{UT}(0,\lambda_2)$ is
constant and does not depend on $\lambda_2.$ The same pattern
occurs in Figure 1(b) but the power functions are always
significantly larger than $\alpha$, in this case larger than 0.4.
If one only considers the size of the test, the PTT is preferred
to RT, though the UT remains as the best choice. However, the RT
is the best choice but the PTT is preferred to UT if the power of
the test at $\lambda_1=2$ is considered. It is impossible to
obtain a test that uniformly minimizes the size and maximizes the
power at the same time. We are looking for a test that is a
compromise between minimizing the size and maximizing the power
(small probabilities of type I and type II errors). The RT is the
best choice for its largest power but the worst choice for its
largest size as $\lambda_2$ grows larger. On the contrary, the UT
is the best choice for its smallest size but the worst choice for
its smallest power. Both RT and UT uniformly minimize or maximize
the size and power at the same time. The PTT has larger power than
the UT for small and moderate values of $\lambda_2$ and it has
significantly smaller size than that of the RT for moderate and
large $\lambda_2$. Therefore, if our objective is to obtain a test
that has better probabilities for both type I and type II errors,
the PTT is suggested as the best option. The PTT is a compromise
between minimizing the size and maximizing the power than the RT
and UT.

The cases for $\bar{c}=0$ and $\bar{c}<0$ are also considered in
this paper, though $\bar{c}>0$ is more important than the other
two because it is more realistic. Setting $\bar{c}=0$ in Figures
1(c) and 1(d) imply all power functions remain the same regardless
of the value of $\lambda_2$ and these constant power functions
increase as $\lambda_1$ increases. Figures 1(e) and 1(f)
illustrate the case when $\bar{c} <0.$ The graphs show that
$\Pi^{RT}<\Pi^{PTT}$ for any $\lambda_2$ and $\Pi^{PTT} \leq
\Pi^{UT}$ for any $\lambda_2$ more than a small positive value,
say $\lambda_0$. The probability of type I error for all test
functions are fairly small. The size and power of the RT is
decreasing to 0 as $\lambda_2$ growing larger (Figures 1(e) and
1(f)) suggesting the RT as the best choice for size but the worst
choice for power. Since $\Pi^{PTT}(2,\lambda_2) \geq
\Pi^{RT}(2,\lambda_2)$ for all $\lambda_2,$ the PTT is preferred
over the RT . Also, $\Pi^{PTT}(2,\lambda_2)\geq
\Pi^{UT}(2,\lambda_2)$ except for some moderate values of
$\lambda_2$ but the difference is relatively small. From the
examination of all the graphs in Figure 1, the PTT is suggested as
the best choice when both probabilities of type I and type II
errors are considered.

\begin{figure}
\begin{center}
\includegraphics[width=0.35\textwidth]{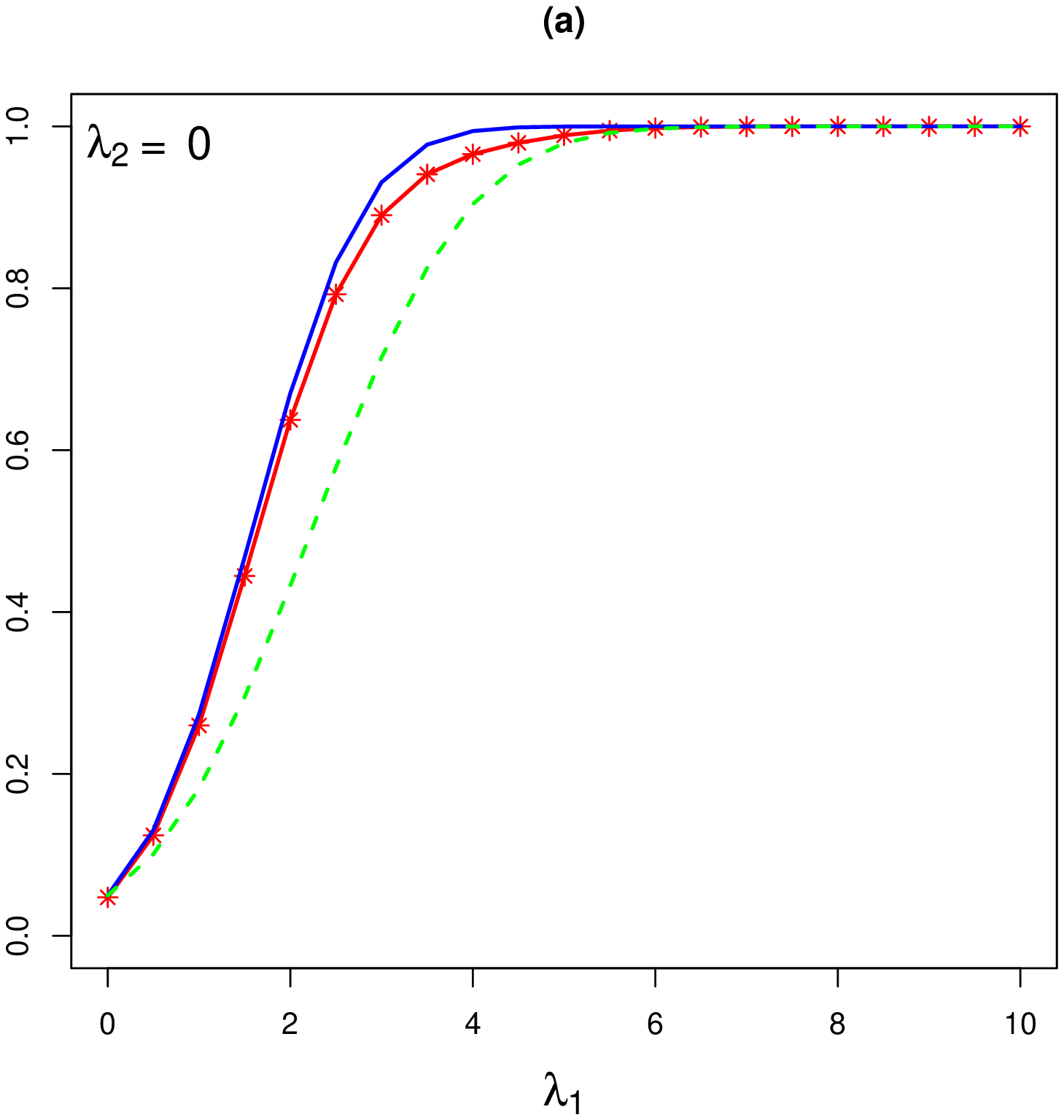} 
\includegraphics[width=0.35 \textwidth]{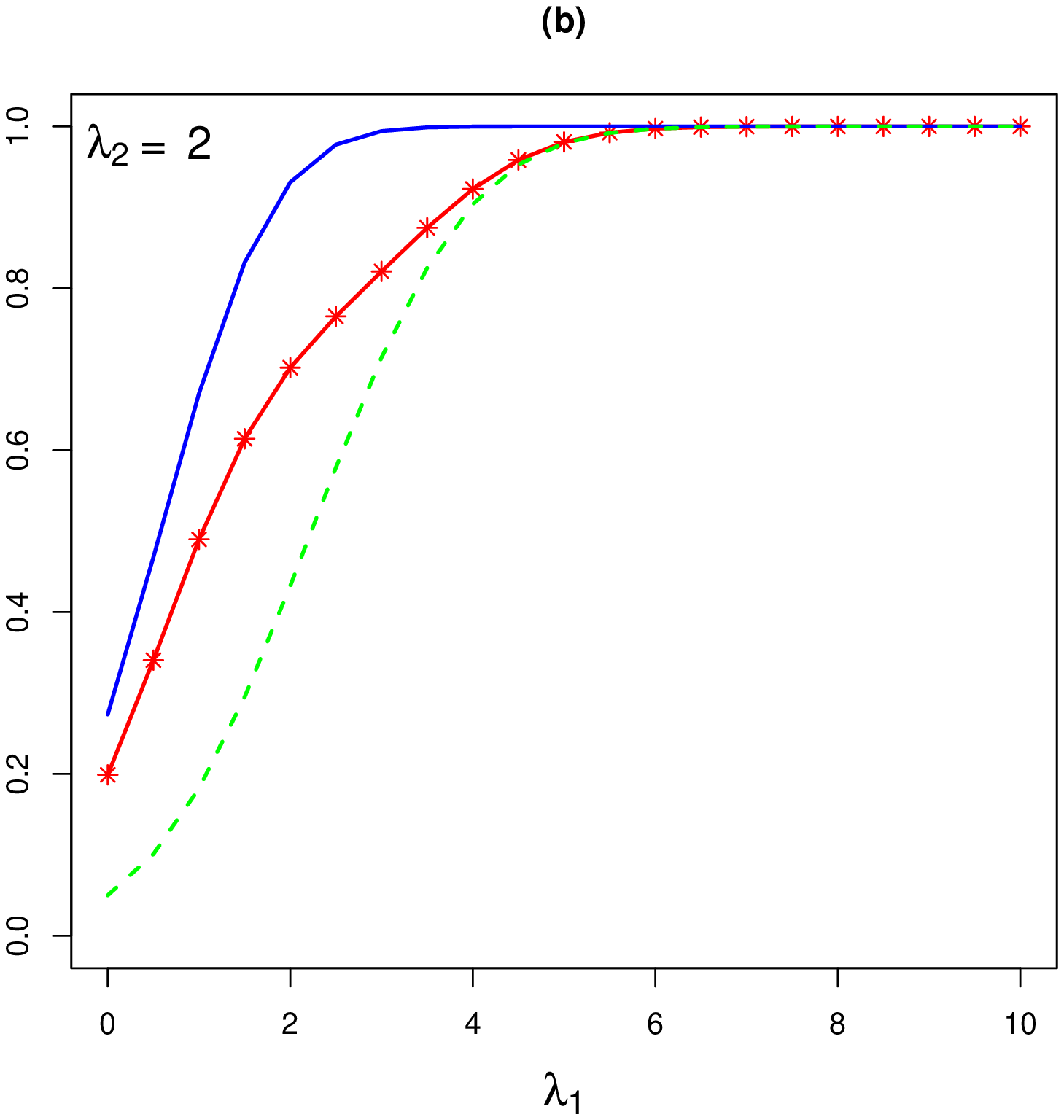} 
\includegraphics[width=0.35\textwidth]{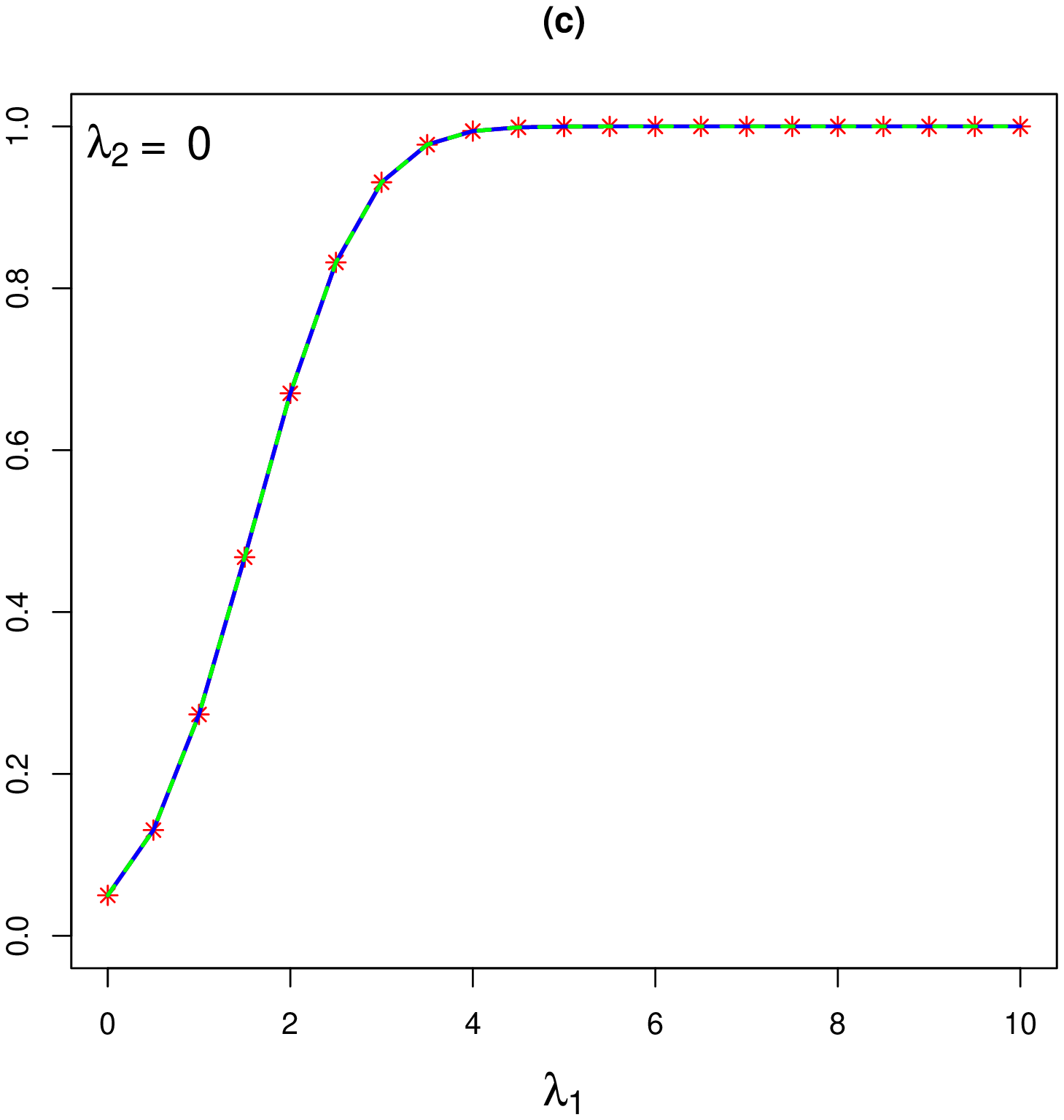} 
\includegraphics[width=0.35 \textwidth]{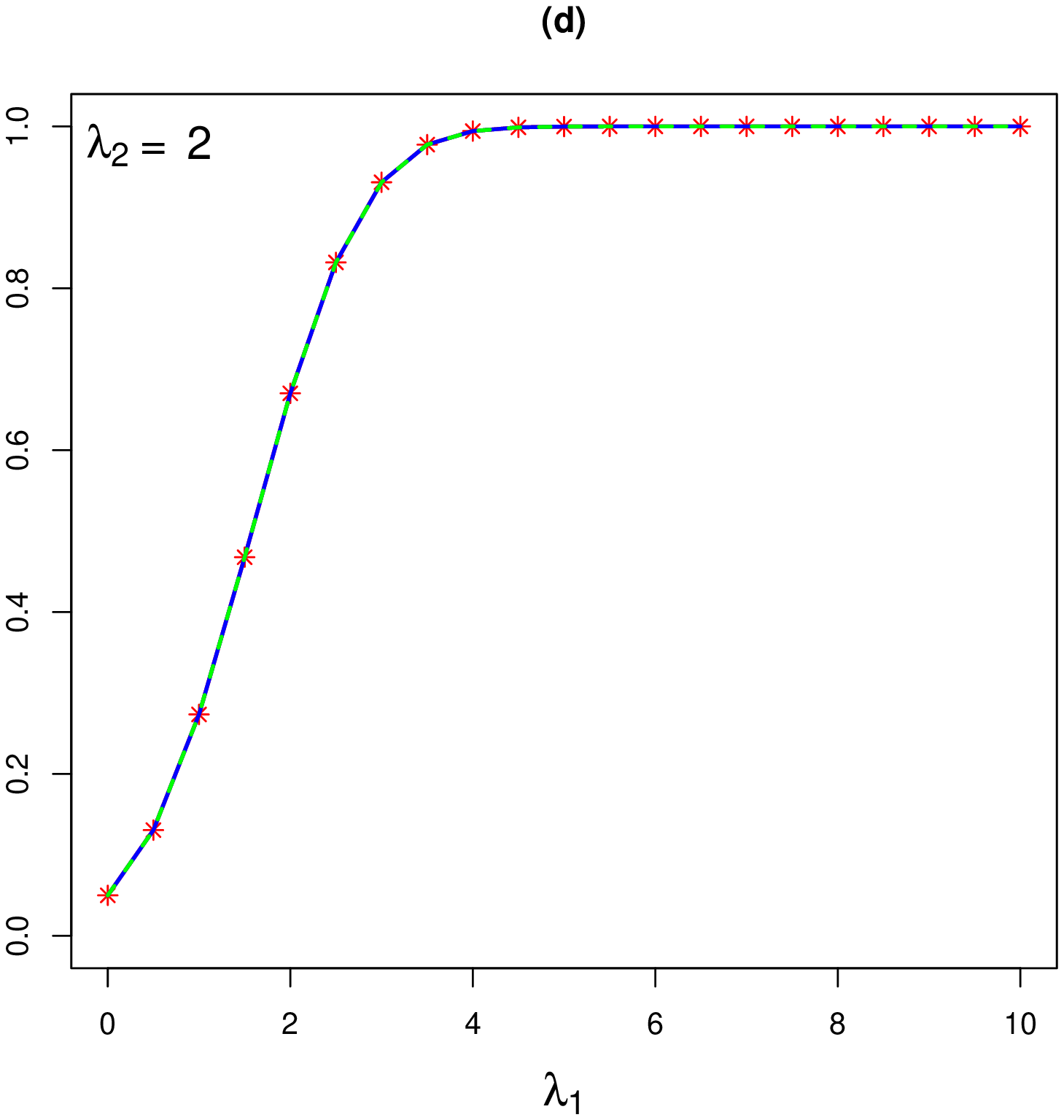}
\includegraphics[width=0.35 \textwidth]{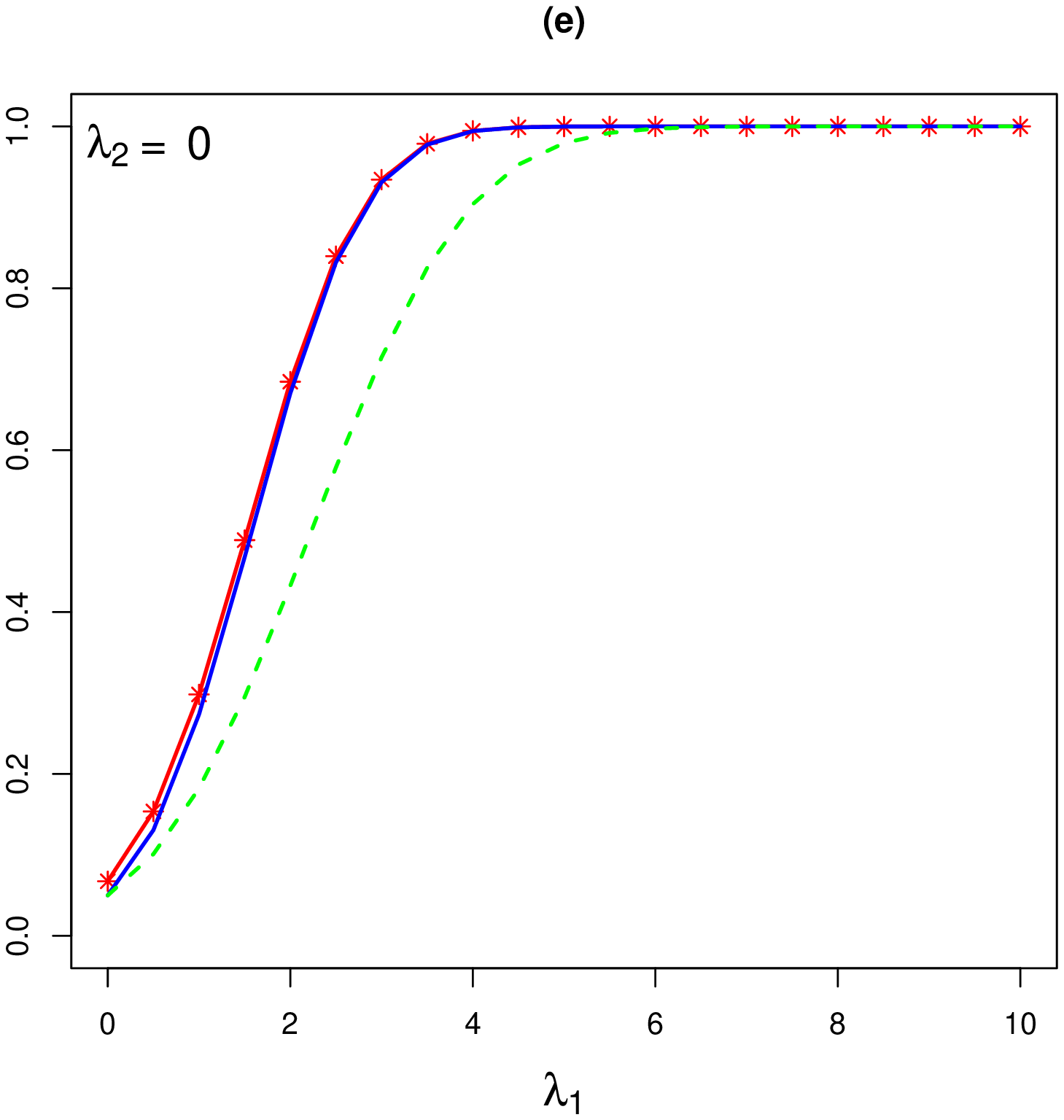}
\includegraphics[width=0.35 \textwidth]{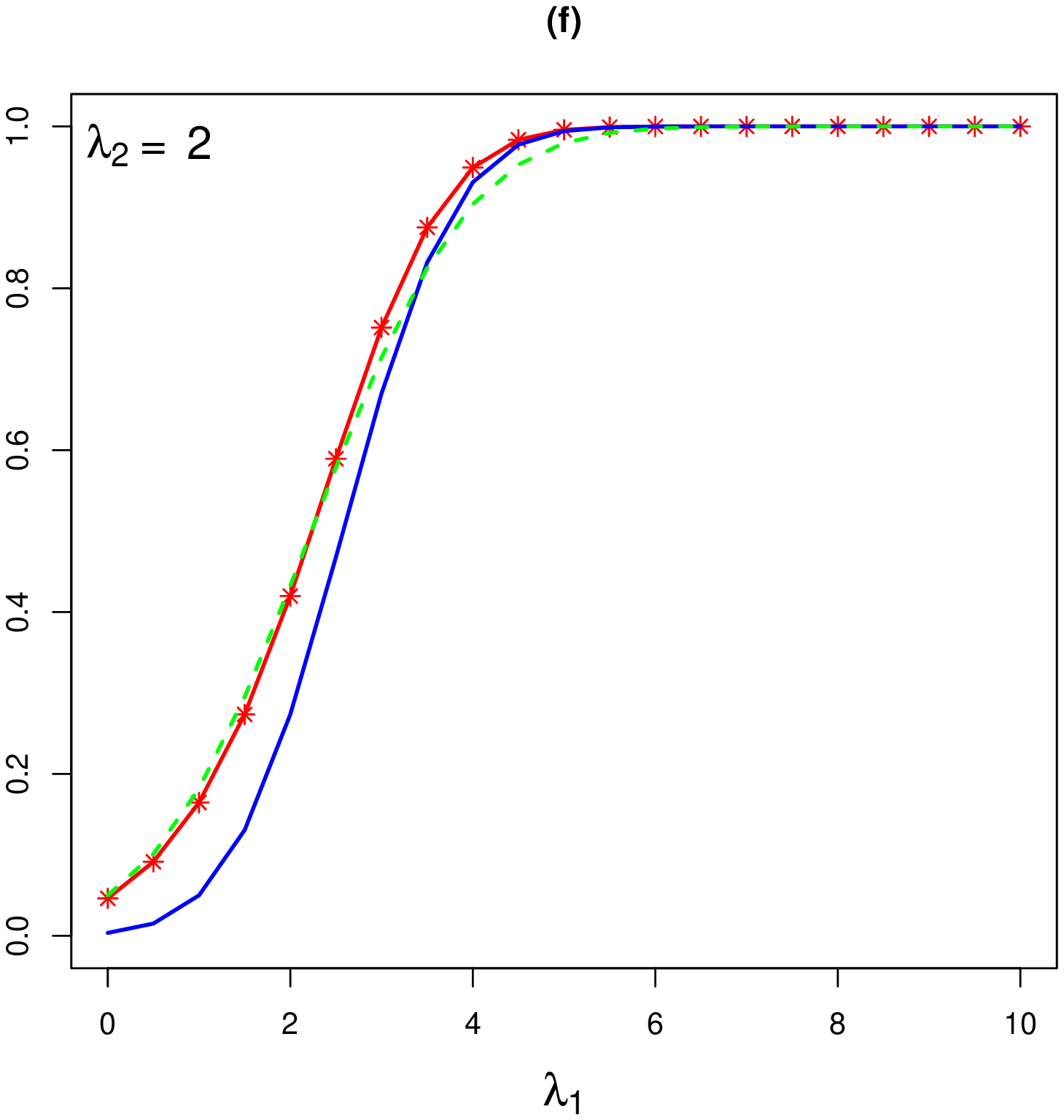}
\caption{Graphs of power functions as a function of \(\lambda_1\)
for selected values of \(\lambda_2\) and
\(\alpha_1=\alpha_2=\alpha_3=\alpha=0.05\). Dotted line, solid
line and line with stars represent
\(\Pi^{UT}(\lambda_1,\lambda_2)\),
\(\Pi^{RT}(\lambda_1,\lambda_2)\) and \(\Pi^{PTT}
(\lambda_1,\lambda_2)\) respectively. Graphs (a) and (b) are for
\(\bar{c}>0 \),
 (c) and (d) are for \(\bar{c}=0 \) and (e) and (f) are for \(\bar{c}<0\).} 
\label{fig2}
\end{center}
\end{figure}

\begin{figure}
\begin{center}
\includegraphics[width=0.35\textwidth]{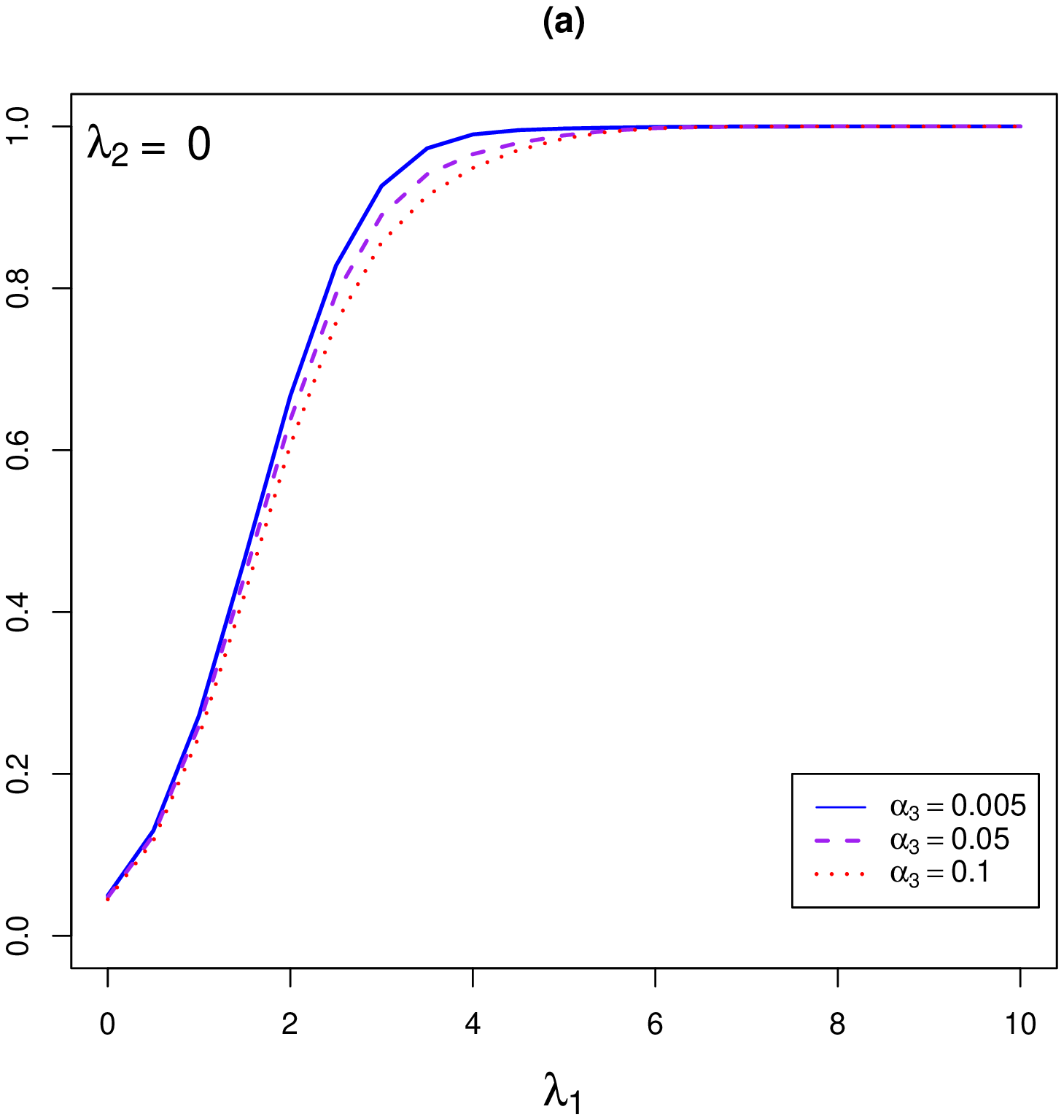} 
\includegraphics[width=0.35 \textwidth]{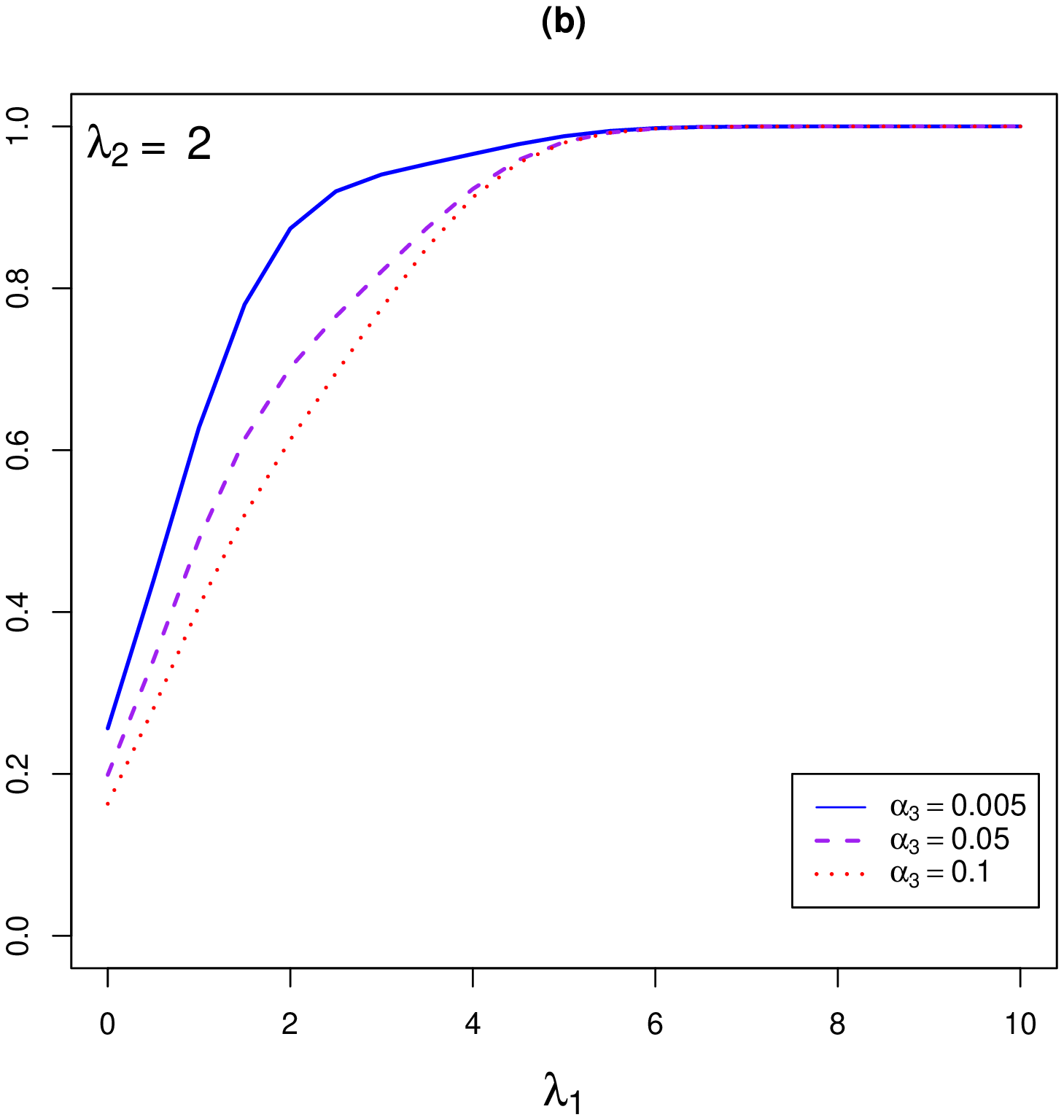} 
\caption{Graphs of power function \(\Pi^{PTT}(\lambda_1,\lambda_2)
\) for nominal sizes \(\alpha_3=0.005\), 0.05 and 0.1. The \(\bar{c} > 0 \) and \(\alpha_2=\alpha_1=\alpha=0.05\) for all graphs.} 
\label{fig3}
\end{center}
\end{figure}

The relation between power functions and $\lambda_1$ is shown in
Fig 2. All power functions are approaching 1 as $\lambda_1$ grows
larger regardless of the value of $\lambda_2.$ This is because the
probability of rejecting $H_0^\star:\theta=0$ increases as
$\lambda_1$ increases. When $\bar{c}>0$, the probability of type
II error for the RT is the smallest, but the PTT is preferable
than the UT for all values of $\lambda_1$. When $\bar{c}<0$, the
PTT is preferable for its comparatively smaller probability of
type II error than the other two tests. When $\bar{c}=0,$ all
tests have the same probability of type II error regardless of the
value of $\lambda_1$ (refer to the equation (\ref{same}) for
analytical result).

Figure 3 illustrates the behavior of the power function
$\Pi^{PTT}(\lambda_1,\lambda_2)$ at three different values of
nominal size $\alpha_3$. The graphs show that the test with
smaller nominal significance level has greater power than that of
larger significance level. The smaller nominal significance level
however increases the probability of type I error as $\lambda_2$
moves away from zero. This is illustrated in Fig 3(b),
$\Pi^{PTT}(\lambda_1,2)$ at $\alpha_3=0.005,\;\;
0.05,\;\mbox{and},\;\; 0.1$ start at different values before
growing larger and converging to 1.

It is of advantage to study the relationship between the size of
the PTT, that is, $\alpha^{PTT}=\Pi^{PTT}(0,\lambda_2)$ and the
nominal significance level of the PT, $\alpha_3.$ One may want to
know what suppose to be the actual level of significance of the PT
that will reject the ultimate test with a predetermined
probability, says 5 percent. Taking $\alpha_1$ and $\alpha_2$ to
be equal, here 0.05, the size of the test depends on $\lambda_2.$
Figure 4 shows the graphs of $\Pi^{PTT}(0,\lambda_2)$ against
$\alpha_3$ for different values of $\lambda_2$ with
$\alpha_1=\alpha_2=0.05$ and $\bar{c}>0.$ For smaller values of
$\lambda_2$, as $\alpha_3$ increases, the size of the PTT
decreases and reaches its minimum at the value of $\alpha_3=
\alpha_3^\prime$ (say), before growing larger and converging to
$\alpha=0.05.$ Let the value of $\alpha_3$ be $\alpha_3^{\prime
\prime}$ when the size of the PTT is 0.05, the value of
$\alpha_3^{\prime \prime}$ increases as $\lambda_2$ increases. As
we consider larger values of $\lambda_2$, the size of the PTT
decreases dramatically then slowly converges (appears as flat in
the graph) to $\alpha$ at some positive value $\alpha_3^{\prime
\prime \prime}$. Table 1 gives the values of the size of the PTT
at $\alpha_3$ for different values of $\lambda_2,$ with
$\alpha_1=\alpha_2=0.05$ when $\bar{c}>0.$ If we want to reject
the ultimate test with significance level 0.05, the nominal
significance level of the PT must be set to 0 when $\lambda_2$ is
0. Then a larger but still acceptable nominal size $\alpha_3$ is
required to achieve 5\% significance level of the PTT as
$\lambda_2$ is a bit larger than 0. But up to some point, we
cannot sacrifice the increases in the probability of type I error
of the PT as $\lambda_2$ grows much than 0. As $\lambda_2$ grows
larger, a larger $\alpha_3$ is required to obtain 5\% significance
level of the PTT (see Table 1). Note: Setting the nominal size of
the PT to 0 is meaningless because this means there is no chance
that $H_0^{(1)}:\beta=0$ is rejected. The size and power of the
PTT converges to the RT when $\alpha_3$ approaches 0 (refer
equation (\ref{app1})), thus supports this result in Table 1.

\begin{figure}
\begin{center}
\includegraphics[width=0.3\textwidth]{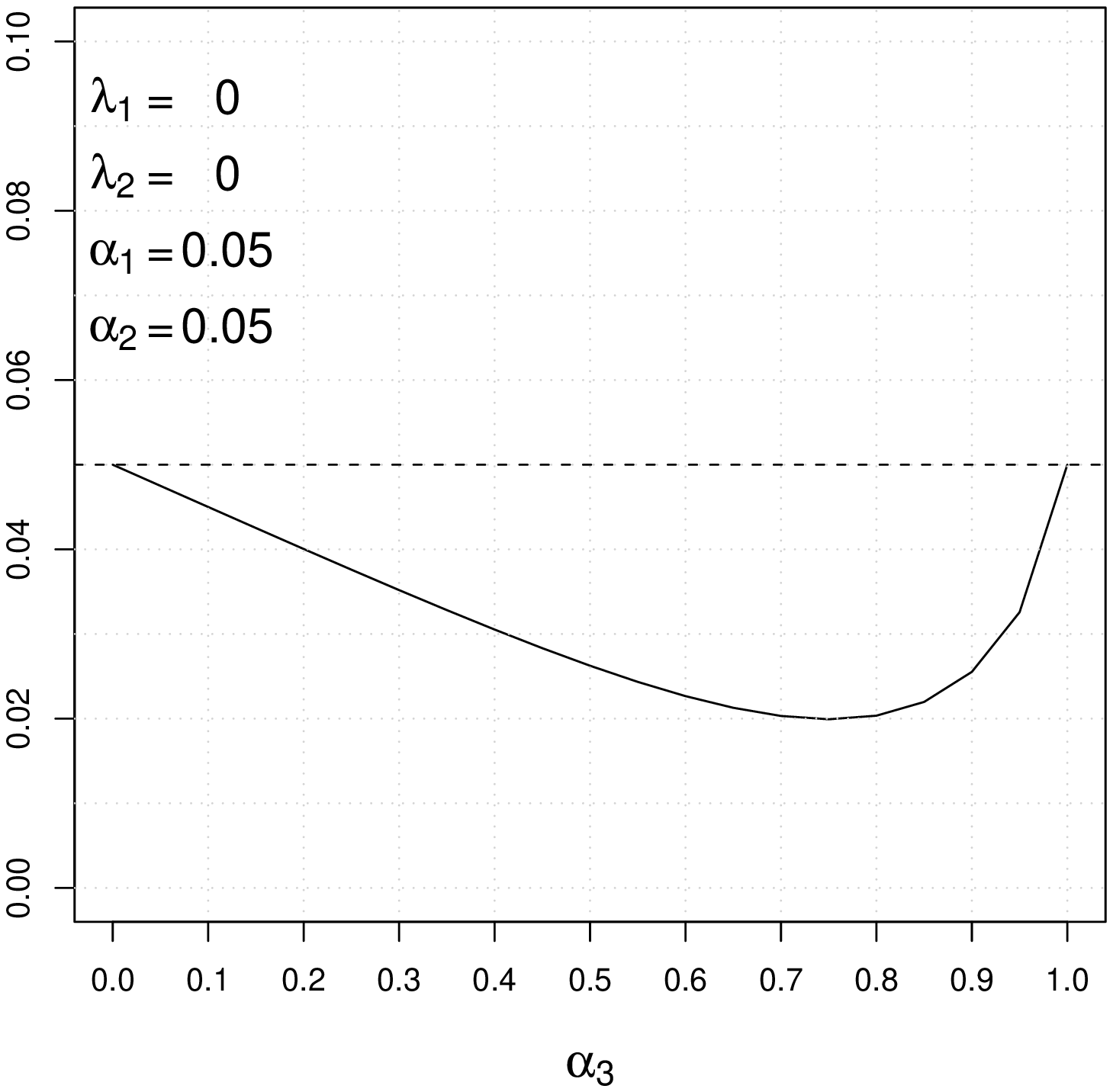} 
\includegraphics[width=0.3 \textwidth]{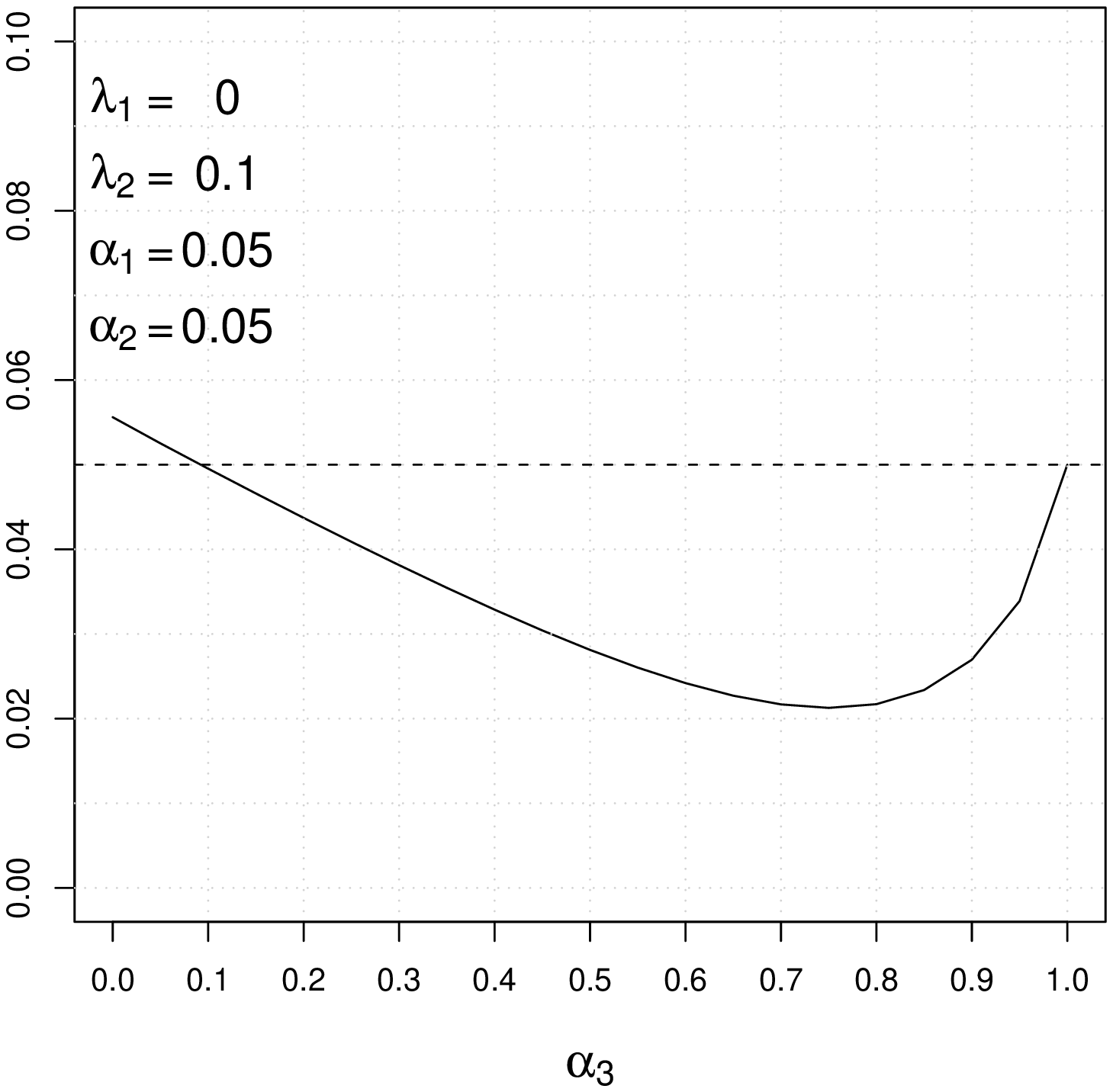} 
\includegraphics[width=0.3 \textwidth]{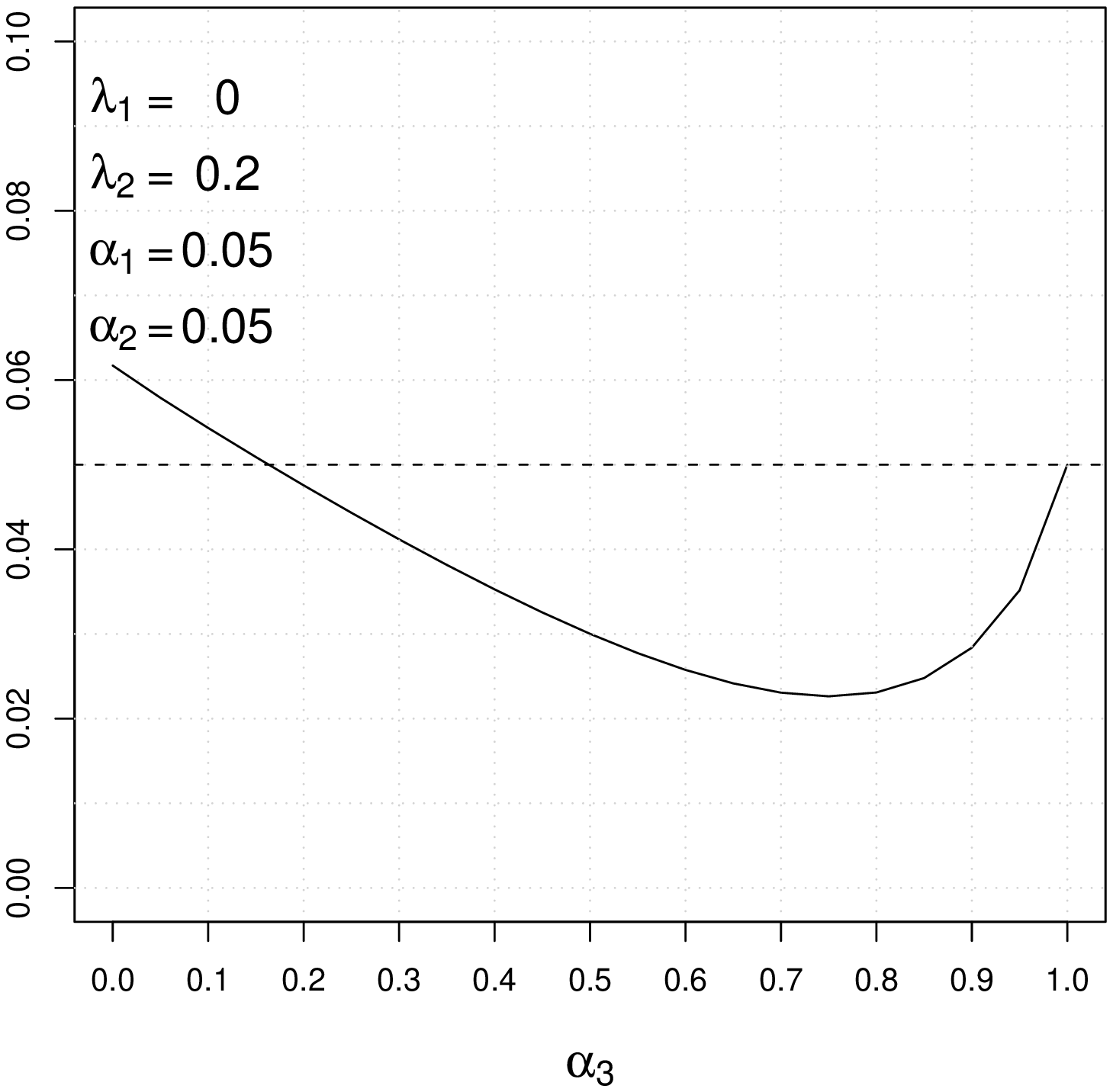} 
\includegraphics[width=0.3 \textwidth]{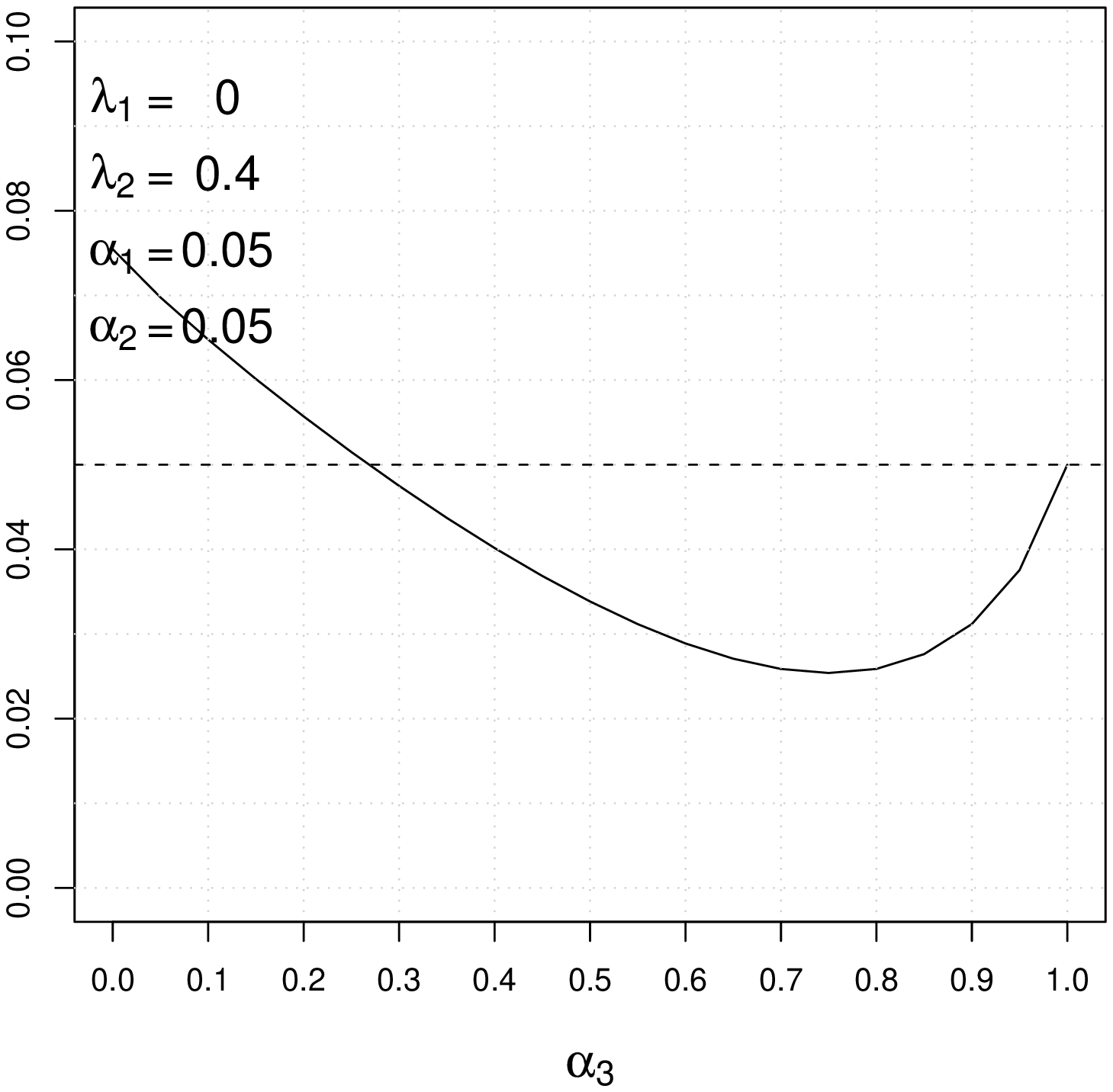} 
\includegraphics[width=0.3 \textwidth]{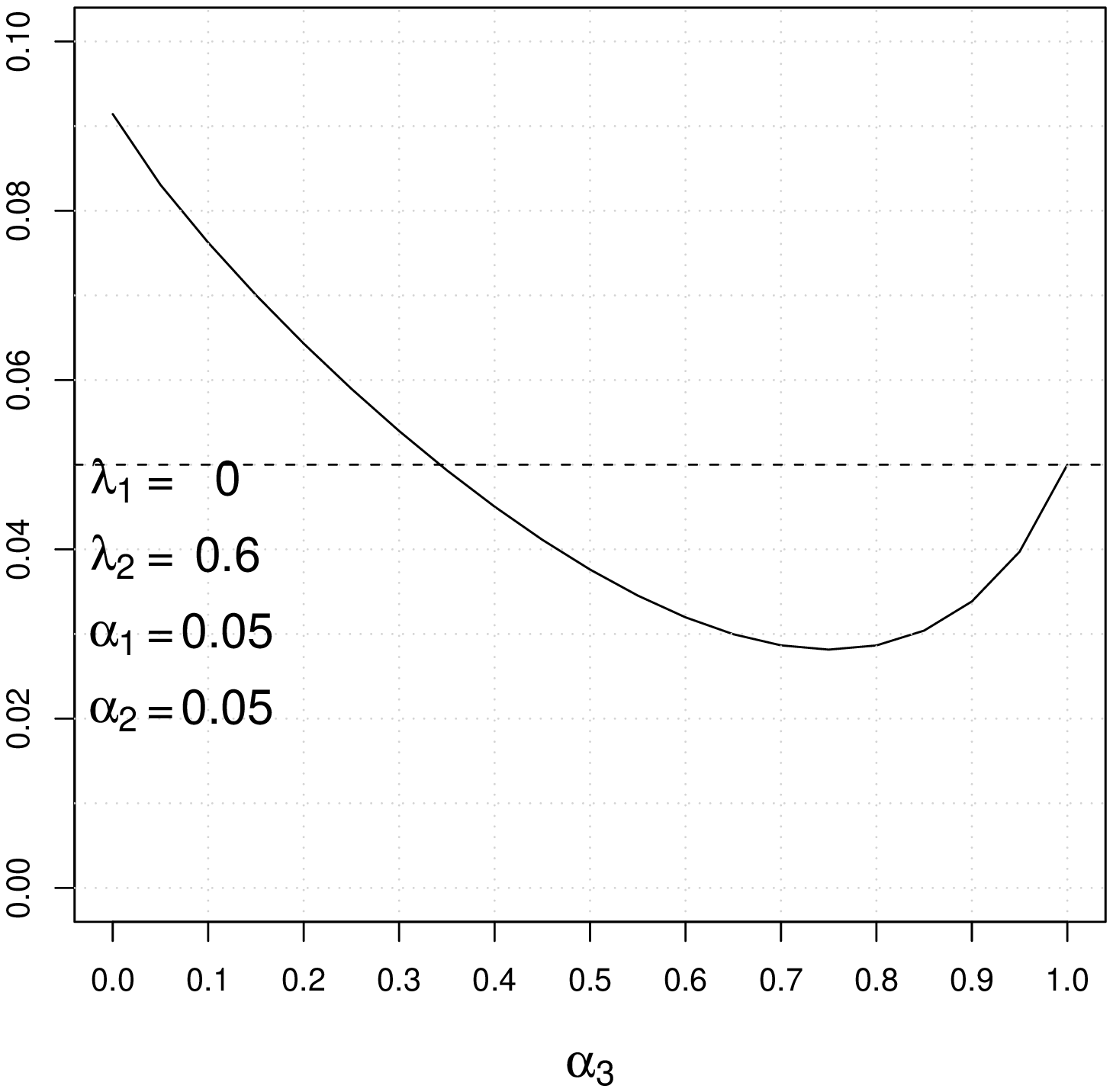}
\includegraphics[width=0.3 \textwidth]{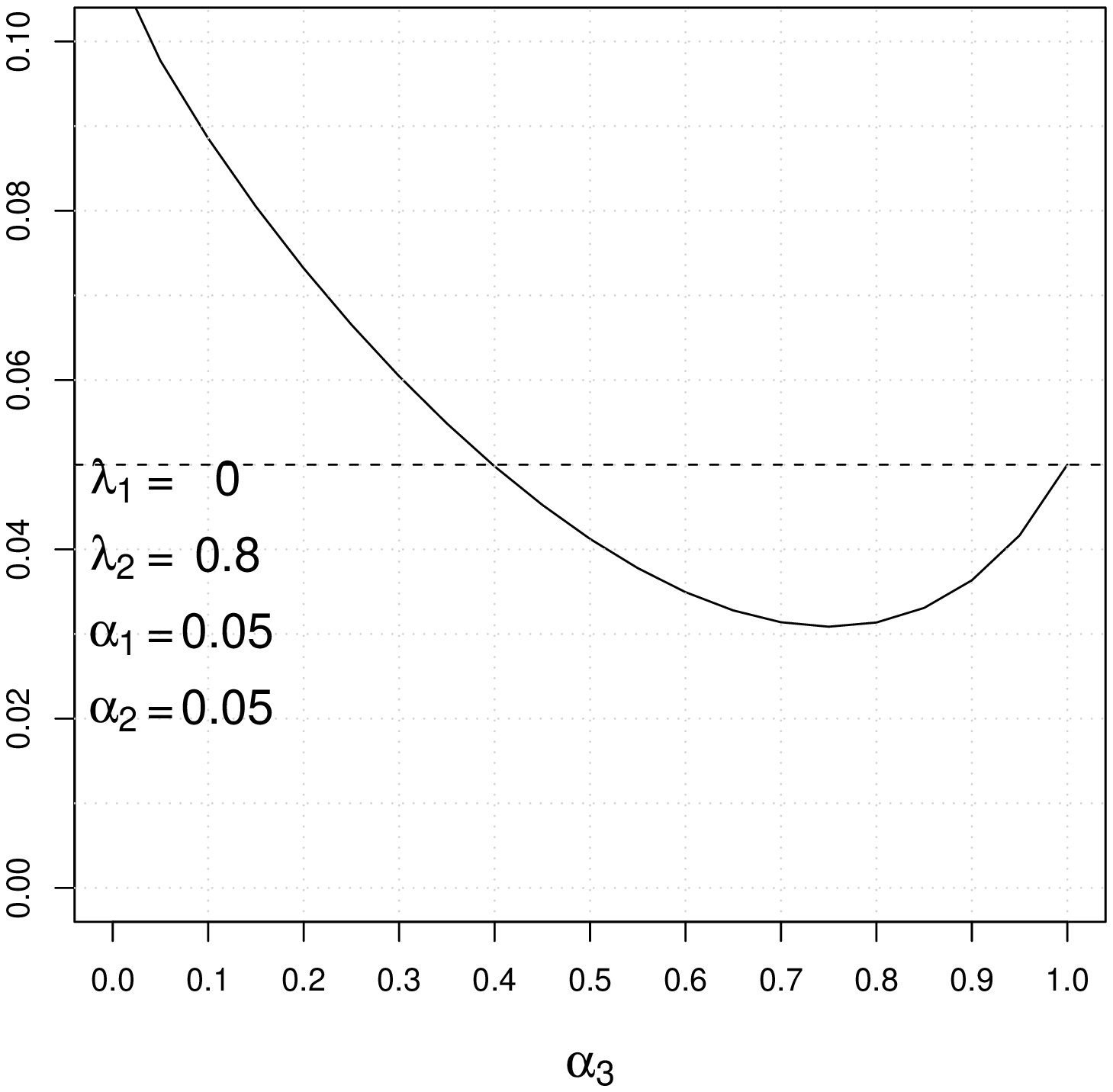} 
\includegraphics[width=0.3 \textwidth]{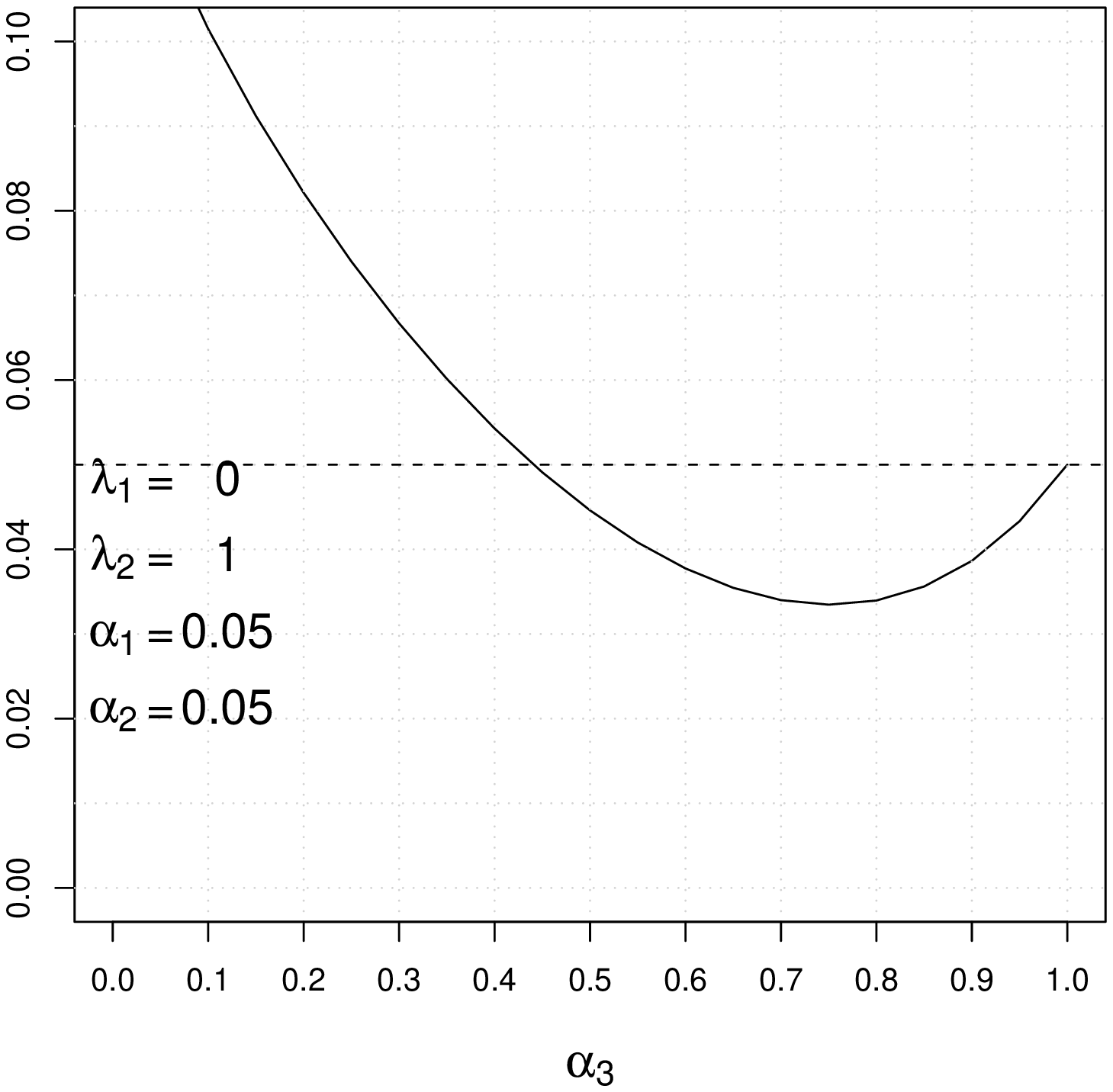} 
\includegraphics[width=0.3 \textwidth]{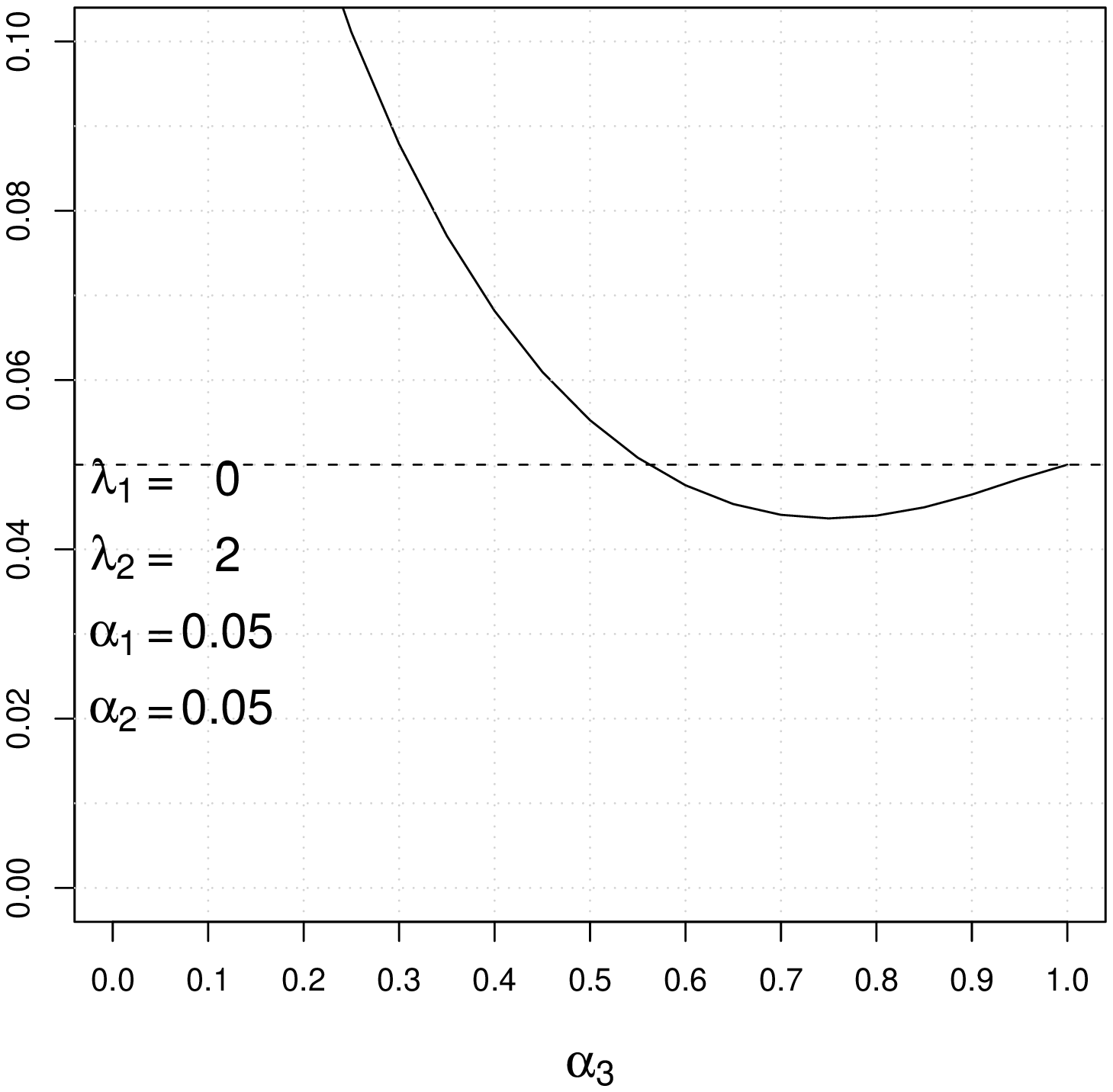} 
\includegraphics[width=0.3 \textwidth]{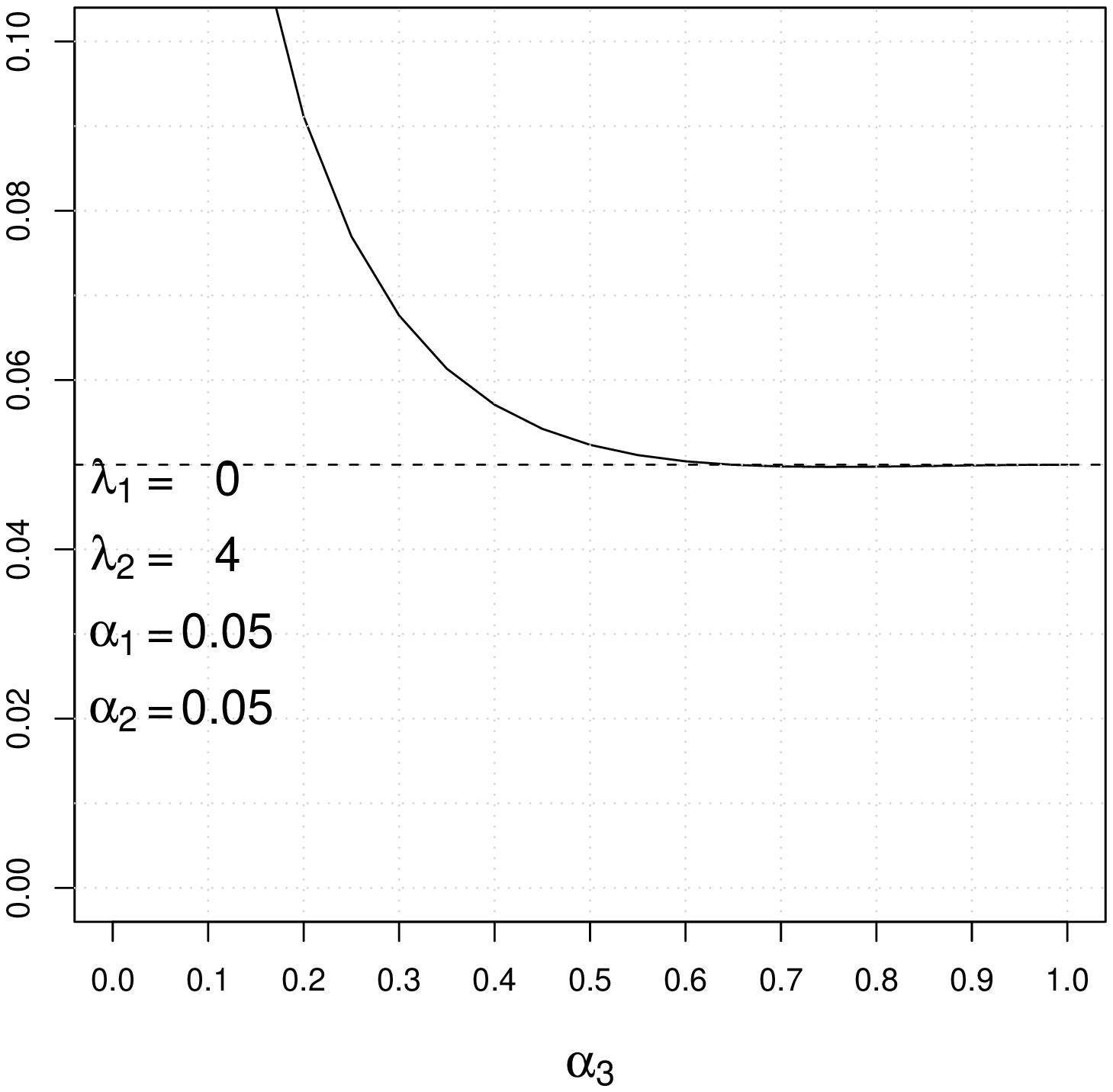} 
\includegraphics[width=0.3 \textwidth]{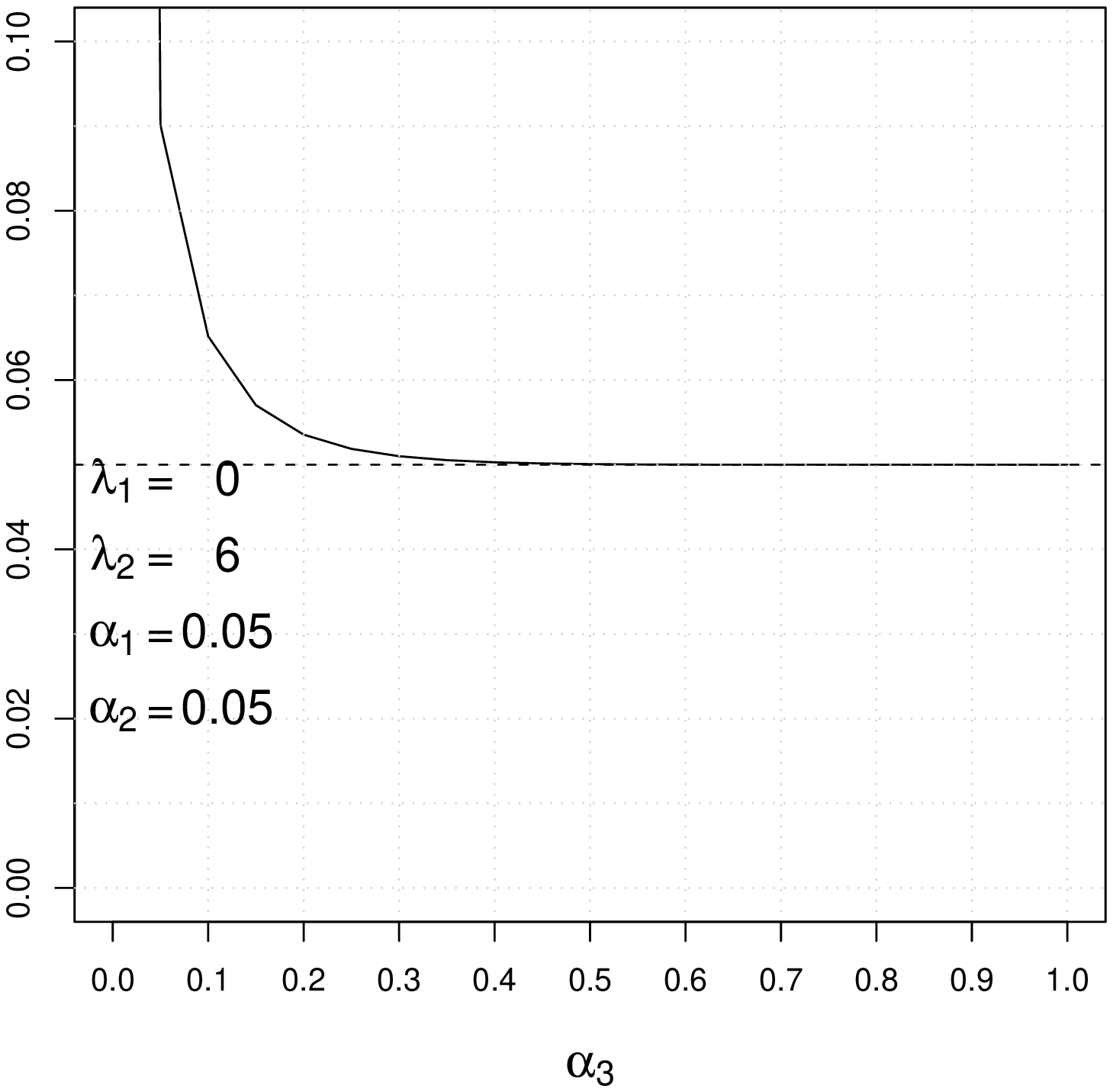} 
\includegraphics[width=0.3 \textwidth]{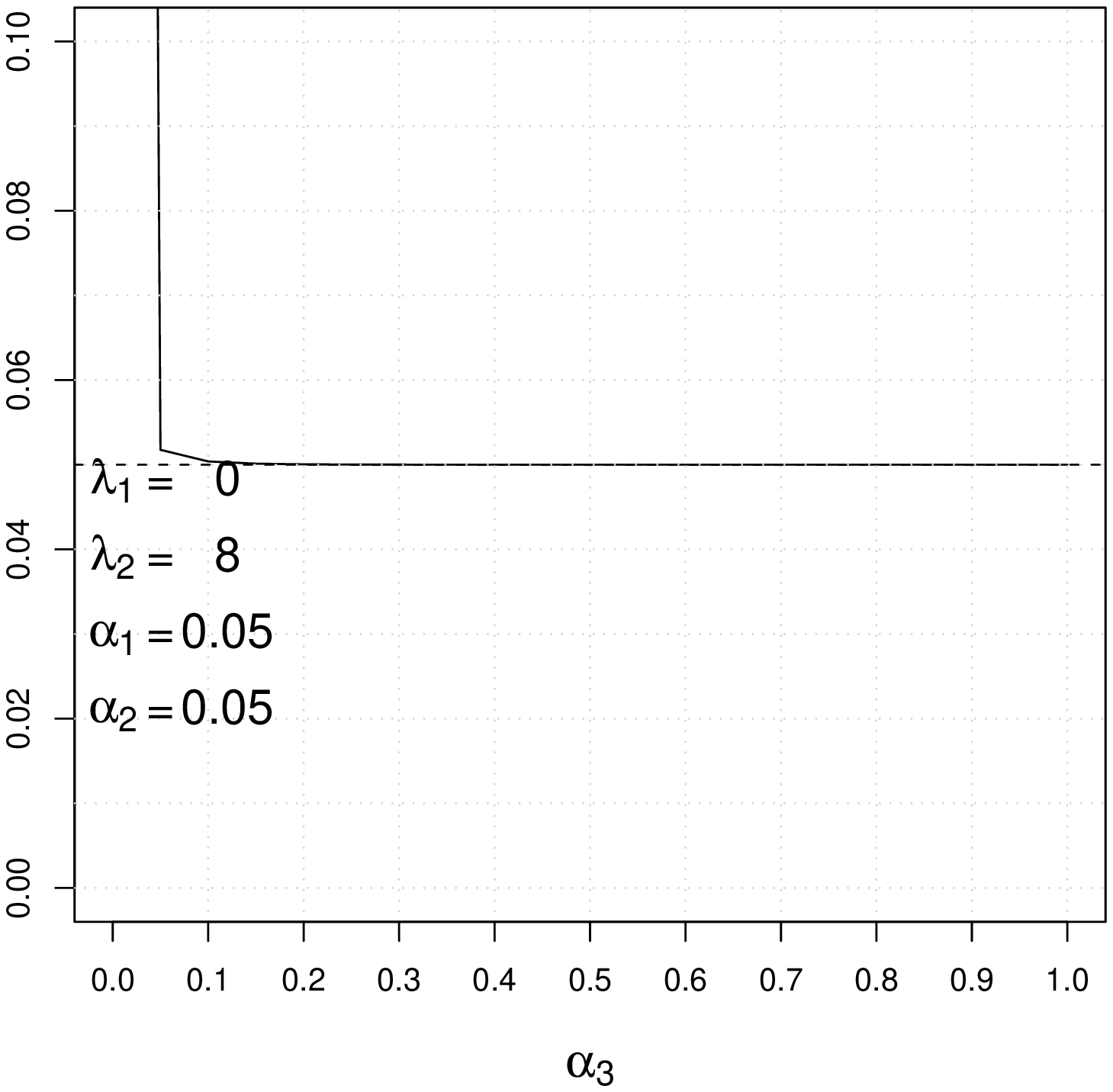}
\includegraphics[width=0.3 \textwidth]{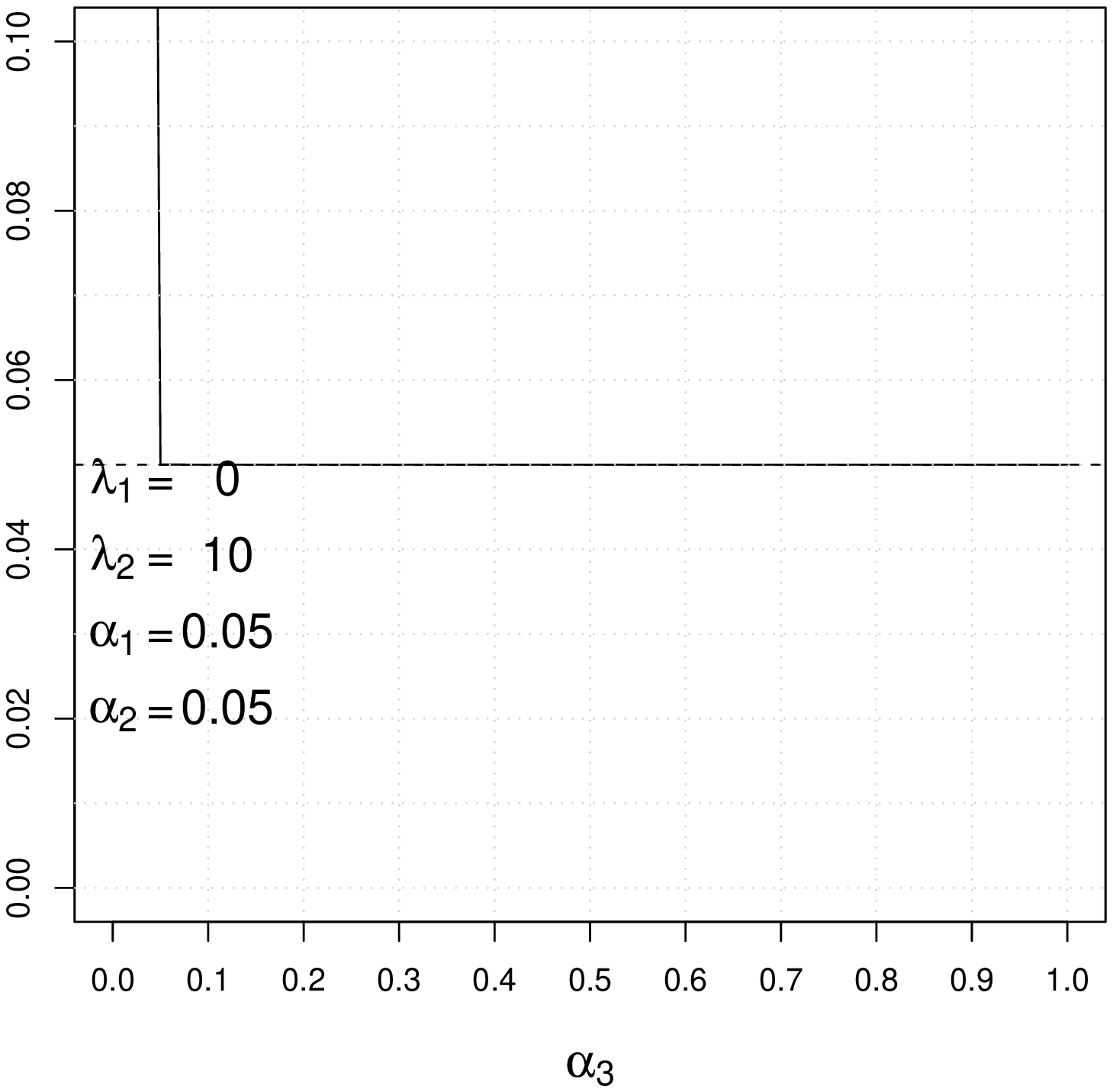} 
\caption{Graphs of size of the PTT
\((\alpha^{PTT}=\Pi^{PTT}(0,\lambda_2)) \) as \(\alpha_3 \)
and \(\lambda_2\) increasing when \(\bar{c} > 0 \) and \(\alpha_1=\alpha_2=0.05\) for all graphs.} 
\label{fig4}
\end{center}
\end{figure}


\begin{table}[!htpb]
\caption{Size of ultimate test \((\alpha^{PTT})\) as a function of
nominal size of PT \((\alpha_3)\) at selected values of
\(\lambda_2\) and \(\alpha_1=\alpha_2=\alpha=0.05\).}
\begin{center}
\begin{tabular}{|c|cc|cc|}\hline \(\lambda_2 \) &\(\alpha_3\)&
\(\alpha^{PTT}\)& \(\alpha_3\)&\(\alpha^{PTT}\)
\\ \hline \hline
0 & 0.05 & 0.0479 & 0.00 & 0.0500 \\ 0.1& 0.10 & 0.0495& 0.05&0.0525 \\
0.2 & 0.20 & 0.0476& 0.15 & 0.0509 \\
0.4 & 0.30 & 0.0475 & 0.25 &0.0515 \\
0.6 & 0.35 &0.0493 & 0.30 & 0.0540 \\
0.8 & 0.40 & 0.0498 & 0.35 & 0.0547\\
1.0 & 0.45 & 0.0491 & 0.40 &0.0543 \\
2.0 & 0.60 & 0.0476 & 0.55& 0.0508 \\
4.0 & 0.65 & 0.0500 & 0.60 & 0.5040 \\
6.0& 0.70 & 0.0500 & 0.65&0.0500\\
8.0 & 0.75 &0.0500 &  & \\ 10.0& 0.75 & 0.0500 &  &\\ \hline
\end{tabular}
\end{center}
\begin{center} \small{The \(\alpha^{PTT}\) is the actual achievable significance
level and \(\alpha_3\) is the nominal PT significance level.}
\end{center}
\end{table}

 Figure 5 shows graphs of
$\alpha^{PTT}=\Pi^{PTT}(0,\lambda_2)$ for $0 \leq \alpha_3 \leq 1$
at selected values of $\lambda_2$, $\alpha_1$ and $\alpha_2$ when
$\bar{c}>0.$ Equations (\ref{app1}) and (\ref{app2}) show that the
size and power of the PTT is approaching the size and power of the
RT as the nominal size of PT is closer to 0 but is approaching the
size and power of the UT as the nominal size of the PT is closer
to 1. From equation (\ref{app1}), setting the nominal significance
level $\alpha_3=0$ implies the size and power of the PTT is
entirely contributed by the size and power of the RT and none from
the UT. The contribution of the size and power of the UT to the
size and power of the PTT is not substantial when the nominal size
of the PT is small. From the graphs, the decreasing in the
contribution of the size of the RT reduces the size of the PTT as
$\alpha_3$ differs from zero. On the contrary, setting the nominal
size $\alpha_3=1$ causes the size of the PTT is totally
contributed by the size of the UT (see equation (\ref{app2})). The
contribution of the size of the RT is not significant when the
nominal size of the PT is large. As the value of $\alpha_3$
differs from 1, lesser contribution from the size of the UT
imposes smaller size of the PTT. The size of the PTT decreases
from both ends and the minimum of the size of the PTT is achieved
at a particular value of $\alpha_3.$

Further, analysis is carried out to investigate the dependence of
the size of the ultimate test to the changes in the nominal sizes
$\alpha_1$, $\alpha_2$ and $\alpha_3.$ From observation of Figures
5(a)-5(f), there is an increase in the percentage of
$\alpha^{PTT}$ in [0,0.1] for $\alpha_3$ in [0,0.2] when we set
smaller nominal size of $\alpha_2$ for a bit larger value of
$\lambda_2$. For example, there is 47.62\% of $\Pi^{PTT}(0,1)$ in
[0,0.10] for $\alpha_3$ in [0,0.2] at nominal size $\alpha_2=0.05$
(see Figure 5(a)) but there is 100\% of $\Pi^{PTT}(0,1)$ when
$\alpha_2=0.03$ (see Figure 5(b)). For some moderate values of
$\lambda_2$, there is an increment in the percentage of
$\Pi^{PTT}(0,\lambda_2)$ when we choose a smaller nominal size
$\alpha_2.$ But only small increment is observed for a larger
value of $\lambda_2$. For example, there is no $\Pi^{PTT}(0,3)$ in
[0,0.10] for $\alpha_3$ in [0,0.2] when we set the nominal size to
be $\alpha_2=0.05$ (see Figure 5(a)) but there is a slightly
4.76\% of $\Pi^{PTT}(0,3)$ when $\alpha_2=0.03$ (see Figure 5(b)).
The small increment suggests setting a much smaller value of
nominal size $\alpha_2$ maybe necessary to achieve a small size of
PTT with small nominal size of pre-test for moderate values of the
slope. However, this rule fails for a large value of $\lambda_2.$

\begin{table}[!htpb]
\caption{Size of ultimate test \((\alpha^{PTT})\) as a function of
nominal size of pre-test \((\alpha_3)\) at selected values of
\(\alpha_2\) and \(\lambda_2\) with \(\alpha_1=0.05\).}
\begin{center}
\begin{tabular}{|c|c|c|c|c|c|c|c|c|c|c|c|} \hline
\multicolumn{3}{|c|}{\(\lambda_2=0.5\)}
&\multicolumn{3}{|c|}{\(\lambda_2=1\)}&\multicolumn{3}{|c|}{\(\lambda_2=3\)}
&\multicolumn{3}{|c|}{\(\lambda_2=6\)} \\ \hline \hline

\multicolumn{1}{|c|}{\(\alpha_2\)}
&\multicolumn{1}{|c|}{\(\alpha_3\)}
&\multicolumn{1}{|c|}{\(\alpha^{PTT}\)}
&\multicolumn{1}{|c|}{\(\alpha_2\)}
&\multicolumn{1}{|c|}{\(\alpha_3\)}
&\multicolumn{1}{|c|}{\(\alpha^{PTT}\)}
&\multicolumn{1}{|c|}{\(\alpha_2\)}
&\multicolumn{1}{|c|}{\(\alpha_3\)}
&\multicolumn{1}{|c|}{\(\alpha^{PTT}\)}
&\multicolumn{1}{|c|}{\(\alpha_2\)}
&\multicolumn{1}{|c|}{\(\alpha_3\)}
&\multicolumn{1}{|c|}{\(\alpha^{PTT}\)}
 \\ \hline \hline
0.03 & 0.03 & 0.0498 & 0.01 & 0.00& 0.0355 & 0.01& 0.08 & 0.0994 &
0.03& 0.04& 0.0983 \\
\cline{2-3}\cline{5-6}\cline{8-9}\cline{11-12}

& 0.02& 0.0507& & 0.01&0.0343&&0.07& 0.1052&&0.03& 0.1145  \\
\hline 0.04& 0.19& 0.0499& 0.02& 0.00& 0.0623& 0.02& 0.15& 0.0987
&
0.04& 0.05& 0.0891\\
\cline{2-3}\cline{5-6}\cline{8-9}\cline{11-12} & 0.18& 0.0508&&
0.01& 0.0605& & 0.14& 0.1029& &0.04& 0.1002\\
\hline

0.05& 0.31& 0.0497& 0.03& 0.00& 0.0870& 0.03& 0.20& 0.0965& 0.05&
0.05& 0.0901\\
\cline{2-3}\cline{5-6}\cline{8-9}\cline{11-12}& 0.30& 0.0507& &
0.01& 0.0839& & 0.19&0.1005& & 0.04& 0.1014\\
\hline

0.06& 0.34& 0.0508& 0.04& 0.02& 0.1026& 0.04& 0.22& 0.1007& 0.06&
0.04& 0.1024\\
\cline{2-3}\cline{5-6}\cline{8-9}\cline{11-12}& 0.04&
0.0499&&0.03&0.0999& &0.23& 0.0972&& 0.05&0.0909\\
\hline

0.07& 0.48& 0.0491& 0.05& 0.11& 0.0993& 0.05&0.25&0.0987& 0.10&
0.05& 0.0093\\
\cline{2-3}\cline{5-6}\cline{8-9}\cline{11-12}& 0.47& 0.0500& &
0.10& 0.1015& & 0.24&0.1021&& 0.04& 0.1045\\
\hline
\end{tabular}
\end{center}
\begin{center} \small{The \(\alpha^{PTT}\) is the actual achievable significance
level and \(\alpha_3\) is the nominal pre-test significance
level.}
\end{center}
\end{table}

\begin{figure}
\begin{center}
\includegraphics[width=0.35 \textwidth]{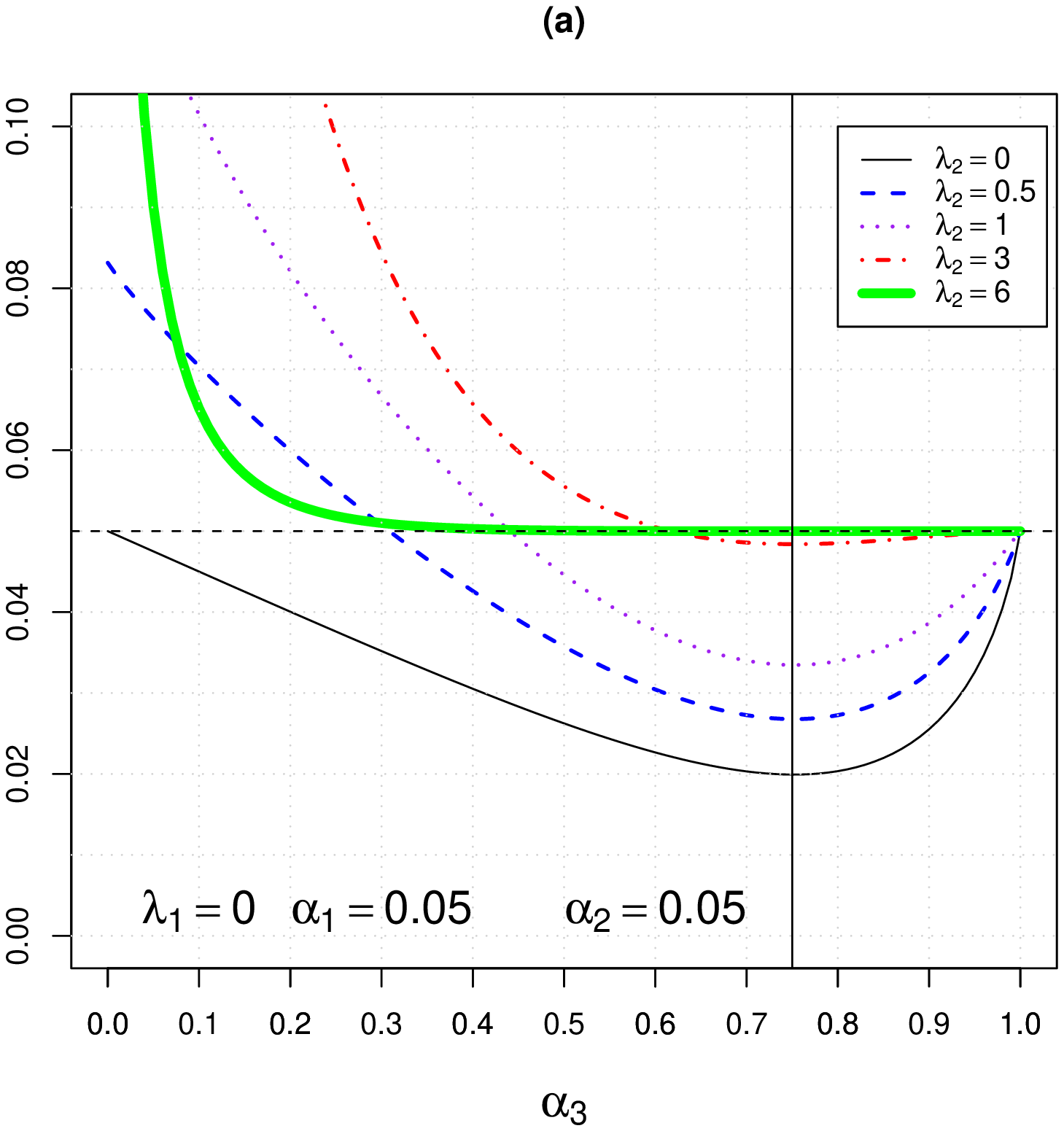} 
\includegraphics[width=0.35 \textwidth]{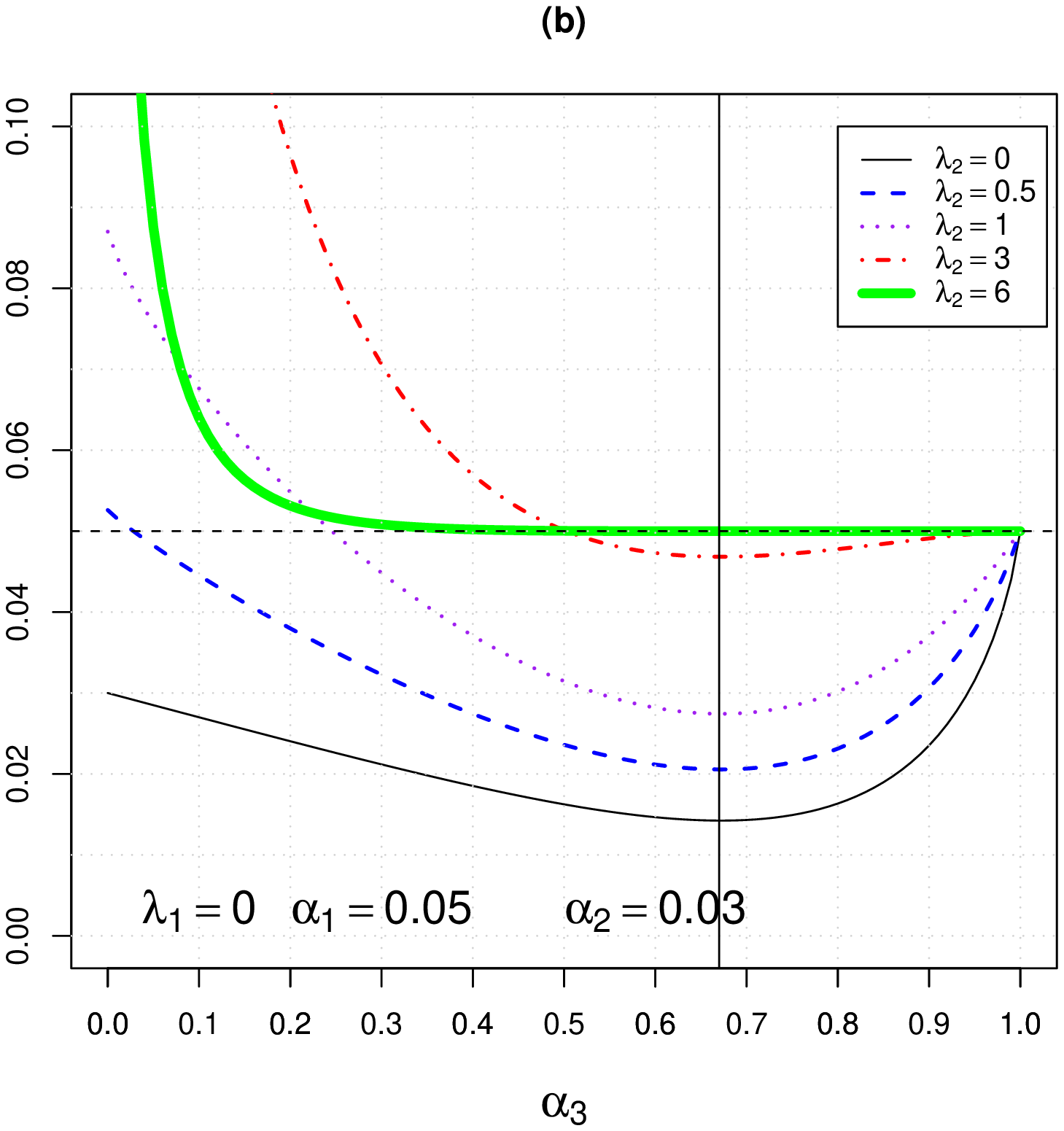} 
\includegraphics[width=0.35 \textwidth]{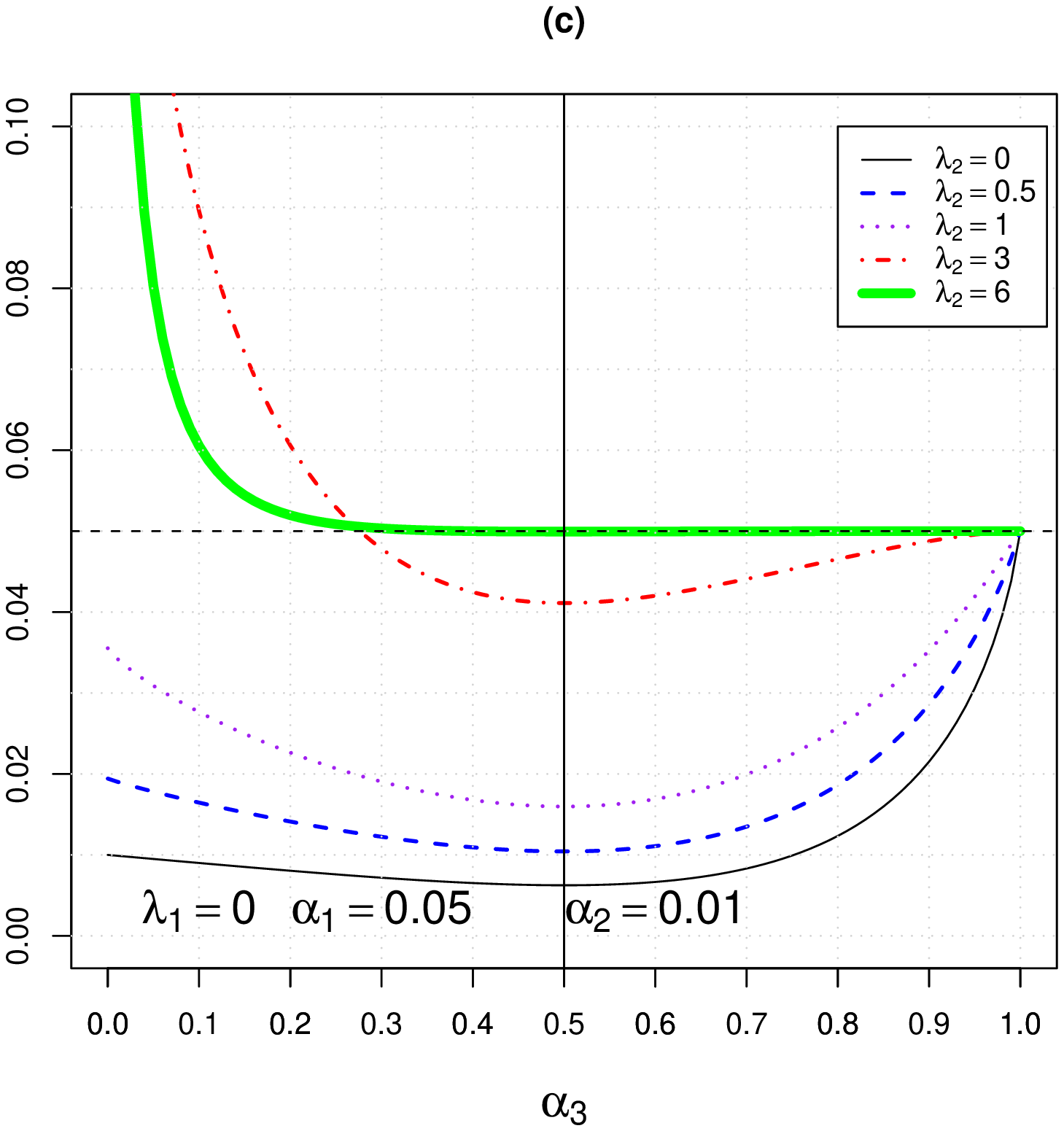} 
\includegraphics[width=0.35 \textwidth]{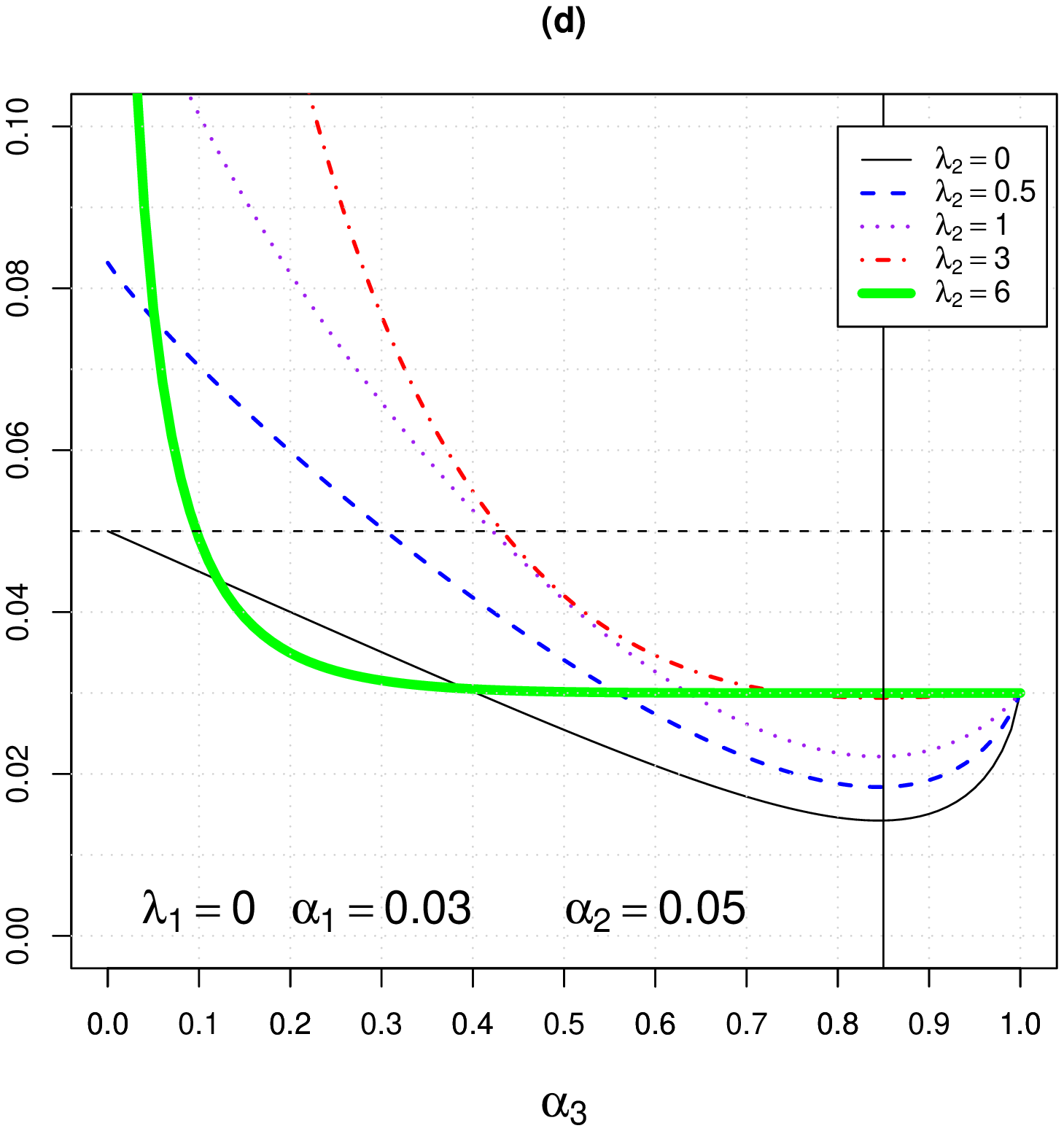} 
\includegraphics[width=0.35 \textwidth]{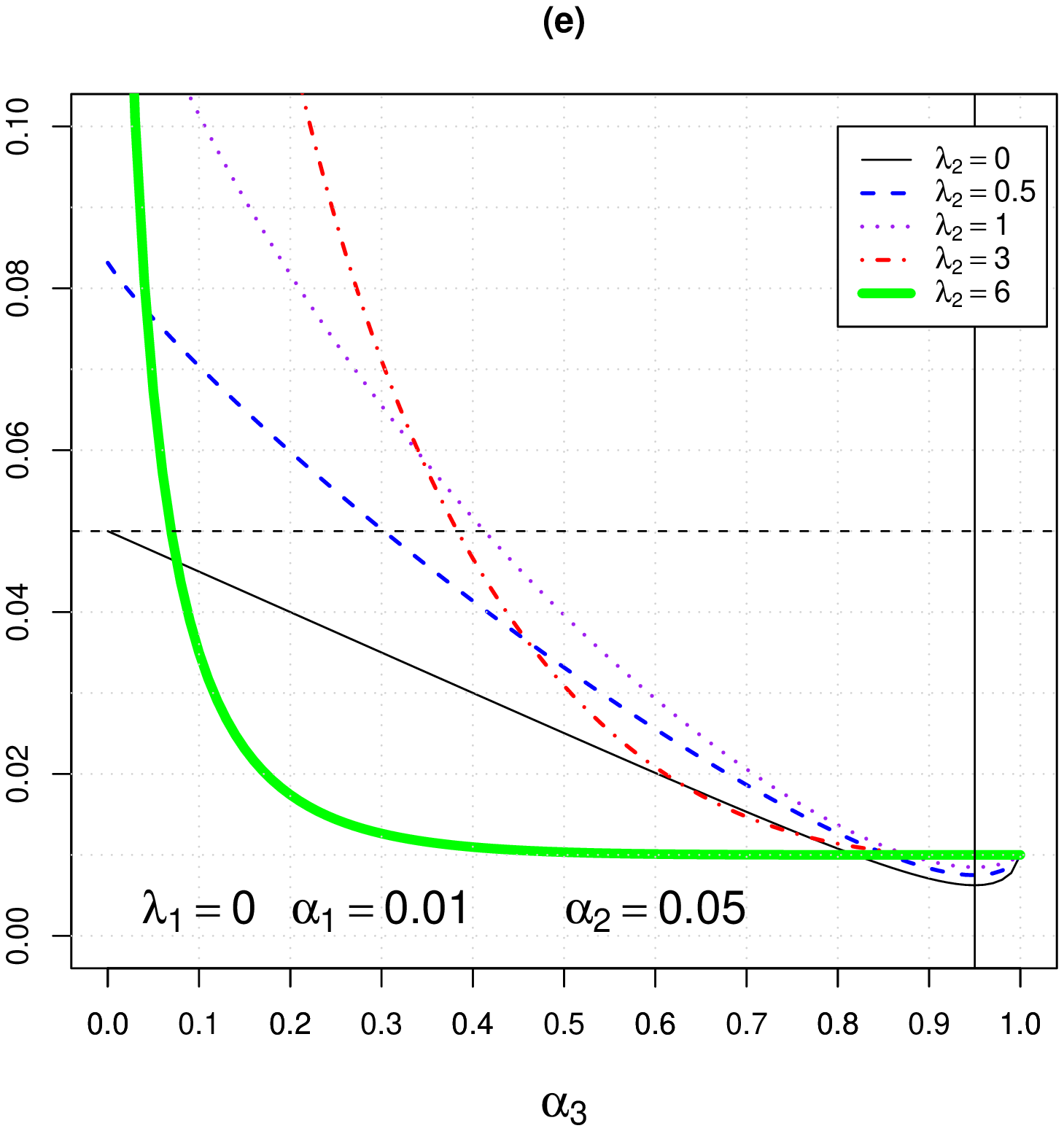} 
\includegraphics[width=0.35 \textwidth]{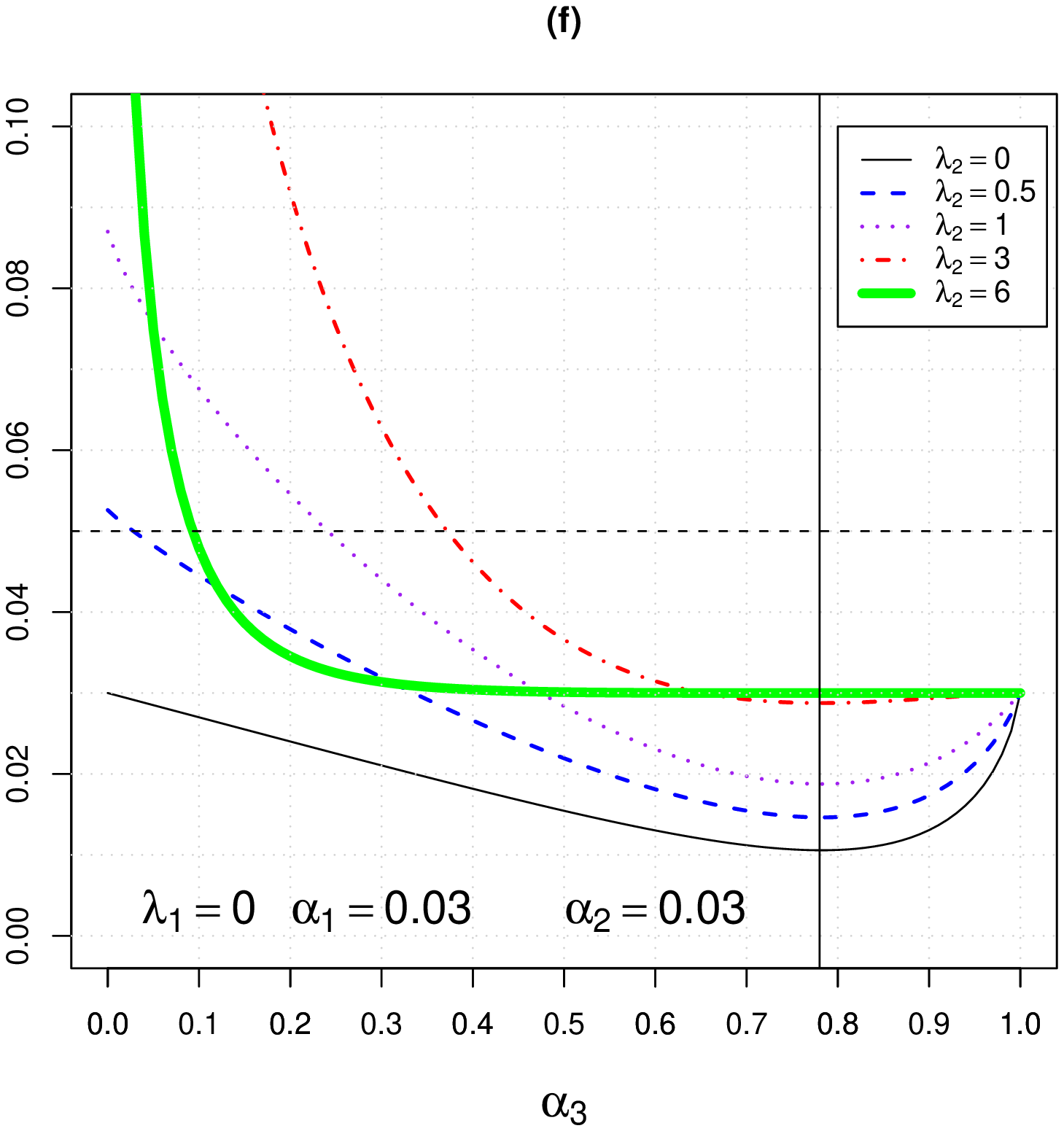} 
\caption{Graphs of the size of ultimate test for increasing
\(\alpha_3 \) selected at different values of nominal sizes of
\(\alpha_1\) and \(\alpha_2\) with $\bar{c}>0$. The intersection with the vertical line represents the minimum.} 
\label{fig6}
\end{center}
\end{figure}

We wish to have small size of the PTT by setting small nominal
sizes of $\alpha_1,\alpha_2$ and $\alpha_3.$ Figure 5 shows that
this could not be achieved when $\lambda_2$ is large and
$\alpha_3$ is very small (close to zero) even if we set a very
small value of $\alpha_2$. For instance, there is less than 100\%
(i.e. 80.95\%) of $\Pi^{PTT}(0,6)$ in [0,0.10] as $\alpha_3$ in
[0,0.2] (see Figures 5(a) and 5(b)) for both nominal sizes
$\alpha_2=0.03$ and $\alpha_2=0.05$. The percentage does not reach
100\% even $0<\alpha_2<0.03$ is chosen.

Since $\alpha^{PTT}$ behaves like $\alpha^{RT}$ when the nominal
size $\alpha_3$ is small, the null hypothesis $H_0^\star:\theta=0$
is rejected more often for small nominal size of $\alpha_3$ when
$\lambda_2$ is large because the nominal size $\alpha_2$ is
smaller than the actual size of the RT. The null hypothesis
$H_0^\star:\theta=0$ should not be rejected if the true value of
$\theta=0.$ In this case, however the possibility of rejection is
large when $\lambda_2$ differs much from 0 because $\beta=0$ is
assumed in the test statistic $T_n^{RT}.$ This fact answers the
reason why very small $\alpha_3$ (close to zero) has a very large
size of the test when $\lambda_2$ is large.

Table 2 shows the size of ultimate test as a function of nominal
size $\alpha_3$ for selected values of $\lambda_2$ and $\alpha_2$
with $\alpha_1=0.05$ and $\bar{c}>0$. The nominal sizes for the RT
and PT are given in the table for the size of ultimate test near
point 0.05 when $\lambda_2=0.5$ and 1 and near point 0.10 when
$\lambda_2=3$ and 6. The table enables us to observe the changes
in the values of the nominal size of PT ($\alpha_3$) as the
nominal size $\alpha_2$ changes and the significance level of the
PTT is around the same value. We wish to have small nominal size
of the PT that allow us to get 5 or 10\% of significance level of
the PTT. From the table, this is achieved by selecting smaller
nominal size of the RT for moderate and small values of slope.
When $\lambda_2=3$ (moderate value), selecting nominal size
$\alpha_2$ as small as 0.01 we have as much as 8\% of nominal size
of the PT to get below than 10\% significance level of ultimate
test (see Table 2, row:1, col:7-9). In column 1-3 of the table,
for $\lambda_2=0.5$ (small), approximately 5\% level of
significance of ultimate test is obtained by setting nominal size
of the RT = 0.05 and nominal size of the PT = 0.3 or by setting
both nominal sizes of the RT and PT = 0.03 but the latter with
smaller nominal sizes of the PT and RT is more preferable. For
larger value of the slope, as the nominal size of the PT closer to
0, the size of the PTT is growing too large. When $\lambda_2=6$
(large), to obtain at most 10\% of significance level of ultimate
test, the least nominal size for the PT that we should set is 5\%
(see Table 2, row:3, col:10-12) when the nominal size $\alpha_2$
is set from 0.05 to 0.10.

\section{Concluding Remarks}

The M-test of the UT, RT and PTT for testing the intercept are
provided in this paper. The asymptotic power functions of the
tests are derived by using the results from the asymptotic
sampling distribution of the statistics.

In the estimation regime, it is well known that the RE has the
smallest MSE if distance parameter (a function of $\beta-\beta_0$)
is 0 or close to 0, but its MSE is unbounded for larger values of
the distance parameter. The UE has a constant MSE that does not
depend on the distance parameter. The PTE has smaller MSE than
that of the RE for moderate and larger values of the distance
parameter. The PTE has smaller MSE than the UE if the value of
distance parameter is close or equal to 0. In the testing context,
the power functions of the UT, RT and PTT demonstrate a similar
behavior as the MSE of the UE, RE and PTE.

For a set of realistic values of the regressor, with mean value
larger than 0, the size of the RT is small when $\beta=0$ or close
to 0, but the size grows large and converges to 1 for larger
values of the slope. The UT has a constant size regardless of the
value of the slope (via $\lambda_2$). The PTT has smaller size
than that of the RT when the slope is 0 and very close to 0, and
significantly smaller than that of the the RT for moderate and
large values of the slope. The PTT has smaller size than the UT
for the value of slope is 0 or very close to 0.

Again for a set of realistic values of the regressor, with mean
larger than 0, the RT is the best choice for having largest power
but the worst choice for having largest size. The size of the UT
is constant regardless of the value of the slope. The UT is the
best choice for having smallest size but the worst choice for
having smallest power. The PTT has smaller size than the RT for
moderate and larger values of the slope and has larger power than
the UT for smaller and moderate values of the slope. Therefore,
the power function of the PTT is found to behave similar to the
MSE of the PTE in the sense that though it is not uniformly the
best statistical test with the smallest size and the largest power
but it protects from the risk of a too large size and a too small
power. Thus, the power function of the PTT is a compromise between
that of the UT and RT. In the face of uncertainty on the value of
the slope, if the objective of a researcher is to minimize the
size and maximize the power of the test, the PTT is the best
choice.

The tables and graphs support the analytical asymptotic comparison
of the UT, RT and PTT as discussed in Section 6. The analysis is
furthered by investigating the relationship between the power
functions and its arguments, namely the slope and the nominal
sizes, of the UT, RT and PT. The chosen values of the nominal
sizes that are set before testing affect the actual size of the
PTT.

In order to get small probability of type I error for the ultimate
test, our investigations concentrate on small nominal sizes of the
UT, RT and PT with a view to achieving  small (actual)
significance level of the PTT. The study revealed that for small
and moderate values of slope, the smaller the nominal size of the
RT, the smaller the size of the PTT when other nominal sizes are
kept fixed and small. For moderate and large values of the slope,
a large size of the PTT is observed when nominal size of PT is set
close to 0. The size of the PTT behaves much like that of the RT
when the nominal size of PT is small, but it behaves more like
that of the UT when the nominal size of the PT is large.


The power of the ultimate test is larger for moderate values of
the slope than for smaller and larger values of the slope. It is
shown analytically that the power of the PTT approaches the power
of the RT when the nominal size of PT is closer to 0 but
approaches the power of the UT when the nominal size of the PT is
closer to 1. In practical applications, size of the PT should be
small (ideally close to 0), and in such cases the power of the PTT
is close to that of the RT (which is much higher than that of the
UT). To avoid the larger size of the RT, practitioners are
recommended to use the PTT as it achieves smaller size (than the
RT) and higher power (than the UT) when the value of the slope is
small or moderate. Even for large values of the slope the PTT has
at least as much power as the UT.\\


\noindent{\bf Acknowledgements}\\
The authors thankfully acknowledge valuable suggestions of
Professor A K Md E Saleh, Carleton University, Canada that helped
improve the content and quality of the results in the paper.\\


\noindent {\bf References}\\
\baselineskip=12pt
\def\ref{\noindent\hangindent 15pt}

\ref Bechhofer, R.E. (1951). The effect of preliminary test of
significance on the size and power of certain tests of univariate
linear hypotheses. Ph.D. Thesis (unpublished), Columbia Univ.

\ref  Bozivich, H., Bancroft, T.A. and Hartley, H. (1956). Power
of analysis of variance test procedures for certain incompletely
specified models. {\em Ann. Math. Statist.} \textbf{27}, 1017 -
1043.

\ref  Caroll, R.J. and Rupert, D. (1988). {\em Transformation and
Weighting in Regression.} Chapman \& Hall, US.

\ref H$\acute{a}$jek, J., $\check{S}$id$\acute{a}$k, Z., and Sen,
P.K. (1999). {\em Theory of Rank Tests.} Academia Press, New York.

\ref  Hoaglin, D.C., Mosteller, F. and Tukey, J.W. (1983). {\em
Understanding Robust and Explanatory Data Analysis}. John Wiley
and Sons, US.

\ref Huber, P.J. (1981). {\em Robust Statistics.} Wiley, New York.

\ref  Jur$\check{e}$ckov$\acute{a}$, J. (1977). Asymptotic
relations of M-estimates and R-estimates in linear regression
model. {\em Ann. Statist.} \textbf{5}, 464-72.

\ref  Jur$\check{e}$ckov$\acute{a}$, J. and Sen, P.K. (1981).
Sequential procedures based on M-estimators with discontinuous
score functions. {\em J. Statist. Plan. Infer.} \textbf{5},
253-66.

\ref  Jur$\check{e}$ckov$\acute{a}$, J. and Sen, P.K. (1996). {\em
Robust Statistical Procedures Asymptotics and Interrelations.}
John Wiley \& Sons, US.

\ref  Khan, S., and Saleh, A.K.Md.E. (2001). On the comparison of
pre-test and shrinkage estimators for the univariate normal mean.
{\em Stat. papers.} \textbf{42}, 451-473.

\ref  Khan, S., Hoque, Z., and Saleh, A.K.Md.E. (2002). Estimation
of the slope parameter for linear regression model with uncertain
prior information. {\em J. Stat. Res.} \textbf{36}, 55-74.

\ref  Saleh, A.K.Md.E. and Sen, P.K. (1982). Non-parametric tests
for location after a preliminary test on regression. {\em
Communications in Statistics: Theory and Methods.} \textbf{11},
639-651.

\ref  Saleh, A.K.Md.E. (2006). {\em Theory of Preliminary test and
Stein-type estimation with applications}. John Wiley \& Sons, New
Jersey.

\ref  Schrader, R. M. and Hettmansperger, T.P. (1980). Robust
analysis of variance based upon a likelihood ratio criterion. {\em
Biometrika.} \textbf{67}, 93-101. \ref{sen82} Sen. P.K. (1982). On
M-tests in linear models. {\em Biometrika.} \textbf{69}, 245-248.

\ref  Sen, P.K. and Saleh, A.K.Md.E. (1987). On preliminary test
and schrinkage M-estimation in linear models. {\em Ann. Statist.}
\textbf{15}, 1580-1592.

\ref  van der Vaart, A.W. (1998). {\em Asymptotic statistics}.
Cambridge University Press, UK. \ref{wilcox05} Wilcox, R.R.
(2005). {\em Introduction to Robust Estimation and Hypothesis
Testing.} Elsevier Inc, US.

\end{document}